\documentclass[3p]{elsarticle}

\usepackage{bm}     
\usepackage{amsmath}
\usepackage[utf8x]{inputenc} 	
\usepackage{lmodern}     
\usepackage[T1]{fontenc} 
\usepackage{graphicx}
\usepackage{latexsym}
\usepackage{amssymb}
\usepackage{amsbsy}
\usepackage{color}
\usepackage{lineno}
\usepackage{ulem}
\usepackage{multirow} 
\usepackage{rotating}
\usepackage{tikz}
\usepackage{url}
\modulolinenumbers[5]
\journal{Journal of Computational Physics}

\def\be{\begin{equation}}
\def\ee{\end{equation}}
\DeclareSymbolFont{msbm}{U}{msb}{m}{n}
\DeclareMathSymbol{\C}{\mathalpha}{msbm}{'103}
\DeclareMathSymbol{\R}{\mathalpha}{msbm}{'122}
\DeclareMathSymbol{\Z}{\mathalpha}{msbm}{'132}
\DeclareMathSymbol{\N}{\mathalpha}{msbm}{'116}
\def\bea{\begin{eqnarray}}
\def\ba{\begin{array}{l}\displaystyle}
\def\eea{\end{eqnarray}}
\def\ea{\end{array}}
\def\E{{\cal E}}

\def\RR{\mathbb R}
\newcommand{\Int}{\displaystyle\int}

\newfont{\numerikEleven}{ecrm1000}
\newfont{\numerikTen}{cmss10}
\newfont{\numerikNine}{cmss9}
\newfont{\numerikEight}{cmss8}
\newcommand{\SSS}{\mathbb{S}}
\DeclareMathOperator{\sinc}{Sinc}
\def\supp{{\rm Supp}}
\def\supp{\mbox{supp }}
\newcommand{\ens}[1]{\mathbb{#1}}
\def\Ball{ {\cal B}}
\def\QL{Q^{R}}

\def\f{\hat f}

\def\bb{\hat \beta}
\def\B{\hat B}

\begin{document}

\begin{frontmatter}

\title{ An efficient numerical method for solving the Boltzmann equation in multidimensions}

\author{Giacomo Dimarco
}
\address{Department of Mathematics and Computer Science, University of Ferrara,
Via Machiavelli 30, 44122 Ferrara, Italy.}

\cortext[mycorrespondingauthor]{Corresponding author}
\ead{giacomo.dimarco@unife.it}


\author{Rapha\"{e}l Loub\`{e}re
}
\address{CNRS and Institut de Math\'{e}matiques de Toulouse (IMT)
 Universit{\'e} Paul-Sabatier, Toulouse, France.
}
\author{Jacek Narski
}
\address{CNRS and Institut de Math\'{e}matiques de Toulouse (IMT)
 Universit{\'e} Paul-Sabatier, Toulouse, France.
}
\author{Thomas Rey
}
\address{Laboratoire Paul Painlev{\'e}, Universit{\'e} Lille 1,
 Lille, France.
}

\begin{abstract}
In this paper we deal with the extension of the Fast Kinetic Scheme (FKS) [J. Comput.
Phys., Vol. 255, 2013, pp 680-698] originally constructed for solving
the BGK equation, to the more challenging case of the Boltzmann equation. The scheme combines a robust and fast method for treating
the transport part based on an innovative Lagrangian technique supplemented with fast spectral schemes to treat the collisional operator by means of an operator splitting approach.
This approach along with several implementation features related to the parallelization of the algorithm permits to construct an efficient simulation tool which is numerically tested
against exact and reference solutions on classical problems arising in rarefied gas dynamic. We present results up to the $3$D$\times 3$D case for unsteady flows for the Variable Hard Sphere 
model which  may serve as benchmark for future comparisons between different numerical methods for
 solving the multidimensional Boltzmann equation. 
For this reason, we also provide for each problem studied details on the computational cost and memory consumption as well as comparisons with the BGK model or the limit model of compressible Euler equations. 
\end{abstract}

\begin{keyword}
\texttt{Boltzmann equation, kinetic equations, semi-Lagrangian schemes, spectral schemes, 3D/3D, GPU, CUDA, OpenMP, MPI.}
\MSC{(82B40, 76P05, 65M70, 65M08, 65Y05, 65Y20}
\end{keyword}

\end{frontmatter}


\tableofcontents


\section{Introduction} \label{sec:introduction}
Kinetic equations consider a representation of a gas as particles undergoing instantaneous collisions
interspersed with ballistic motion \cite{cercignani, ACTA}. Nowadays, these models appear in a variety of sciences and applications such as astrophysics, aerospace and nuclear engineering, semiconductors, plasmas related to fusion processes as well as biology, medicine or social sciences. The common structure of such equations consists in a combination of a linear transport term with one or more interaction terms which furnishes the time evolution of the distribution of particles in the phase space. For its nature, the unknown distribution lives in a seven dimensional space, three dimensions for the physical space and three dimensions for the velocity space, plus the time. 
This makes the problem a real challenge from the numerical point of view, since the computational cost becomes immediately prohibitive for realistic multidimensional problems \cite{ACTA}. Aside from the curse of dimensionality problem, there are many other difficulties which are specific to kinetic equations. We recall two among
the most important ones. The computational cost related to the evaluation of the collision operator involving multidimensional integrals which should be solved in each point of the physical space \cite{FiMoPa:2006, PR2}. The second challenge is represented by the presence of multiple scales which requires the development of adapted numerical schemes to avoid the resolution of the stiff dynamics \cite{Dimarco_stiff1, Dimarco_stiff2, Jin_review, Jin, Jin2, BLM, degondrev} typically arising when dealing with problems with multiple regimes. 

Historically, there exists two different approaches which are generally used to tackle kinetic equations from a numerical point of view: deterministic numerical methods such as
finite volume, semi-Lagrangian and spectral schemes \cite{ACTA}, and, probabilistic numerical methods such as Direct Simulation Monte Carlo (DSMC) schemes \cite{bird, Cf}.
Both methodologies have strengths and weaknesses. While the first could normally reach high order of accuracy, the second are often faster, especially
for solving steady problems but, typically, exhibit lower convergence rate and difficulties in describing non stationary and slow motion flows.

In this work, we deal with deterministic techniques. In particular, we focus on semi-Lagrangian approaches \cite{CrSon,CrSon1, Gu, Shoucri, FilbetRusso, Filbet2} for the transport part coupled with spectral methods \cite{canuto:88} for the interaction part. Our main goal is to tackle the challenges related to the high dimensionality of the equations and with the difficulties related to the approximation of the collision integral. More in details, we deal with the extension of our recent works \cite{FKS, FKS_HO, FKS_DD, FKS_GPU} which were based on an innovative semi-Lagrangian technique for discretizing the transport part of a kinetic model (FKS method). This technique has been applied  solely to the solution of a simple kinetic equation with a relaxation type collision operator, \textit{i.e.} the BGK (Bhatnagar-Gross-Krook) operator \cite{Gross}. Here we extend it to the case of the full Boltzmann operator \cite{bird, cercignani} and we test his performances up to the six dimensional case. The FKS is based on the classical discrete velocity models (DVM) approach \cite{bobylev, Pal, Pal1, Mieussens}. Successively, in order to overcome the problem of the excessive computational cost, we propose to use a Lagrangian technique which exactly solves the transport step on the entire domain, without reprojecting the solution on the grid at each time step. The FKS approach was shown to be an efficient way to solve linear transport equations, and, it has permitted the simulation of full six dimensions problems on a single processor machine \cite{FKS}. Unfortunately the solutions computed with this method are limited to a first order in space and time precision. Extension of the method to high order reconstruction is under consideration \cite{RFKS}. Concerning the discretization of the collision operator we rely on Fourier techniques. For the resolution of the Boltzmann integral, these techniques have been first introduced independently by L. Pareschi and B. Perthame in \cite{pareschi:1996} and by A. Bobylev and S. Rjasanow in \cite{bobylev_rjasanow97}. Since then, this approach has been investigated by a many authors \cite{FiMoPa:2006, PR2, gamba:2010, gamba, gambaL, wu2013deterministic, bobylev_rjasanow97, Rjasanow, Villani, filbet:2011Conv, wu2015influence, Aleks, PaRu:stab:00, wu2015fast}. In this work we will make use in particular of the fast method described in \cite{MoPa:2006,FiMoPa:2006} which has a complexity of the order of $\mathcal O(N^{d_v} \log(N^{d_v}))$ where $N$ is the number of point in which the velocity space is discretized in one direction and $d_v$ the dimension of the velocity space. The method preserves mass, and approximates with spectral accuracy momentum and energy.

Combining opportunely the FKS method with the fast spectral approach we have developed a method for solving the Boltzmann equation up to the six dimensional case for unsteady flows. In order to obtain such result we  have constructed a parallel version of our algorithm taking advantages of Graphical-Processor-Unit (GPU) under CUDA language. The results presented in this work show that we are nowadays ready and able to use kinetic equations to simulate realistic multidimensional flows. Up to our knowledge this is one of the first examples in literature of solution of the full multidimensional Boltzmann equation by means of deterministic schemes. Our main limitation to run extensive simulations remains at the present moment the lack of memory capacity due to the use of shared memory machines. However, in the present paper, we also report some preliminary results about the extension of our method to deal with distributed-memory computers but we postpone for a future research the detailed analysis of the Message Passive Interface (MPI) version of the algorithm as well as the analysis of its performances.

The article is organized as follows. In Section \ref{sec:FKS} we recall the Boltzmann equation. 
In Section \ref{ssec:FKScheme} we present the Fast Kinetic Scheme and the fast spectral scheme. In Section \ref{sec:implementation} we detail the aspects related to the implementation of the resulting
algorithm necessary to realize an efficient parallel numerical tool. Several tests starting from the 0D$\times$ 2D up to the 3D$\times$3D case are studied in details in Section \ref{sec:numerics}. These tests assess the validity of our approach as well as detail all the computational aspects. 
Conclusions and future works are finally exposed in Section \ref{sec:conclusion}.


%
%
\section{The Boltzmann equation} \label{sec:FKS}
In this section we briefly recall the Boltzmann equation and its main properties, we also recall some related models, i.e. the BGK model and the compressible Euler equations, which will be used for numerical comparisons. We refer the reader to \cite{ACTA, cercignani} for an exhaustive description. 

In the kinetic theory of rarefied gases, the non-negative function $f(x,v,t)$ characterizes the state of the system and it defines 
the density of particles having velocity $v\in\R^{d_v}$ in position $x\in\R^{d_x}$ at time $t\in\R^+$, where $d_x$ is the physical dimension
and $d_v$ the dimension of the velocity space. The time evolution of the particle system is obtained through the equation
\be
\frac{\partial f}{\partial t} + v\cdot \nabla_x f= Q(f).
\label{eq:i1}
\ee
The operator $Q(f)$, on the right hand side in equation (\ref{eq:i1}), 
describes the effects of particle interactions and its form depends on the details of the microscopic dynamic. 
Independently on the type of microscopic interactions considered, typically the operator is characterized by some conservation properties of the physical system.
They are written as
\be
\int_{\R^{d_v}} Q(f)\phi(v)\,dv=0,
\label{eq:i5}
\ee
where $\phi(v)=(1,v,|v|^2)$ are commonly
called the collision invariants.
We denote by $$U(x,t)=\Int_{\R^{d_v}} f(x,v,t) \phi(v)\,dv\in\R^{2+d_v}$$ the first three moments of the distribution function $f$, namely,
$U(x,t)=(\rho,\rho u,E)$, where $\rho$ is the density, $\rho u$ the momentum and $E$ the energy.
Integrating (\ref{eq:i1}) against $\phi(v)$ yields a system of macroscopic conservation laws
\be
\frac{\partial}{\partial t}\int_{\R^{d_v}}f\phi(v)\,dv + \int_{\R^{d_v}} v \cdot \nabla_x f\,\phi(v)\,dv =0.
\label{eq:i4}
\ee
The above moment system is not closed since the second term involves higher order moments of the distribution function $f$. 
However, using the additional property of the operator $Q(f)$ that the functions belonging to its kernel satisfy \be Q(f)=0\quad\hbox{iff}\quad f=M[f],\label{eq:i6}\ee 
where the Maxwellian distributions $M[f]=M[f](x,v,t)$ can be expressed in terms of the set of moments $U(x,t)$ by
\be M[f]=\frac{\rho}{(2\pi T)^{d_v/2}}e^{\frac{-(v-u)^2}{2T}}, \qquad \frac{1}{2}d_v\rho T=E-\frac{1}{2}\rho|u|^2  \ee
with $T$ the temperature, one can get a closed system
by replacing $f$ with $M[f]$ in (\ref{eq:i4}).
This corresponds to the set of compressible Euler equations
which can be written as
\be
\frac{\partial U}{\partial t} + \nabla_x \cdot F(U) =0,
\label{eq:i7}
\ee
with $F(U)=\Int_{\R^{d_v}} M[f]v\phi(v)\,dv$.
The simplest operator satisfying (\ref{eq:i5}) and (\ref{eq:i6}) is the relaxation operator~\cite{Gross}
\be
Q_{BGK}(f)=\nu(M[f]-f),
\label{eq:ibgk}
\ee 
where $\nu=\nu(x,t) > 0$ defines the so-called collision frequency. Its values will be specified in the numerical test Section in order for the model to be as close as possible to the Boltzmann model
described next. The classical Boltzmann operator reads
\be
Q_B(f)=\int_{\R^{d_v}}
\int_{\SSS^{d_v-1}} B(|v-v_*|,\omega) \left( f(v')f(v'_*)-f(v)f(v_*) \right) dv_* d \omega,
\label{eq:bolt}
\ee
where $\omega$ is a vector of the unitary sphere
$\SSS^{d_v-1} \subset \R^{d_v}$. The post-collisional velocities $(v',v'_*)$ are given by the
relations \be v'=\frac{1}{2}(v+v_*+|q|\omega),\quad
v'_*=\frac{1}{2}(v+v_*+|q|\omega), \ee
where $q=v-v_*$ is the
relative velocity. The kernel $B$ characterizes the
details of the binary interactions, it has the form
\begin{equation}
\label{defVHSKernel}
B(|v-v_*|,\cos\theta)=|v-v_*|\sigma(|v-v_*|,\cos\theta)
\end{equation}
where the {scattering cross-section} $\sigma$, in the case of inverse $k$-th power forces between particles, can be written as 
\be \sigma(\vert v - v_{\ast} \vert,
\cos\theta) = b_{\alpha}(\cos\theta) \vert v - v_{\ast}
\vert^{\alpha-1}, 
\ee with $\alpha=(k-5)/(k-1)$. 
The special situation {$k=5$} gives the so-called {Maxwell pseudo-molecules model} with  \be B(v,v_*,\omega)=
b_{0}(\cos\theta). \ee
For the Maxwell case the collision kernel is independent of the relative velocity. For numerical purposes, a
widely used model is the {variable hard sphere} (VHS) 
model introduced by Bird~\cite{bird}. The model corresponds to {$b_{\alpha}(\cos\theta)=C_\alpha$,}
where {$C_\alpha$} is a positive constant, and hence
\be \sigma(\vert v - v_{\ast} \vert, \cos\theta) =
C_{\alpha} \vert v - v_{\ast} \vert^{\alpha-1}. \ee
In the numerical test Section we will consider the Maxwell molecules case when dealing with a velocity space of dimension $d_v=2$ and with the more realistic case of VHS molecules when
dealing with the three dimensional case in velocity space: $d_v=3$. For comparison purposes we will also consider the simplified BGK model (\ref{eq:ibgk}) and the compressible Euler system 
case (\ref{eq:i7}).


%
%

\section{The Fast Kinetic Scheme and the Fast Spectral Scheme}
\label{ssec:FKScheme}
In this section we detail the numerical scheme. 
In the first part we discuss the Fast Kinetic Scheme (FKS) and in the second part the Fast Spectral method
for the Boltzmann equation. The two solvers are connected by using splitting in time approaches as the ones described in \cite{Des}.

\subsection{The Fast Kinetic Scheme (FKS)}
\label{ssec:FKScheme1}

The Fast Kinetic Scheme \cite{FKS, FKS_HO} belongs to the family 
of so-called semi-Lagrangian schemes \cite{CrSon, CrSon1, Filbet}
which are typically applied to a Discrete Velocity Model (DVM) \cite{bobylev, Mieussens} approximation of the original kinetic equation.

In order to introduce the scheme, let us truncate the velocity space by fixing some given bounds and set a cubic grid in velocity space of $N$ points with $\Delta v$ the
grid step which is taken equal in each direction. The continuous distribution function $f$ is then replaced by a vector whose components are assumed to be approximations of the 
distribution function $f$ at locations $v_k$:
\be
\tilde f_{k}(x,t) \approx f(x,v_{k},t).
\ee
The discrete velocity kinetic model consists then of a set of $N$ evolution
equations in the velocity space for $\tilde f_k$, $1\leq k \leq N$, of the form
\be
\partial_t \tilde f_{k} + v_{k} \cdot\nabla_{x} \tilde f_{k} = Q(\tilde f_k), 
\label{eq:DM1}
\ee
where $Q(\tilde f_k)$ is a suitable approximation of the collision operator $Q(f)$ at location $k$ discussed in Section (\ref{sec:Boltz}). 
Observe that, due to the truncation of the velocity space and to the finite number of points with which $f$ is discretized, the moments of the discrete distribution function $\tilde f_k$ are such that
\be
\widetilde{U}(x,t) = \sum_{k}  \phi_{k} \, \tilde f_k(x,t)\, \Delta v \ne U(x,t),
\label{eq:DM}
\ee
with $\phi_k=(1,v_k,|v_k|^2)$ the discrete collision invariants, i.e. they are no longer those given by the continuous distribution $f$. This problem concerns all 
numerical methods based on the discrete velocity models and different strategies can be adopted to restore the correct macroscopic physical quantities 
\cite{Mieussens,Pal,Pal1,gamba}. At the same time, the spectral approach does not preserve the energy and the momentum of the system even if it approximates these quantities
with spectral accuracy \cite{PR2}. In order to solve this problem of loss of conservation we adopt an $L_2$ projection technique for the discretized distribution and the discretized collision 
operator which permits to project the discretized $\tilde f_k$ and $ Q(\tilde f_k)$ 
in the space of the distributions for which the moments are exactly the continuous macroscopic quantities 
we aim to preserve. We detail this approach in the following paragraph. For lightening the notations we omit the tilde in the rest of the paper, we only use it in the next paragraph in order to introduce the
$L_2$ projection procedure. We will then suppose from now on that all the operations preserves the macroscopic moments if not differently stated.

\begin{figure}
  \hspace{-2cm}
  \begin{center}
    \includegraphics[width=1.0\textwidth]{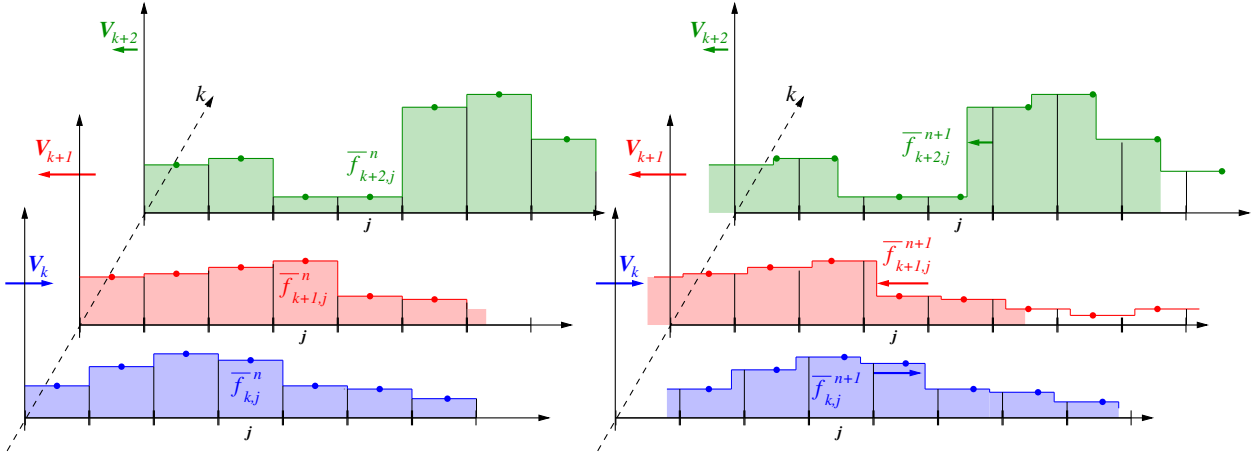}
    \caption{ \label{fig:transport}
      Illustration of the transport scheme for the FKS scheme.
      Left panel before transport step, right panel after transport step. Each discrete velocity (index $k$) drives its own
      transport equation with velocity $v_k$. The representation
      of $f$ is made by means of a piecewise constant function.
      The shape of the entire function has not changed during the transport
      but the cell-centered values (bullets) may have.}
  \end{center}
\end{figure}
Let us now introduce a Cartesian uniform grid in 
the three dimensional physical space made of $M$ points with $\Delta x$ a scalar which represents the grid step (the same in each direction) in the
physical space. Further we define a time discretization $t^{n+1}= t^n + \Delta t$ starting at $t^0$,
where $\Delta t$ is the time step defined by an opportune CFL condition discussed next.
The time index $n$ varies between $0$ and $N_t$ so that the final time is $t_{\text{final}} = t^{n_t} = t^0 + N_t \Delta t$.

Each equation of system (\ref{eq:DM1}) is solved by a time splitting procedure. 
We recall here a first order splitting approach: 
first a transport step exactly solves the left-hand side, whereas a collision stage solves 
the right-hand side using the solution from the transport step as initial data:
\bea
\label{eq:transport} \textit{Transport stage}  & \longrightarrow & \partial_t f_{k} + v_{k} \cdot\nabla_{x} f_{k} = 0, \\
\label{eq:collision}\textit{Collisions stage} & \longrightarrow & \partial_t f_{k}  = Q(f_k).
\eea
\paragraph{Transport step}
Let $f^{0}_{j,k}$ be the point-wise values at time $t^0$ of the distribution $f$,
${f}^{0}_{j,k}=f(x_{j},v_{k},t^0)$.
The idea behind the fast kinetic scheme is to solve the transport stage (\ref{eq:transport}) continuously in space, see
Fig.~\ref{fig:transport} for a sketch in the one dimensional setting. To this aim we define at the initial time the function $\overline{f}_k^0(x)$ as a piecewise continuous function for all
$x\in \Omega_j$, where $\Omega_{j} = [x_{j-1/2}; x_{j+1/2}]$ and $\Omega=\bigcup_j \Omega_{j}$.
Hence starting from data $\overline{f}_{k,j}^0$ at time index $0$, the exact solution of (\ref{eq:transport})
is simply
\bea
\overline{f}_k^{*,1} = \overline{f}_k^0(x -v_{k} \Delta t), \quad \quad \forall x \in \Omega.
\eea
In other words, the entire function $\overline{f}_k^0$ is advected with velocity $v_k$ 
during $\Delta t$ unit of time and the $*$ superscript indicates that only the transport step has been solved so far. 
The extension of this procedure to the generic time step $n$ gives 
\bea
\overline{f}_k^{*,n+1} = \overline{f}_k^n(x -v_{k} \Delta t), \quad \quad \forall x \in \Omega , \label{eq:f_bar}
\eea
where now, the key observation is that the discontinuities of the piecewise function $\overline{f}_k^n(x)$ do not lie on the interfaces
of two different cells. Instead, the positions of the discontinuities depend entirely on the previous advection step and thus they may be 
located anywhere in the physical space. 
This means that if only the linear transport equation has to be solved, 
this approach gives the exact solution to the equation if the initial data is truly a piecewise constant function  initially 
centered on the spacial mesh.

\paragraph{Collision step}
The effect of the collisional step is to change the amplitude of $\overline{f}_k(x)$.
The idea is to solve the collision operator locally on the grid points, and, successively, extend these computed values to the full domain $\Omega$. 
Thus we need to solve the following ordinary differential equation
\be
\partial_t f_{j,k} =Q(f_{j,k}), 
\label{eq:relax}
\ee
where $f_{j,k}=f(x_j, v_k, t)$,
for all velocities of the lattice $k=1,\ldots,N$ and grid points $j=1,\ldots,M$.
The initial data for solving this equation is furnished by the result of the
transport step obtained by (\ref{eq:f_bar}) at points $x_j$ of the mesh at time $t^{n+1}=t^{n}+\Delta t$, i.e.
$\overline{f}^{*,n+1}_{k}(x_{j})$, for all $k=1,\ldots,N$, and $j=1,\ldots,M$.
Then, the solution of (\ref{eq:relax}), locally on the grid points, reads if, for simplicity, a forward Euler scheme is used as
\be \label{eq:f_coll}
f^{n+1}_{j,k}= f^{*,n+1}_{j,k} \, + \, \Delta t \, Q(f^{*,n+1}_{j,k}),
\ee
where $f^{*,n+1}_{j,k}=\overline{f}^{*,n+1}_{k}(x_{j})$. Many different time integrators can be employed
to solve this equation.
In particular special care is needed in the case in which the equation becomes stiff, refer to \cite{Dimarco_stiff1, Dimarco_stiff2} for 
alternative strategies. Since the time integration of the collision term is not the issue considered in this paper, we considered the simplest possible scheme, but the
FKS technique remains the same when other time integrators are employed. Equation (\ref{eq:f_coll}) furnishes a new value for the distribution $f$ at time
$t^{n+1}$ only in the cell centers of the spacial cells for each velocity $v_k$. However, one needs also the value of the distribution $f$ in all points of
the domain in order to perform the transport step at the next time step. 
Therefore, we define a new piecewise constant function $\overline{Q}_k$ for each velocity of the lattice $v_k$ as
\be
\overline{Q}^{n+1}_{k}(x)=
Q(f^{*,n+1}_{j,k}),
\; \; \forall x \; \mbox{such that} \; \; \overline{f}^{*,n+1}_{k}(x)=\overline{f}^{*,n+1}_{k}(x_j). 
\ee
Said differently we make the fundamental assumption that the shape of $Q(f_k)$ in space coincides
with the one of $f_k$.
Thanks to the above choice one can rewrite the collision step in term of
spatially reconstructed functions as
\be \label{eq:expdt}
\overline{f}^{n+1}_{k}(x) = \overline{f}_{k}(x, t^n+\Delta t) =
  \overline{f}^{*,n+1}_{k}(x)
  + \, \Delta t \, \overline{Q}^{n+1}_{k}(x).
\ee
This ends one time step of the FKS scheme. \\

Concerning the transport part of the scheme, the time step $\Delta t$ is constrained by a CFL like condition of type
\be
\Delta t \max_{k} \left( \frac{|v_{k}|}{\Delta x} \right) \leq 1 = \text{CFL}.
\label{eq:Time}
\ee
The time step constraint for the collision step depends on the choice of the operator $Q$. Since in the numerical test section
we used both a BGK and a Boltzmann operator, the time step constraint for the interaction part has been chosen as the minimum time step which gives stability in the
solution of the ODE (\ref{eq:f_coll}) independently on the type of collision operator employed.
As observed in \cite{FKS} the transport scheme is stable for every choice of the time step, being the solution
for a given fixed reconstruction performed exactly. Nonetheless 
the full scheme being based on a time splitting technique, the error is of the
order of $\Delta t$ in the case of first order splitting or of order $(\Delta t)^q$ for a splitting of order $q$. 
This suggests to take the usual CFL condition for the transport part in order to maintain the time error small enough. \\

To conclude this Section let us observe that time accuracy can be increased by high order time splitting methods,
while spacial accuracy can be increased close to the fluid limit to a nominally second-order accurate scheme by the use
of piecewise linear reconstructions of state variables, see the details in \cite{FKS_HO}. Spacial accuracy for all regime can 
also be increased by using high order polynomial reconstruction for the distribution $f$. This work is in progress and results are discussed in \cite{RFKS}.

\subsection{Conservation of macroscopic quantities} \label{ssec:conservation}
In order to preserve mass, momentum and energy in the scheme, we employ the strategy proposed in \cite{gamba}. We consider one space cell, the same renormalization of $f$ should be considered for all spatial cells. This step is performed at the beginning, i.e. $t=0$. for the distribution $\tilde f(x_j,v_k,t=0)$ and after each collision step (\ref{eq:f_coll}) which causes loss of momentum and energy
due to the spectral discretization. Let $\widetilde{f} = \left(\widetilde{f}_1, \widetilde{f}_2,\ldots,\widetilde{f}_N \right)^{T} $ be the distribution function vector at $t=t^n$ at the center of the cell and let $f = \left(f_1, f_2,\ldots, f_N\right)^{T}$ be the unknown corrected distribution vector which fulfills the conservation of moments. 
Let \[ 
C_{(d_v+2)\times N}=\left(
 \begin{array}{ll}
 & (\Delta v)_k^{d_v}\\
 &  v_k(\Delta v)_k^{d_v} \\
 & |v_k|^{2}(\Delta v)_k^{d_v}\\
 \end{array}
 \right) \] 
be a matrix of coefficients depending on the discretization parameters and $U_{(d_v+2)\times 1} = (\rho \ \rho u \ E)^{T}$ be the vector of conserved quantities, namely density, momentum and energy. Conservation can be imposed solving a constrained optimization formulation:
 \begin{eqnarray}
   &\nonumber  \mbox{ Given } \widetilde{f}\in \mathbb{R}^{N}, \ C \in 
 \mathbb{R}^{(d_v+2)\times N}, \mbox{ and } U \in\mathbb{R}^{(d_v+2)\times 1},\label{eq:minim}\\
   & \mbox{ find } f \in \mathbb{R}^{N} \mbox{ such that } \\
   &\nonumber \|\widetilde{f} -f\|^{2}_{2} \mbox{ is minimized subject to the constrain } Cf = U.
 \end{eqnarray}
The solution of this minimization problem can be analytically obtained by employing the Lagrange multiplier method. 
Let $\lambda\in \mathbb{R}^{d_v+2}$ be the Lagrange multiplier vector. The
corresponding scalar objective function to be minimized is given by
 \be L(f, \lambda) = \sum_{k=1}^{N} |\widetilde{f}_{k}-f_{k} |^{2} +
 \lambda^{T} (Cf- U) . \ee 
Then, by nullifying the derivative of $L(f, \lambda)$ with
respect to $f_{k}$ we get
 \be \label{eq:legrange}
 f = \widetilde{f} + \frac{1}{2} C^{T}\lambda 
 \ee
 while the Lagrange multipliers are obtained by solving
\be
CC^{T}\lambda = 2(U- C\widetilde{f}) .
\ee 
In particular, the above expression says that the value of $\lambda$ is uniquely determined by 
$\lambda= 2(CC^{T})^{-1}(U-C \widetilde{f})$. 
Back substituting $\lambda$ into (\ref{eq:legrange}) finally provides 
\be 
f = \widetilde{f} + C^{T} (CC^{T})^{-1}(U-C\widetilde{f}). \label{eq:minim1} 
\ee 
Note that matrix $C^{T} (CC^{T})^{-1}$ is pre-computed and stored as being constant for each simulation.
If a discretization of the Maxwellian distribution $M[f]$ is needed, as for example when using the BGK model, the same procedure should be applied to the function $M(x_j,v_k,t^n)$ in order
to assure conservation of the first three moments at each instant of time at which the Maxwellian function $M(x_j,v_k,t^n)$ is invoked. In this case by defining with $\E[f]$ the approximated Maxwellian and by $\tilde \E[f]=M(x_j,v_k,t^n)$ the pointwise value of the equilibrium
distribution, we get mimicking (\ref{eq:minim1})
\be \E[U] = \tilde\E[f] + C^{T}
(CC^{T})^{-1}(U-C\tilde\E[f]).\label{eq:minimMax}\ee 
This ends the description of the procedure which permits to conserve macroscopic quantities.

 \subsection{Fast Spectral Scheme (FSS) to discretize the Boltzmann collision operator} \label{sec:Boltz}
The fast spectral discretization of the Boltzmann operator employed in this work is described in this section. 
To this aim, we focus again on a given cell $x_j$ at a given instant of time $t^n$. The same computation is repeated for all cells $x_j, \ j=1,..,M$ and times $t^n, \ n=0,..,N_t$. 
Moreover, since the collision operator acts only on the velocity variable, 
to lighten the notation in this paragraph, 
only the dependency on the velocity variable $v$ is considered for the distribution function $f$, i.e. $f=f(v)$.

The first step to construct for our spectral discretization is to truncate the integration domain of the Boltzmann integral (\ref{eq:bolt}) as done for the distribution $f$. 
As a consequence, we suppose the distribution function $f$ to have compact support on the ball $\Ball_0(R)$ of
radius $R$ centered in the origin. Then, since one can prove that $\supp (Q(f)(v)) \subset \Ball_0({\sqrt 2}R)$, in order to write a spectral approximation which avoid aliasing, it is sufficient that the distribution function $f(v)$ is restricted on the cube $[-T,T]^{d_v}$ with $T \geq (2+{\sqrt 2})R$. Successively, one should assume $f(v)=0$ on $[-T,T]^{d_v} \setminus \Ball_0(R)$ and extend $f(v)$ to a periodic function on the set $[-T,T]^{d_v}$. Let observe that the lower bound for $T$ can be improved. For instance, the choice $T=(3+{\sqrt 2})R/2$ guarantees the absence of intersection between periods where $f$ is different from zero. However, since in practice the support of $f$ increases with time, we can just minimize the errors due to aliasing \cite{canuto:88} with spectral accuracy. 

To further simplify the notation, let us take $T=\pi$ and hence $R=\lambda\pi$ with $\lambda = 2/(3+\sqrt{2})$ in the following. 
We denote by $\QL_B(f)$ the Boltzmann operator with cut-off. Hereafter, using one index to denote the $d_v$-dimensional sums, we have that the approximate function $f_N$ can be represented as the truncated Fourier
series by
\begin{equation}
f_N(v) = \sum_{k=-N/2}^{N/2} \f_k e^{i k \cdot v},
\label{eq:FU}
\end{equation}
\begin{equation}
\f_k = \frac{1}{(2\pi)^{d_v}}\int_{[-\pi,\pi]^{d_v}} f(v)
e^{-i k \cdot v }\,dv.
\label{eq:FC}
\end{equation}
We then obtain a spectral quadrature of our collision operator by projecting (\ref{eq:bolt}) on the space
of trigonometric polynomials of degree less or equal to $N$, i.e.
\begin{equation}
{\hat Q}_k=\int_{[-\pi,\pi]^{d_v}}
\QL_B(f_N)
e^{-i k \cdot v}\,dv, \quad k=-N/2,\ldots,N/2. 
\label{eq:VAR}
\end{equation}
Finally, by substituting expression (\ref{eq:FU}) in (\ref{eq:VAR}) one gets after some computations
\begin{equation}
{\hat Q}_k = \sum_{\substack{l,m=-N/2\\l+m=k}}^{N/2} \f_l\,\f_m
\bb(l,m),\quad k=-N,\ldots,N,
\label{eq:CF1}
\end{equation}
where $\bb(l,m)=\B(l,m)-\B(m,m)$ are given by
\begin{equation}
\B(l,m) = \int_{\Ball_0(2\lambda\pi)}\int_{\ens{S}^{d_v-1}} 
|q| \sigma(|q|, \cos\theta) e^{-i(l\cdot q^++m\cdot q^-)}\,d\omega\,dq. \label{eq:KM}
\end{equation}
with \begin{equation}
q^{+} = \frac12(q+\vert q\vert \omega), \quad
q^{-} = \frac12(q-\vert q\vert \omega).
\label{eq:VV2}
\end{equation}
Let us notice that the naive evaluation of (\ref{eq:CF1}) requires $O(n^2)$ operations, where $n=N^3$. This causes the spectral method to be computationally very expensive,
especially in dimension three. In order to reduce the number of operations needed to evaluate the collision integral, the main idea is to use another representation of \eqref{eq:bolt}, the so-called Carleman representation \cite{Carl:EB:32} which is obtained by using the  following identity
  \begin{equation*}
    \frac{1}{2} \, \int_{\SSS^{d_v-1}} F(|u|\sigma - u) \, d\sigma	 = \frac{1}{|u|^{d-2}} \, \int_{\RR^{d_v}} \delta(2 \, x \cdot u + |x|^2) \, F(x) \, dx.
 \end{equation*} 
This gives in our context for the Boltzmann integral
\begin{equation}
\label{defQBCarleman}
 Q_B (f)= \int_{\R^{d_v}} \int_{\R^{d_v}} {\tilde B}(x,y) 
 \delta(x \cdot y) 
 \left[ f(v + y) \, f(v+ x) - f(v+x+y) \, f(v) \right] \, dx \,
 dy,
 \end{equation} 
with 
  \begin{equation}\label{eq:Btilde}
  \tilde{B}(|x|,|y|) =
  2^{d_v-1} \, \sigma\left(\sqrt{|x|^2+|y|^2}, \frac{|x|}{\sqrt{|x|^2+|y|^2}} \right) \, (|x|^2+|y|^2)^{-\frac{d_v-2}2}.
  \end{equation}
This transformation yields the following new spectral quadrature formula 
 \begin{equation}\label{eq:ode}
 \hat{Q}_k  =
 \sum_{\underset{l+m=k}{l,m=-N/2}}^{N/2} {\hat{\beta}}_F(l,m) \, \hat{f}_l \, \hat{f}_m, \ \ \
 k=-N,...,N
 \end{equation}
where ${\hat{\beta}}_F(l,m)=\B_F(l,m)-\B_F(m,m)$ are now 
given by
 \begin{equation}
 \B_F(l,m) = \int_{\Ball_0(R)} \int_{\Ball_0(R)}
 \tilde{B}(x,y) \, \delta(x \cdot y) \, 
 e^{i (l \cdot x+ m \cdot y)} \, dx \, dy.
 \label{eq:FKM}
 \end{equation}
Now, in order to reduce the number of operation needed to evaluate~\eqref{eq:ode}, we look for a convolution structure. 
The aim is to approximate each
${\hat{\beta}}_F(l,m)$ by a sum
 \[ {\hat{\beta}}_F(l,m) \simeq \sum_{p=1} ^{A} \alpha_p (l) \alpha' _p (m), \]
where $A$ represents the number of finite possible directions of collisions.
This finally gives a sum of $A$ discrete convolutions and, consequently, the algorithm can be computed in $O(A \, N \log_2 N)$ operations by means of
standard FFT technique~\cite{canuto:88}. 

In order to get this convolution form, we make the decoupling assumption
 \begin{equation}\label{eq:decoup}
 \tilde{B}(x,y) = a(|x|) \, b(|y|).
 \end{equation}
This assumption is satisfied if $\tilde B$ is constant. This is the case of Maxwellian molecules in dimension two, and hard spheres in dimension three, the two cases treated in this paper. 
Indeed, using kernel \eqref{defVHSKernel} in \eqref{eq:Btilde}, one has
\[
  \tilde{B}(x,y) = 2^{d_v - 1} C_\alpha (|x|^2+|y|^2)^{-\frac{d_v-\alpha-2}2},
\] 
so that $\tilde{B}$ is constant if $d_v = 2$, $\alpha = 0$ and $d_v = 3$, $\alpha = 1$.

We start by dealing with dimension $2$ and $\tilde{B} =1$, i.e. Maxwellian molecules. Here we write $x$ and $y$ in spherical coordinates $x = \rho e$ and $y = \rho' e'$ to get
 \begin{equation*}
 \B_F(l,m) = \frac14 \, \int_{\ens{S}^1} \int_{\ens{S}^1}
 \delta(e \cdot e') \,
 \left[ \int_{-R} ^R e^{i \rho (l \cdot e)} \, d\rho \right] \,
 \left[ \int_{-R} ^R e^{i \rho' (m \cdot e')} \, d\rho' \right] \, de \,
 de'.
 \end{equation*}
Then, denoting $ \phi_R ^2 (s) = \int_{-R} ^R e^{i \rho s} \, d\rho,$ for $s \in \R$,  
we have the explicit formula
 \[ \phi_R ^2 (s) = 2 \, R \,{\sinc} (R s), \]
 where ${\sinc}(x)=\frac{\sin(x)}{x}$.
This explicit formula is further plugged in the expression of $\B_F(l,m)$ and using its parity property, this yields
 \begin{equation*}
 \B_F (l,m) =  \int_0 ^{\pi} \phi_R ^2 (l \cdot e_{\theta})\, \phi_R ^2 (m \cdot e_{\theta+\pi/2}) \,
 d\theta.
 \end{equation*}
Finally, a regular discretization of $A$ equally spaced points, which is spectrally accurate because of the periodicity of the function, gives
 \begin{equation}
 \B_F (l,m) = \frac{\pi}{M} \, \sum_{p=1} ^{A} \alpha_p (l) \alpha' _p (m),
 \label{eq:dfor0}
 \end{equation}
with
 \be \alpha _p (l) = \phi_R ^2 (l \cdot e_{\theta_p}), \hspace{0.8cm} \alpha' _p (m) = \phi_R ^2 (m \cdot e_{\theta_p+\pi/2}) 
 \ee
where $\theta_p = \pi p/A$.

Now let us deal with dimension $d=3$ with $\tilde{B}$ satisfying
the decoupling assumption~\eqref{eq:decoup}. First we switch to
the spherical coordinates for $\B_F(l,m)$:
 \begin{equation*}
 \B_F(l,m) = \frac14 \int_{\ens{S}^2\times\ens{S}^2}
 \delta(e \cdot e') 
  \left[\int_{-R} ^R \rho  a(\rho) e^{i \rho (l \cdot e)}  d\rho  \right]
  \left[\int_{-R} ^R \rho' b(\rho') e^{i \rho' (m \cdot e')} d\rho' \right]  de  de'.
 \end{equation*}
Then, integrating first $e'$ on the intersection of the unit
sphere with the plane $e^\bot$ gives
 \begin{equation*}
 \B_F (l,m) = \frac14 \, \int_{e \in \ens{S}^2} \phi_{R,a} ^3 (l \cdot e) \, \left[
 \int_{e' \in \ens{S}^2 \cap e^\bot} \phi_{R,b} ^3 (m \cdot e') \, de' \right] \, de
 \end{equation*}
where
 \begin{equation*}
 \phi_{R,a} ^3 (s) = \int_{-R} ^R \rho  \, a(\rho) \, e^{i \rho s} \, d\rho.
 \end{equation*}
This leads to the following decoupling formula 
 \begin{equation*}
 \B_F (l,m) = \int_{e \in \ens{S}^2 _+} \phi_{R,a} ^3 (l \cdot e) \,
 \psi_{R,b} ^3 \big(\Pi_{e^\bot}(m)\big) \, de
 \end{equation*}
where $\ens{S}^2_+$ denotes the half sphere and
 \begin{equation*}
 \psi_{R,b} ^3 \big(\Pi_{e^\bot}(m)\big) = \int_0 ^\pi \sin\theta \, \phi_{R,b} \big(|\Pi_{e^\bot}(m)|
 \, \cos\theta \big)  \, d\theta.
 \end{equation*}
Now, in the particular case where $\tilde{B}=1$, i.e. the hard sphere
model, we can explicitly compute the functions $\phi_R ^3$. These are
 \begin{equation*}
 \phi_R ^3 (s) = R^2 \, \left[ 2 \, \mbox{Sinc} (R s) - \mbox{Sinc} ^2 (R s /2) \right], \hspace{0.8cm}
 \psi_R ^3 (s) = 2 \, R^2 \, \mbox{Sinc} ^2 (R s /2).
 \end{equation*}
 Taking a
spherical parametrization $(\theta,\varphi)$ of $e \in \ens{S}^2
_+$ and uniform grids of respective size $A_1$ and $A_2$ for
$\theta$ and $\varphi$ (again spectrally accurate because of the periodicity of the function) leads to the following quadrature formula for $\B_F(l,m)$
 \begin{equation*}
 \B_F (l,m) = \frac{\pi^2}{A_1 A_2} \, \sum_{p,q=0} ^{A_1,A_2} \alpha_{p,q} (l) \alpha' _{p,q} (m)
 \end{equation*}
where
 \[ \alpha_{p,q} (l) = \phi_{R,a} ^3 \left(l \cdot e_{(\theta_p,\varphi_q)} \right),
 \hspace{0.8cm} \alpha' _{p,q} (m) = \psi_{R,b} ^3 \left(\Pi_{e_{(\theta_p,\varphi_q)} ^\bot} (m) \right), \]
 \[ \phi_{R,a} ^3 (s) = \int_{-R} ^R \rho \, a(\rho) \,  e^{i \rho s} \, d\rho,
 \hspace{0.8cm} \psi_{R,b} ^3 (s) = \int_0 ^\pi \sin \theta \, \phi_{R,b} ^3 (s \cos \theta) \, d\theta, \]
and for all $p$ and $q$
 \[ (\theta_p,\varphi_q) = \Big(\frac{p \, \pi}{A_1}, \frac{q \, \pi}{A_2} \Big). \]
       

%
%
\section{Numerical implementation} \label{sec:implementation}
We discuss in this part some aspects relative to the numerical implementation of the method. 

\subsection{Algorithm for the transport part}  \label{ssec:algorithm}
The method described in the previous sections can be resumed into two main actions: transport and interaction.
In our implementation we have adopted a particle like interpretation of the FKS scheme which helps in reducing the computational effort and the memory requirement.
In this interpretation, each point of the quadrature of the phase space is represented by a particle with a given mass, velocity and position. The transport phase causes the particle to move in the physical space, while the interaction phase causes the mass of each particle to change. 
In this setting, the distribution function $f$ can be expressed as
\begin{eqnarray}
& & f(x,v,t) = \sum_{i=1}^{N_p} \mathfrak{m}_i(t) \, \delta(x-x_i(t))\delta(v-v_i(t)), \quad v_i(t)=v_k,
\label{particle2}
\end{eqnarray}
where $x_i(t)$ represents the particle position, $v_i(t)$ the
particle velocity, $\mathfrak{m}_i$ the particle mass. Moreover, the particle velocity corresponds 
to the quadrature point chosen to discretize the 
velocity space. The effect of the transport is a simple shift of particles to their new positions
according to \be x_i(t+\Delta
t)=x_i(t)+v_i(t)\Delta t. \label{transport}\ee 
The collision step acts only locally and changes the velocity
distribution. The particle interpretation of this part of the scheme consists in changing the mass of each particle
through the formula (\ref{eq:f_coll}) which reads \be \label{eq:part_coll}
\mathfrak{m}_i(t+\Delta t)= \mathfrak{m}_i(t) \, + \, \Delta t \, Q(v_i),
\ee
where $Q(v_i)$ corresponds to the approximation of the collision integral in the center of the cell evaluated at location $v_i$ by means of the fast spectral method described before. 
Let observe that thanks to the uniform grid in velocity and physical space, the number of particles is the same in each cell and remains constant in time. 
This permits to consider the motion of only a fixed subset of particles which belongs to only one chosen cell. 
The relative motion of all other particles in their cells being the same \cite{FKS}. 
This gives great computational advantage since the number of particle to track in time is greatly reduced. 
The only information which should be tracked remains the mass
of each individual particle which is different for every velocity $v_i$ and position $x_i$. 
This reinterpretation of the scheme permits to strongly reduce the computational cost related to the transport part
as shown in the numerical test section.

\subsection{Profiling and parallelization strategy}  \label{ssec:parallel}
The Fast Kinetic scheme for the BGK equation was shown to be very efficient and extremely parallelizable (with a shown acceleration very close to the ideal scaling) in \cite{FKS_GPU} on mild parallel architecture, namely using OpenMP on maximum $64$ threads and on two Graphical Processing Units (GPU) using CUDA framework. 
These light parallel infrastructures are unfortunately not sufficient when Boltzmann operator is to be simulated. 
To this aim, let us briefly discuss the profiling of the collision algorithm. 
As we have already seen in the description of the method in the previous section, the resolution of Boltzmann operator requires several passages 
from the physical space to the Fourier space, which implies several calls to (inverse) Fourier transforms 
within nested loops which are due to the convolutive form of the integral. 
This interaction step covers about $98\%$ of the total cost of a 2D$\times$2D simulation on a serial machine, see table~\ref{tab:cost}. 
The situation is  alike in the case of 3D$\times$3D simulations, and, thus, results are not reported. 
From table~\ref{tab:cost} we deduce that the collision step is the part of the code which demands to be dealt with great attention
if any gain in performance is expected.
\begin{table}
  \begin{center}
  \numerikNine
  \begin{tabular}{|c||c|cc|}
    \hline
    \textbf{Main routines} &  \textbf{Cost CPU \%}  &  \textbf{Sub-routines} &   \textbf{Sub-cost CPU \%} \\
    \hline
    \hline
    Transport      &  $\simeq$2  \% & --- &  $\simeq$2  \%\\
    \hline
    Update moments &  $\simeq$0  \% & --- & $\simeq$0  \% \\
    \hline
    \multirow{4}{*}{Collision}     &  \multirow{4}{*}{$\geq$  98 \%}
    & Compute $\mathcal{F}(f)$                   &   $\simeq$5   \%  \\
    & & $\theta_p$-Loop: compute coeff.            &   $\simeq$5   \%  \\
    & & $\theta_p$-Loop: compute $\mathcal{F}^{-1}$ &  $\simeq$87   \%  \\
    & & Update $Q$                               &   $\simeq$3   \%  \\
    \hline
    \hline
    Total & 100 \%& & 100\% \\
    \hline
  \end{tabular}
  \caption{ \label{tab:cost} Schematic table of the average cost
  for each routine of the code on a cylindrical Sod like
  problem in 2D$\times$2D. Simulation is run on a serial machine. 
  The collision part demands almost all the CPU resources and, within 
  this step, the calls to the (inverse) Fourier transform routine 
  covers more than $85\%$ of the total cost of one single step.
  }
  \end{center}
\end{table}
More precisely, if we give a closer look to the computational cost related to the computation of the collision term $Q(f)$ in one fixed spacial cell, 
the major cost (beyond $90\%$) is spent by the Fourier transform routines. In our implementation, the collision is performed without
any communication with the local neighborhood. In other words, once the $N$ values of the distribution functions $f_{j,k}, \ k=1,..,N$ are stored then
$Q(f)$ is computed locally. Because $N$ is not too large, typically the number of points chosen to discretize the velocity space goes from $N=16^3$ to
$N=64^3$, then the memory requirement is not extremely demanding. Finally, as seen in the previous Section the computation of $Q(f)$ is based on
a loop over the number of angles $\theta_p$ with $p=1,..,A$ in dimension two, and, a double loop over the angles $\theta_p$ and $\phi_q$ 
in dimension three with $p=1,..A_1$
and $q=1,..,A_2$ with a typical number of chosen angles of eight. After a first Fourier transform of the distribution $f$, for each iterate of these loops, 
the coefficients of the quadrature formula are computed and then three calls to the inverse  fast Fourier transform are done. 
The results of the FFT are gathered into temporary complex variables later arranged to form the solution of $Q(f)$. 
The costs related to the four main routines needed for the evaluation of the operator namely: 
Fourier transform, coefficients computation, inverse Fourier transform and sum of the obtained values to form $Q(f)$, are detailed in table ~\ref{tab:cost}.


We discuss now some aspects related to the parallelization strategy. The framework described is extremely well suited for parallelization. The strategy adopted is
to divide the space domain in several subdomains because the interaction step, the most expensive one, does not demand any communication between spatial cells. 
Each cell can be therefore treated independently from the others. Then, the idea employed has been to use the OpenMP on shared memory systems and make each computational core responsible for a subset of a space mesh. This approach requires very little modification of a sequential code and gives strong scaling close to ideal as shown in \cite{FKS_GPU} for the BGK
case. A different possibility is to compute the collision kernel, and all the related Fast Fourier Transforms on a GPU. Even if, this involves substantial amount of slow communications between CPU and GPU, the computational complexity of this part is so large that the communication time does not play a fundamental role. This approach requires more programming effort, compared to the OpenMP parallelization strategy as several routines have to be completely rewritten.
In particular, the code has to be adapted to the Single Instruction,
Multiple Threads (SIMT) parallel programming mode imposed by the GPU
architecture. That is to say, parts of the code that can be run in
parallel on large number of threads have to be identified and
rewritten as so called CUDA kernels - functions that are executed
multiple times in parallel by different CUDA threads. Moreover, the
GPU-CPU memory data transfer has to be taken into account.
The payoff is however higher and, at least for computations which does not require too large memory storage, this is probably the best choice. 

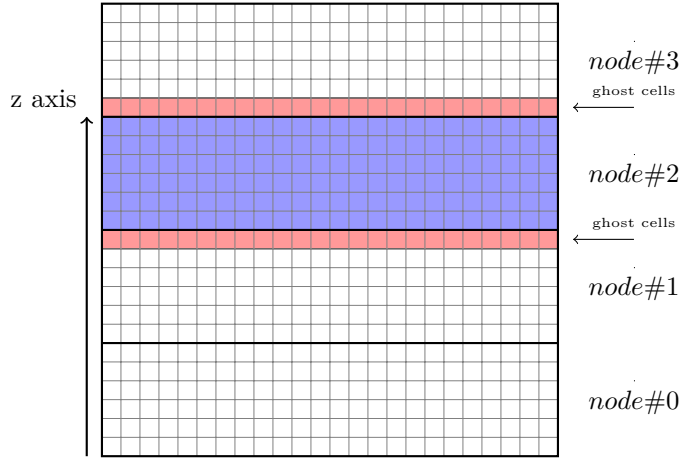
\begin{figure}
  \begin{center}
    \def\xmax{6}
    \def\xwidth{6/4}
    \def\xswidth{1/3*\xmax/8}
    
    \begin{tikzpicture}
      
      \filldraw[fill=red!40!white, draw=black] (0,2*\xwidth - \xswidth) rectangle (\xmax,2*\xwidth);
      \filldraw[fill=blue!40!white, draw=black] (0,2*\xwidth) rectangle (\xmax,3*\xwidth);
      \filldraw[fill=red!40!white , draw=black] (0,3*\xwidth) rectangle (\xmax,3*\xwidth+\xswidth);
      
      \draw[step=\xswidth,gray,very thin] (0,0) grid (\xmax,\xmax);
      \draw[thick] (0,0) rectangle (\xmax,\xmax);
      
      \foreach \y in {1,2,3}
      \draw[thick] (0,\y *\xwidth) -- (\xmax,\y *\xwidth);
      
      \draw[thick,->] (-0.2,0) -- (-0.2,4.5) node[anchor=south east] {z axis};
      
      \foreach \y in {0,1,2,3}
      \draw (\xmax+1,\xwidth*\y +1) -- (\xmax+1 ,\xwidth*\y +1) node[anchor=north] {$node \#\y$};
      
      \draw[thin,<-] (\xmax +0.2,3*\xwidth+\xswidth/2) -- (\xmax+1,3*\xwidth+\xswidth/2) node[anchor=south] {\tiny ghost cells};
      \draw[thin,<-] (\xmax +0.2,2*\xwidth-\xswidth/2) -- (\xmax+1,2*\xwidth-\xswidth/2) node[anchor=south] {\tiny ghost cells};
      
    \end{tikzpicture}
  \end{center}

  \caption{Domain decomposition for the MPI parallelization. The
    domain slices are distributed over computational nodes. After
    the collision step the new information brought by the particles on the cells located on the
    subdomain boundary are communicated to the ghost cells of the
    neighboring nodes.}
  \label{fig:mpi-decomp}
\end{figure}

Unfortunately, when the number of points needed to perform a simulation grows in the case the Boltzmann operator is used, the above implementation strategies on shared memory systems are not sufficient any more. This is true particularly in the full three dimensional case. Luckily, the method described in this work can be easily adapted for its application in distributed memory systems. The idea adopted in this case is to distribute spatial degrees of freedom over computational nodes, keeping on every node a complete set of velocity space points. The computational domain is then
decomposed into slices along the $Z$ direction (see Fig.\ref{fig:mpi-decomp}). The slices are successively distributed over different MPI processes. Every node will therefore compute the
collision term for a part of computational domain. This can be done on each node using the OpenMP or GPU parallelization strategy, depending on the architecture at disposal. As the update of the primitive variables requires an exchange with neighboring spatial cells, some of the particle masses, which contain all the information related to the distribution function $f$ in our implementation, need to be broadcast to other computational nodes. This data exchange, which is typically a bottleneck in the MPI computations, and, usually demands large efforts in order to be minimized the internodal communications, in the case of the Boltzmann equation, takes only a small part of the total runtime, even if, in our first MPI implementation, the domain decomposition is far from being optimal. We postpone to a future work the discussion related to a more efficient MPI parallelization strategy and all details related to this type of implementation.


%
%
\section{Numerics} \label{sec:numerics}
In this section, we validate the proposed SFKS (Spectral-FKS) method which couples the fast semi-Lagrangian  and the fast spectral scheme. 
The testing methodology is divided into four parts. 
\begin{description}
\item
  \textbf{Part 1.} Sanity checks. We only consider the space homogeneous Boltzmann equation in two and three dimensions in velocity space.
  Two problems from \cite{PR2} for which the analytical solution is known are simulated. The results show that our implementation
  of the fast spectral discretization of the Boltzmann operator is correct. 
\item
  \textbf{Part 2.} Using the SFKS method in one dimension of physical space and two or three dimensions in velocity space,
we show on Riemann like problems the differences between the BGK and the Boltzmann models. 
The reported results justify the use of the more complex and costly collisional Boltzmann operator especially in situations far from the thermodynamical equilibrium.
We also report the details of the computational costs as well as the numerical convergence of the method for an increasing
number of points in the physical or velocity space.
\item
  \textbf{Part 3.} Using the SFKS method in the two dimensional case in space and velocity, we show that using the Boltzmann or the BGK operator leads to notably different
results. We consider a regular vortex like problem and a reentry like problem with evolving
angle of attack.
Performances and profiling study of the SFKS approach supplemented with BGK or Boltzmann collisional operator are provided for refined grids in space and time.
Details about the computational costs of the different part of the scheme, scalability with  respect to the number of cells
are also furnished in order to characterize the behavior of the method as precisely as possible.
\item
  \textbf{Part 4.} We test the SFKS method in the full three dimensional case in space and velocity for an unsteady test problem.
The interaction of a flying object with the surrounding ambient gas is simulated. 
We compare the results obtained with the Boltzmann collisional operator with those given by a simpler BGK model. 
Moreover, we measure the cost of such a 3D$\times$3D simulation in terms
of CPU, cost per degree of freedom and scalability with respect to the number of cells. Some code profiling is also provided to measure the cost of the main components of the scheme.
\end{description}
Finally, the ability of the whole numerical code to run on several parallel environments is tested:
CPU (OpenMP, MPI) and GPU (CUDA) types of parallelism frameworks are considered.

%
%
\subsection{Part 1. Numerical results for the space homogeneous case} \label{sec:act0}
In this part, we validate our spectral discretization implementation of the Boltzmann collisional operator. 
To this aim, we re-employ the test cases considered in \cite{PR2} and \cite{FilbetRusso} for the homogeneous two dimensional and three dimensional Boltzmann equation. 
For both situations, exact solutions are available and briefly recalled in the following. 

\subsubsection{Test 1.1. Convergence to equilibrium for the Maxwell molecules in dimension two.}
\begin{figure}
\begin{center}
    \includegraphics[width=0.48\textwidth]{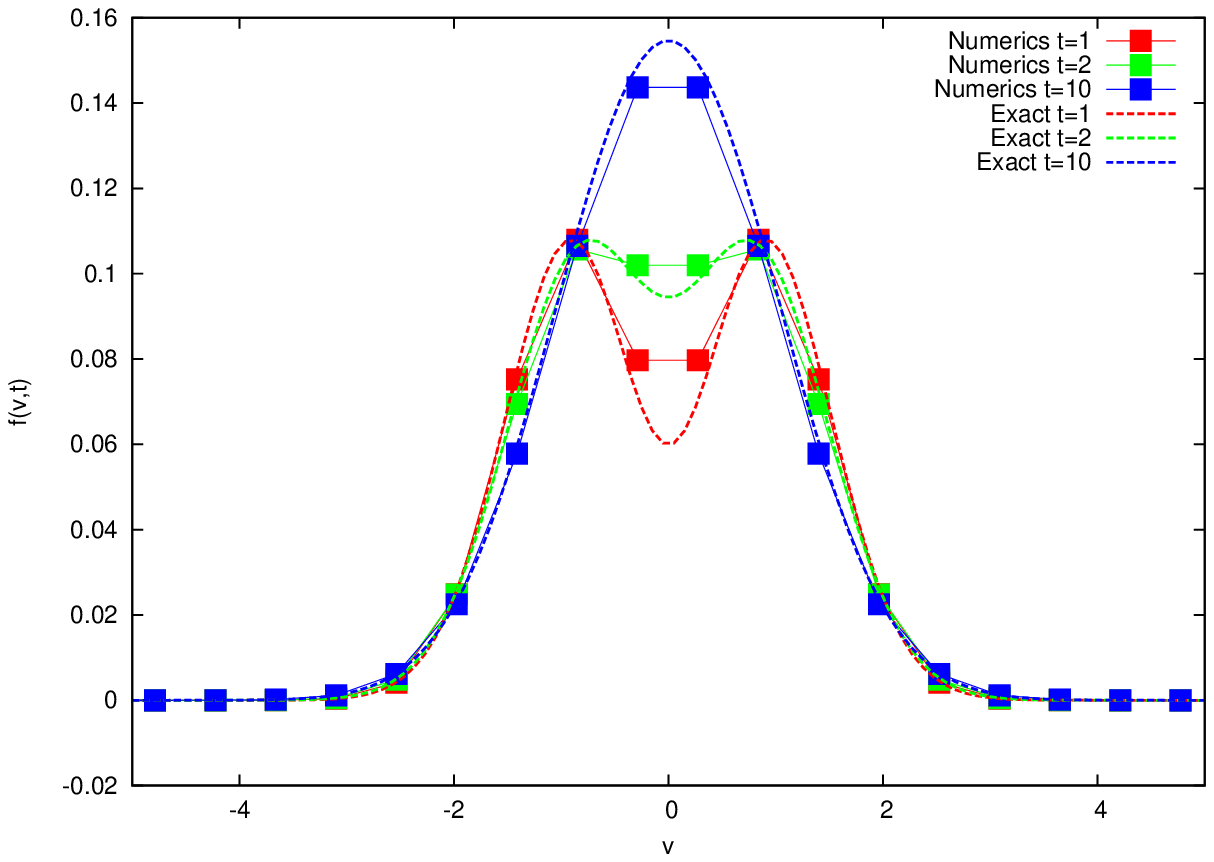} 
    \includegraphics[width=0.48\textwidth]{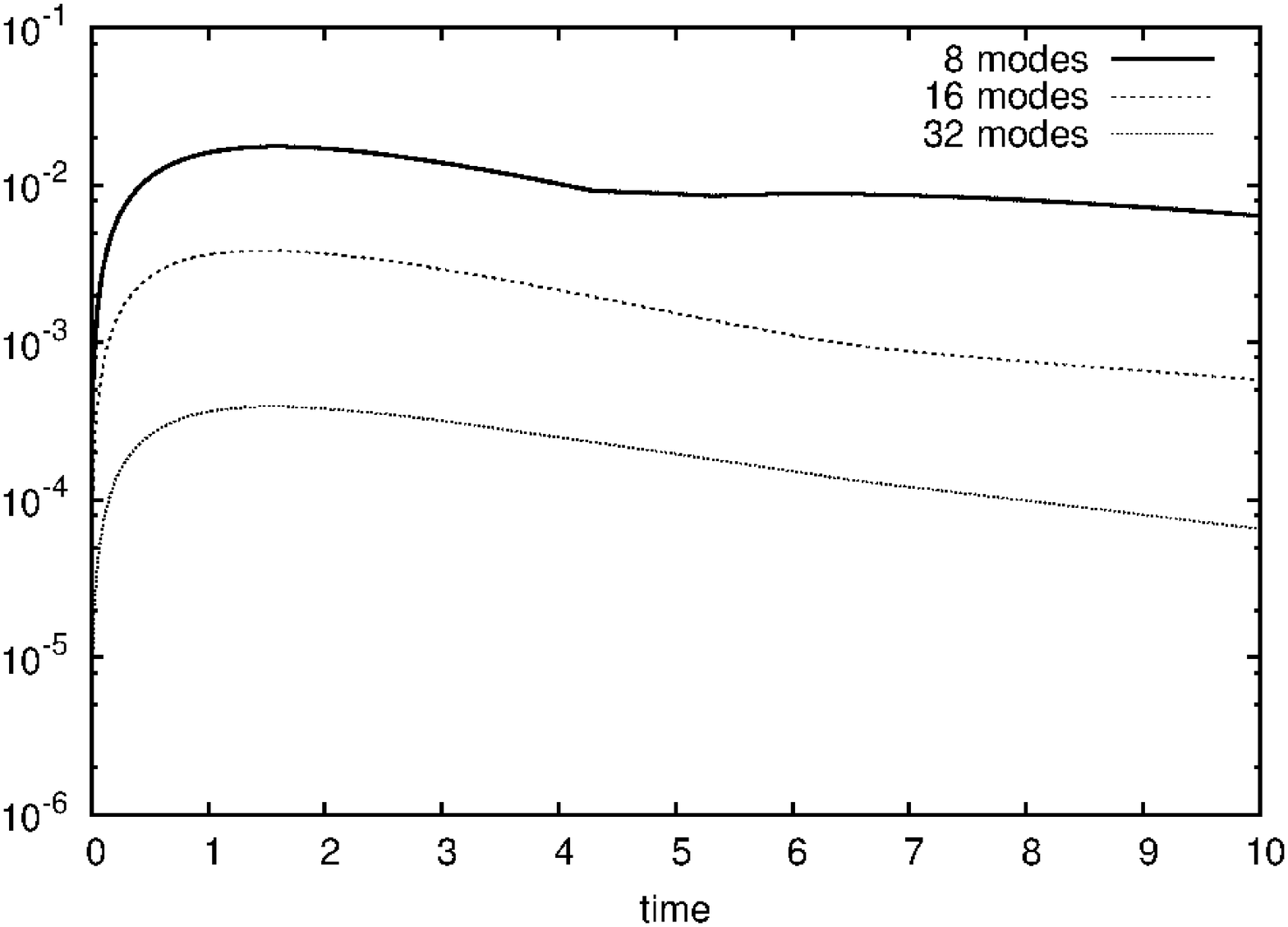}
  \caption{Test 1.1. Sanity check --- 
    Left: $f(v_x,v_y=0,t)$ (symbols and straight line) as a function of velocity $v_x$
    and time for $t=1$, $2$ and $t_{\text{final}} =10$ versus the exact solution (dashed line).
    Right: 
    $L_1$ error as a function of time for $8$, $16$ and
    $32$ points (or modes) in each direction. Maxwellian molecules in dimension two.}
  \label{fig:test1}
  \end{center}  
\end{figure}

We consider two dimensional in velocity Maxwellian molecules and the following space homogeneous initial condition 
\begin{gather}
  f(v,t=0) = \frac{v^2}{\pi} \exp (-v^2) .
  \label{eq:BKW}
\end{gather}
The analytical solution (the so-called BKW one) is given for all times $t$  by \cite{Bobylev:75,KrookWu:1977}:
\begin{gather}
  f(v,t)=
  \frac{1}{2\pi S^2}
  \exp(-v^2/2S)
  \left[
    2S-1+ \frac{1-S}{2S}v^2
  \right], \quad \forall t>0,
  \label{eq:Jngb}
\end{gather}
with $S=1-\exp (-t/8)/2$. The final time is set to $t_{\text{final}} = 10$ and the
time step is equal to $\Delta t = 0.02$. 
The test is performed for $N=8^2, 16^2, 32^2$ points and for eight discrete angles $\theta_p$. 
In order to keep the aliasing error smaller than the spectral error, the velocity domain is set to $[-4,4]^2$ for $N=8^2$, $[-6,6]^2$ for $N=16^2$ and
to $[-9,9]^2$ for $N=32^2$ points. This choice helps fighting back the aliasing error behavior since it increases with the number of points. 
Therefore, to have comparable
aliasing errors for the three meshes, we enlarge the size of the velocity space when the number of points increases. 
We report the discrete $L_1$ and $L_2$ norms of the error for the distribution function $f$ in Table \ref{tab:test1}. Moreover, in Figure \ref{fig:test1}, 
we report the distribution function 
$f(v_x,v_y=0,t)$ at different times ($t=1$, $2$ and $10$) when $N=32$. The solution is plotted  versus the exact solution.
The results clearly show the convergence of the method towards the exact solution.
Last in Figure~\ref{fig:test1} right panel, we present the time evolution  of the $L_1$ error for the three configurations. Consistently the
errors decrease by about one order of magnitude. The last recorded results, at $t_{\text{final}} = 10$ corresponds to the Figures reported in table~\ref{tab:test1}.
\begin{table}
  \centering
  \begin{tabular}{|c||c|c|}
    \hline
    $\#$ of points & $L_1$ error & $L_2$ error \\
    \hline\hline
    $8^2$ &  $6.3\times 10^{-3}$ & $3.8\times 10^{-3}$  \\
    \hline
    $16^2$ &  $5.7\times 10^{-4}$ & $2.9\times 10^{-4}$ \\
    \hline
    $32^2$ &  $6.5\times 10^{-5}$ & $7.2\times 10^{-5}$  \\
    \hline
  \end{tabular}
  \caption{Test 1.1. Sanity check ---
    $L_1$ and $L_2$ relative errors for the fast spectral method for
    different number of points at $t_{\text{final}} = 10$. Maxwellian molecules in dimension two.}
  \label{tab:test1}
\end{table}

\subsubsection{Test 1.2. Space homogeneous comparison between the BGK and the Boltzmann model. Maxwellian molecules.}
\begin{figure}
  \begin{center}
    \includegraphics[width=0.6\textwidth]{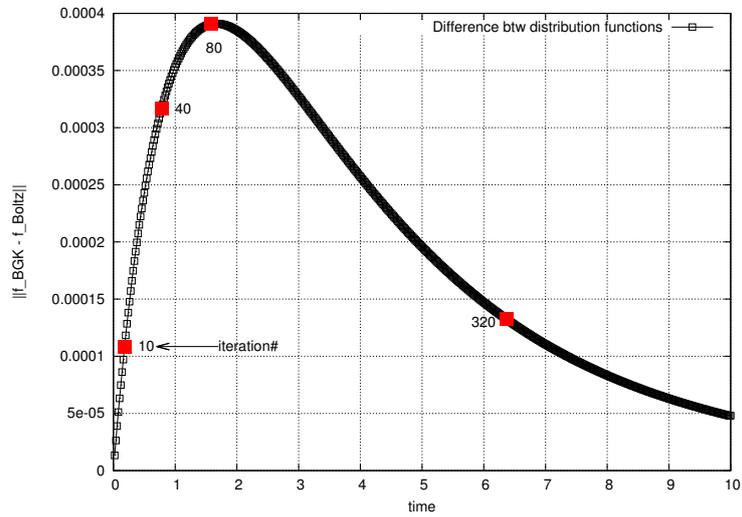}
  \end{center}
  \caption{Test 1.2. Differences between the distribution functions obtained with the BGK and the Boltzmann models. Maxwellian molecules in two dimensions.
    The distribution function values for the red marked iterations are depicted in Figure~\ref{fig:test1_comp2}.
  }
  \label{fig:test1_comp}
\end{figure}

\begin{figure}
\begin{center}
    \includegraphics[width=0.49\textwidth]{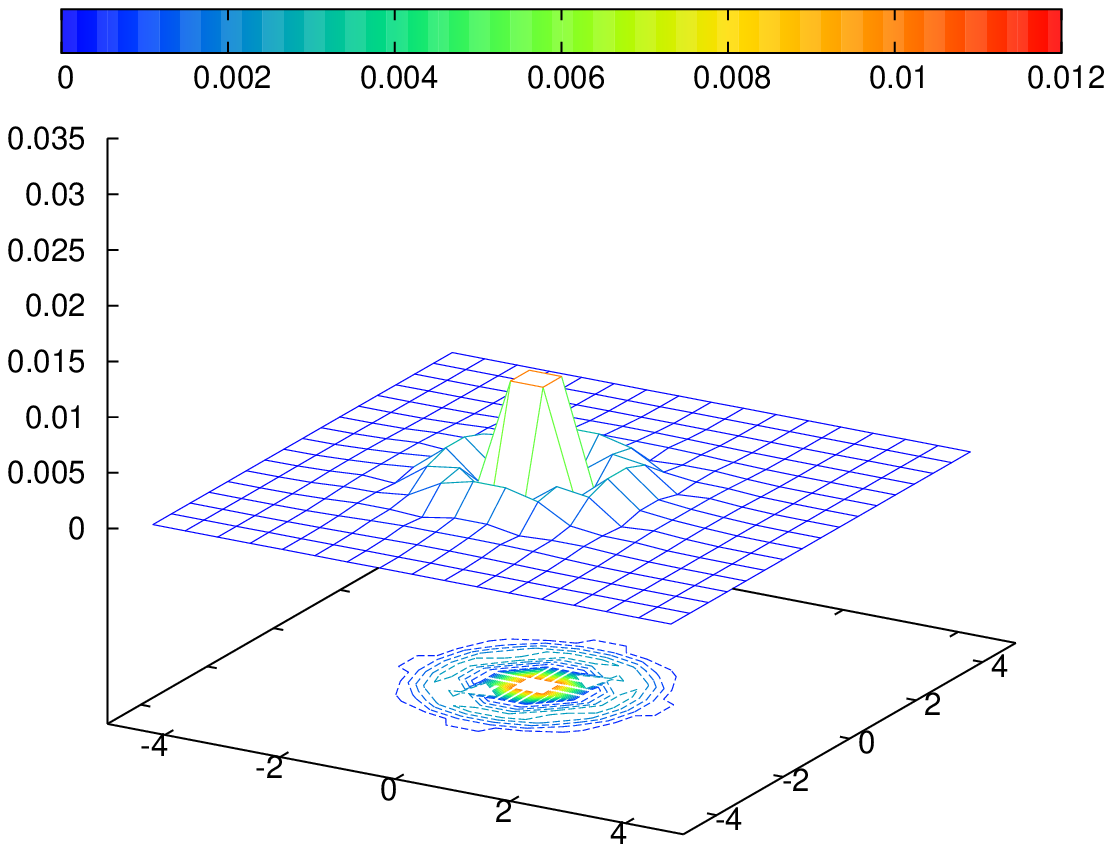}
    \includegraphics[width=0.49\textwidth]{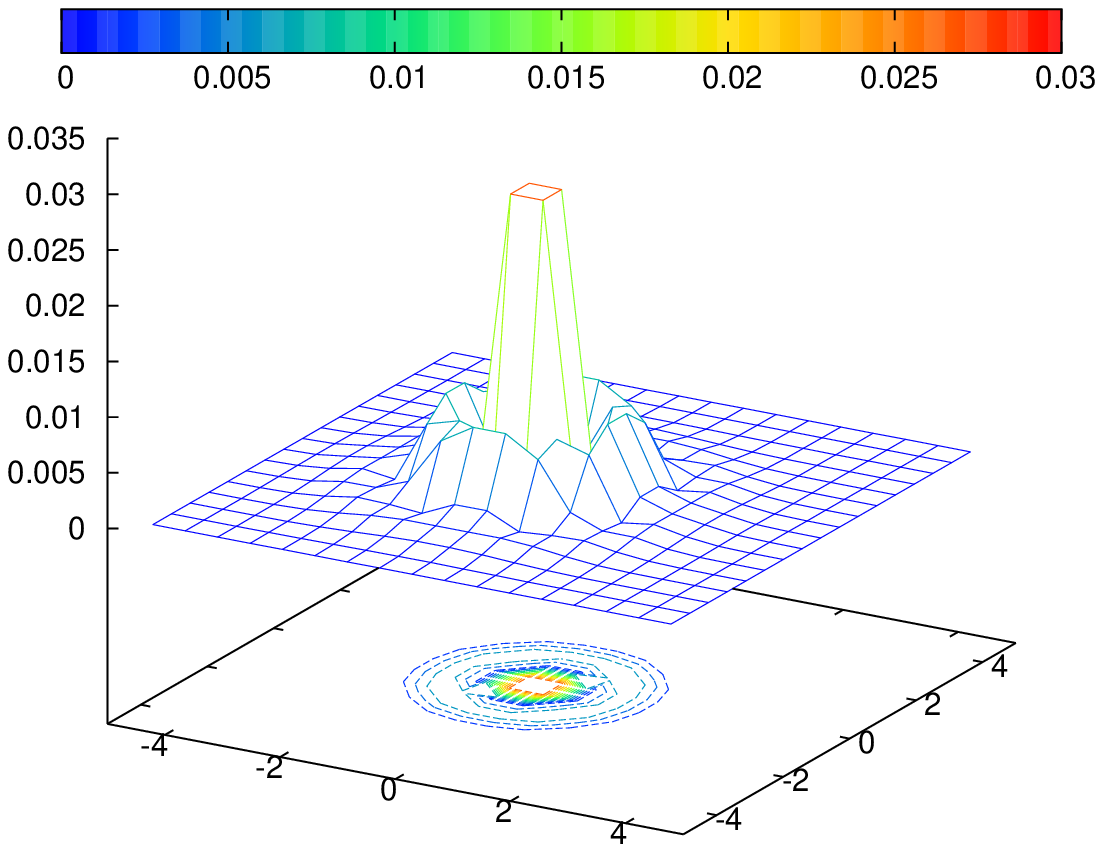}\\
    \includegraphics[width=0.49\textwidth]{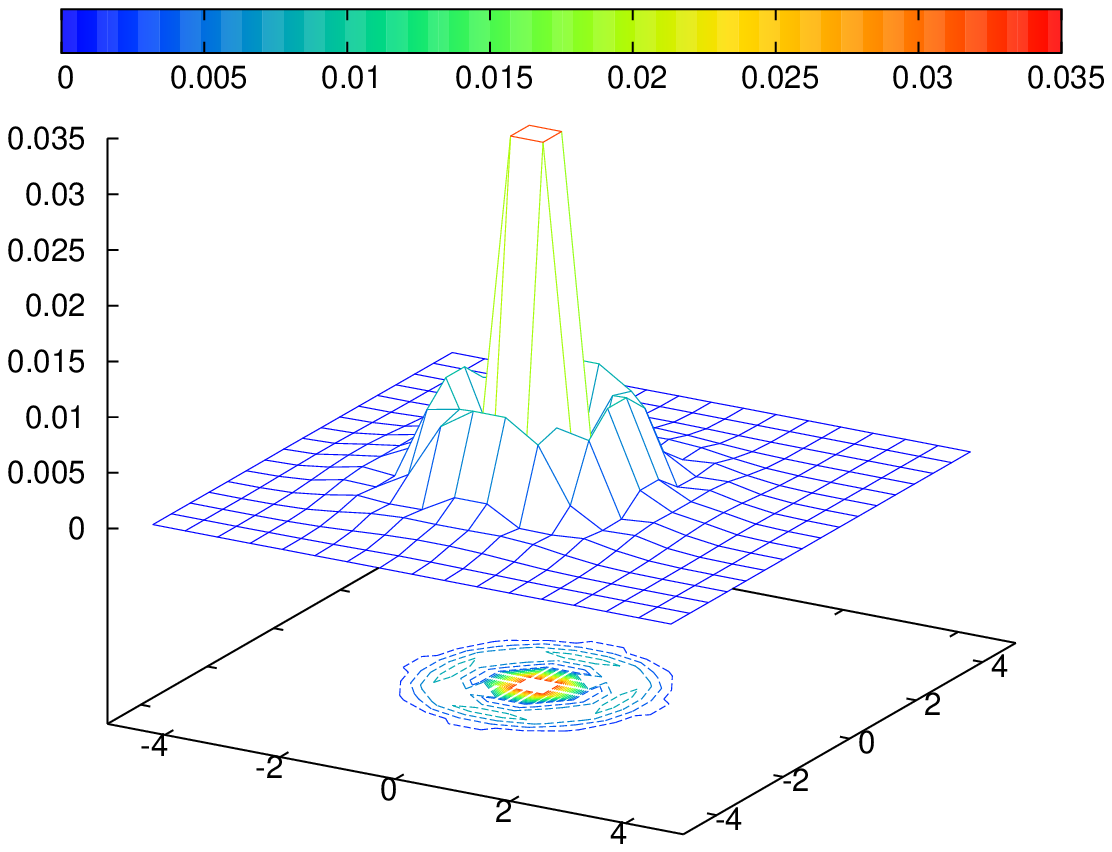}
    \includegraphics[width=0.49\textwidth]{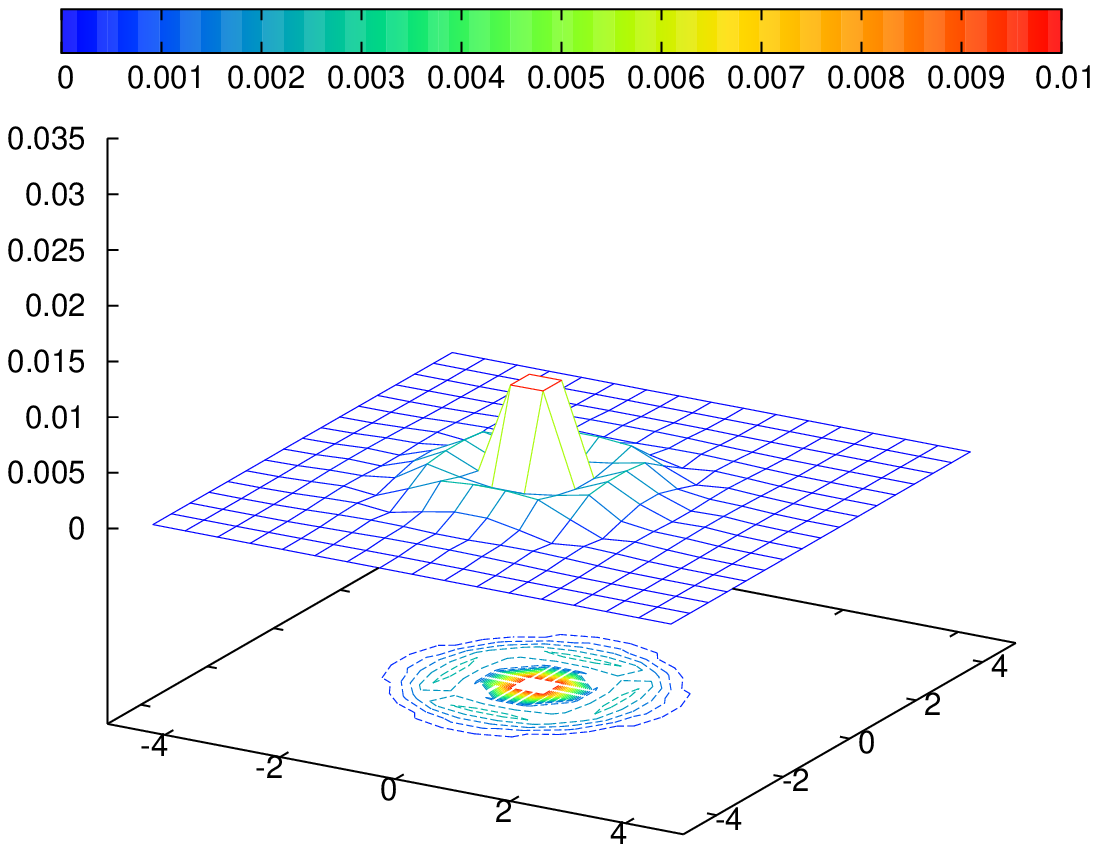}
\end{center}
  \caption{Test 1.2. Differences in the distribution functions obtained by using the BGK and the Boltzmann models 
    at  $t=0.2$, $t=1$, $t=2$ and $t=6.4$ from top-left to bottom-right. Maxwellian molecules in two dimensions.}
  \label{fig:test1_comp2}
\end{figure}

In Figure~\ref{fig:test1_comp}, we compare the convergence to equilibrium for the BGK and Boltzmann models
on the test case described in the previous paragraph using $N=32^2$ points on a domain $[-9,9]^2$. 
The $L_1$ norm of the difference between the two distribution functions $f_{\text{BGK}}(v,t)$ and $f_{\text{Boltz}}(v,t)$ as
a function of time is shown. As expected the differences are increasing at early stages of the relaxation towards the equilibrium to reach a maximum difference around
time $t\simeq 1.66$. Then, the two models slowly converge towards the same equilibrium solution. 
The red marked iterations in Figure~\ref{fig:test1_comp} are further depicted in Figure~\ref{fig:test1_comp2},
where we show the details of the difference in time between the two distribution functions.
The same azimuthal scale is employed to ease the comparison.

\subsubsection{Test 1.3. Convergence to equilibrium  in dimension three. Hard sphere molecules.}
\begin{figure}
\begin{center}
    \includegraphics[width=0.4\textwidth]{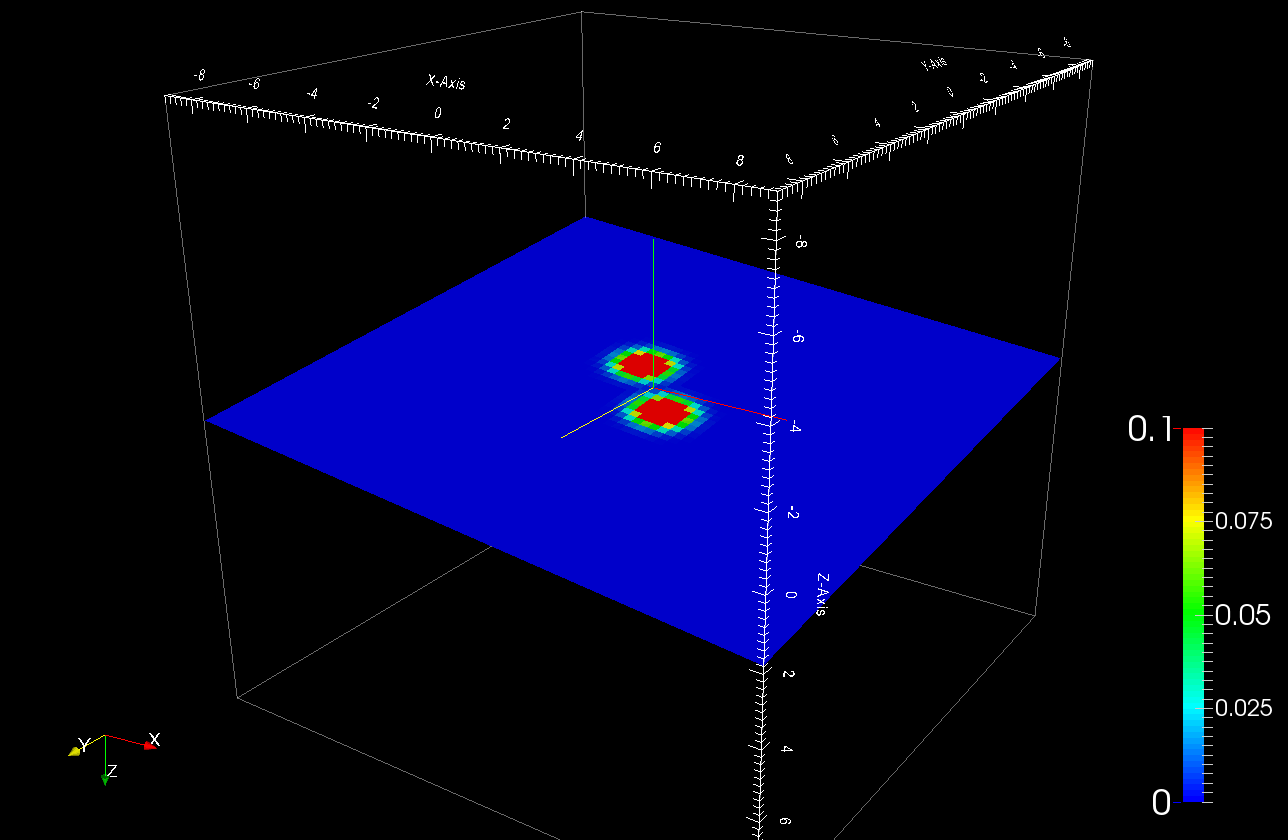} 
    \includegraphics[width=0.4\textwidth]{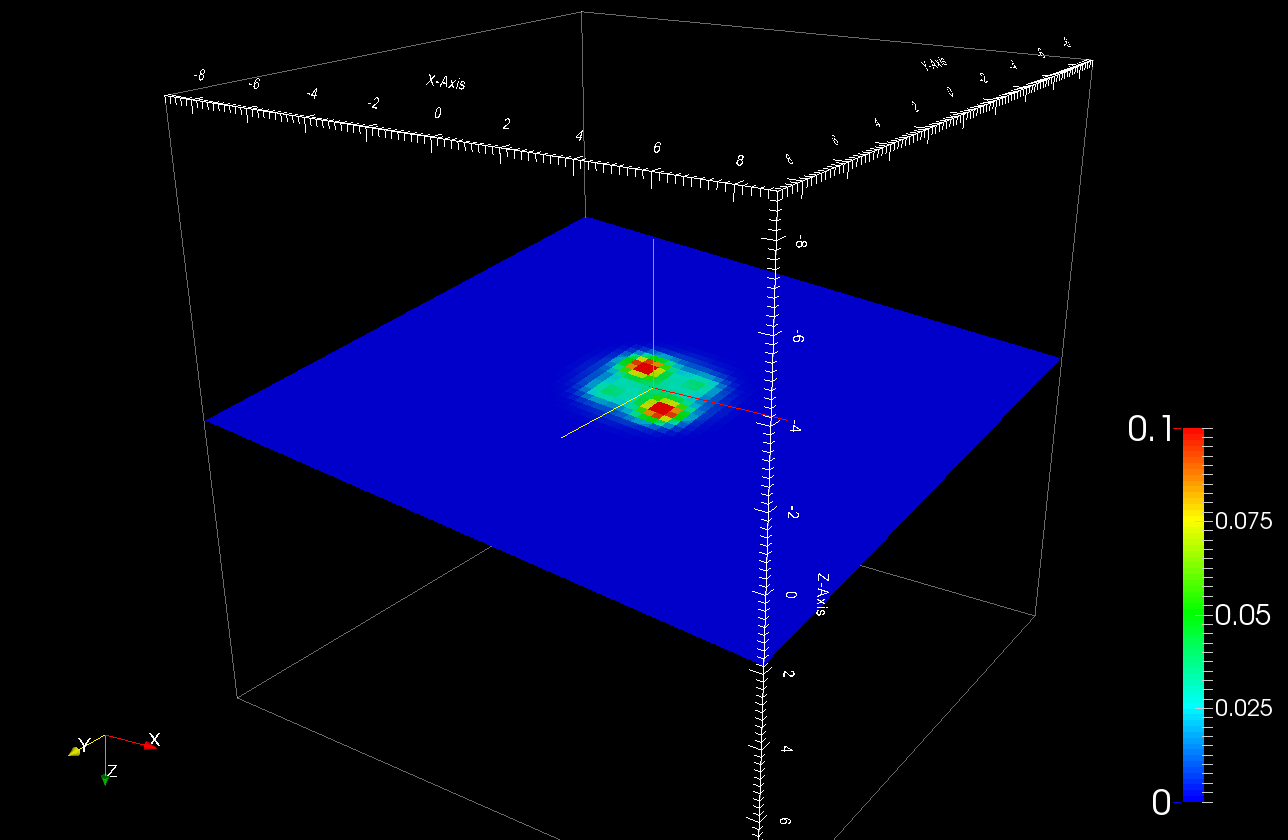} \\
    \includegraphics[width=0.4\textwidth]{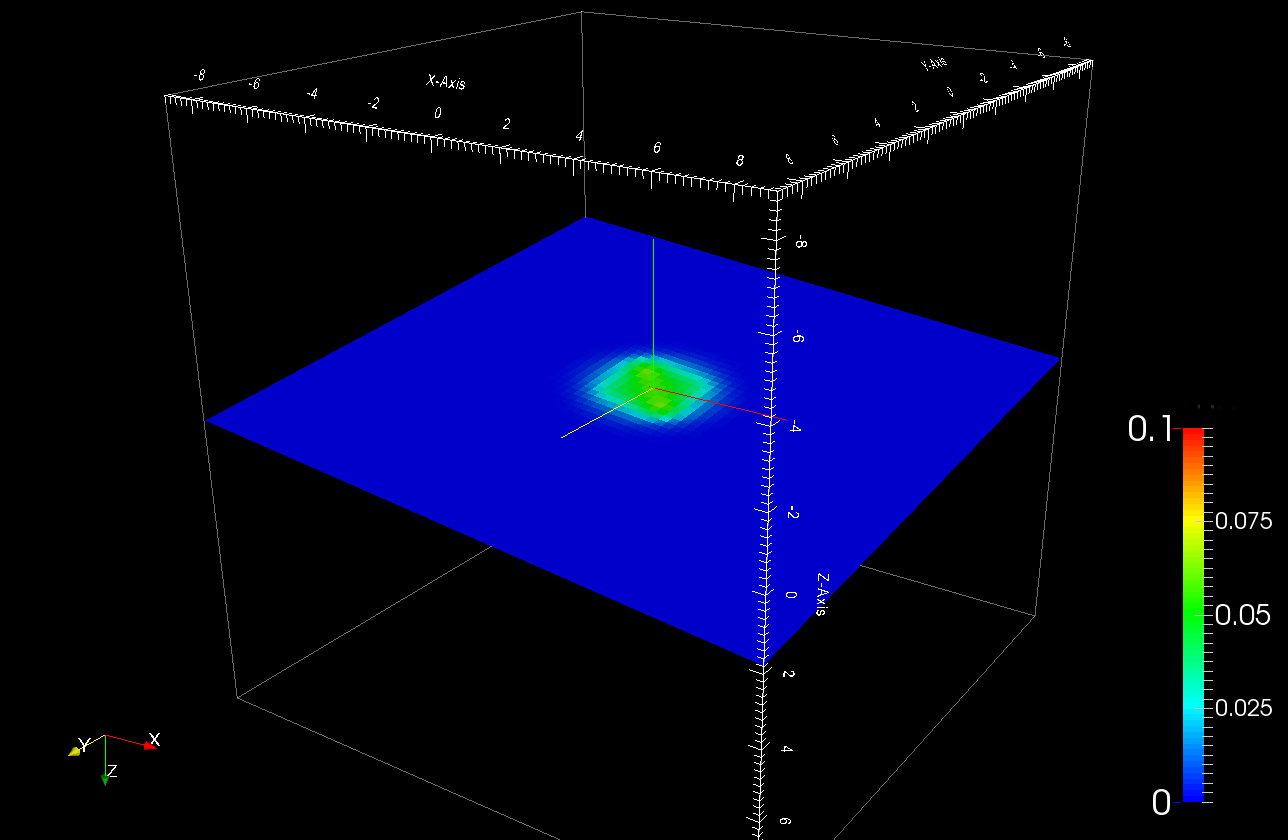} 
    \includegraphics[width=0.4\textwidth]{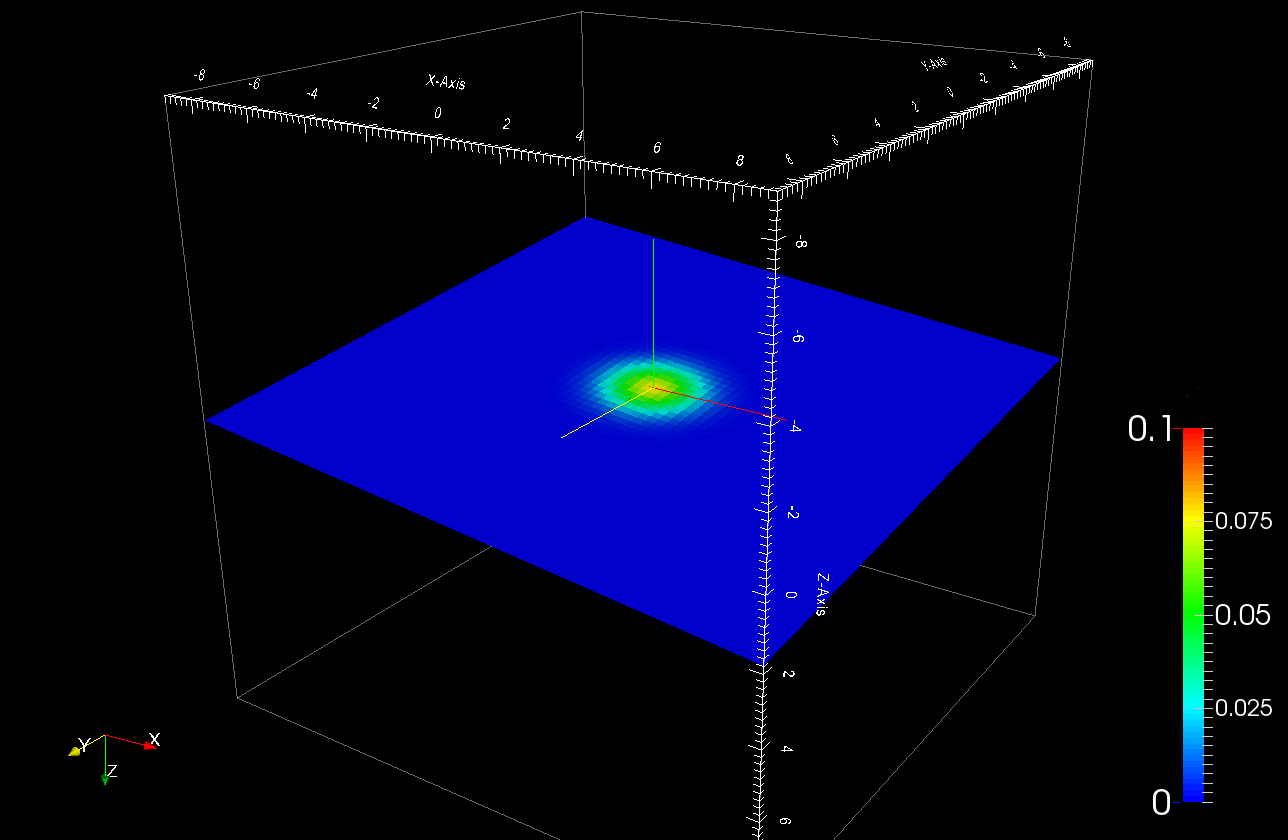} 
    \end{center}    
  \caption{Test 1.3.  Sanity check for the hard sphere case 
    at time $t=0, 0.25, 0.5$ and $t_{\text{final}} =2$ for
    $32$ points in each direction.
    The Figure shows the relaxation to the Maxwellian state for the distribution function $f(v_x,v_y,v_z,t)$ on a plane which passes from the points $\pm(−1,−1,−0.25)$ and of normal $(−0.1,−0.1,1)$.
  }
  \label{fig:test2_bis}
\end{figure}

\begin{figure}
  \begin{center}
    \includegraphics[width=0.6\textwidth]{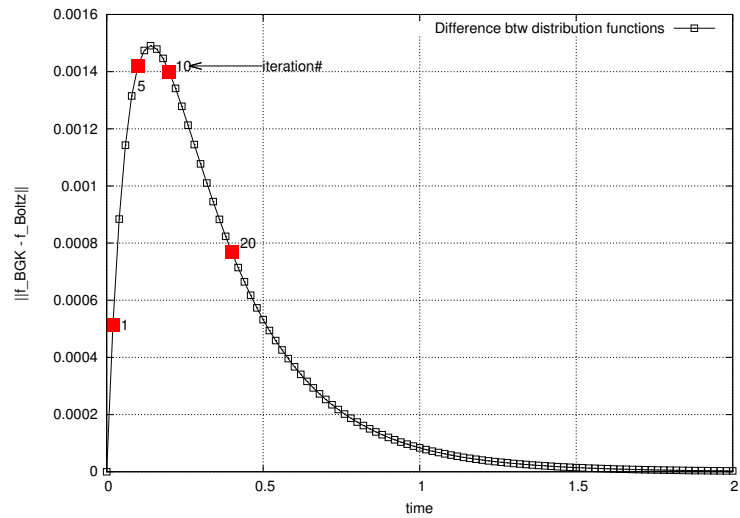}
  \end{center}
  \caption{Test 1.3. Differences between the distribution functions obtained with the BGK and the Boltzmann models. Hard sphere molecules in three dimensions.
    The distribution function values for the red marked iterations are depicted in Figure~\ref{fig:test2_comp2}.
  }
  \label{fig:test2_comp}
\end{figure}

\begin{figure}
  \begin{center}
    \includegraphics[width=0.45\textwidth]{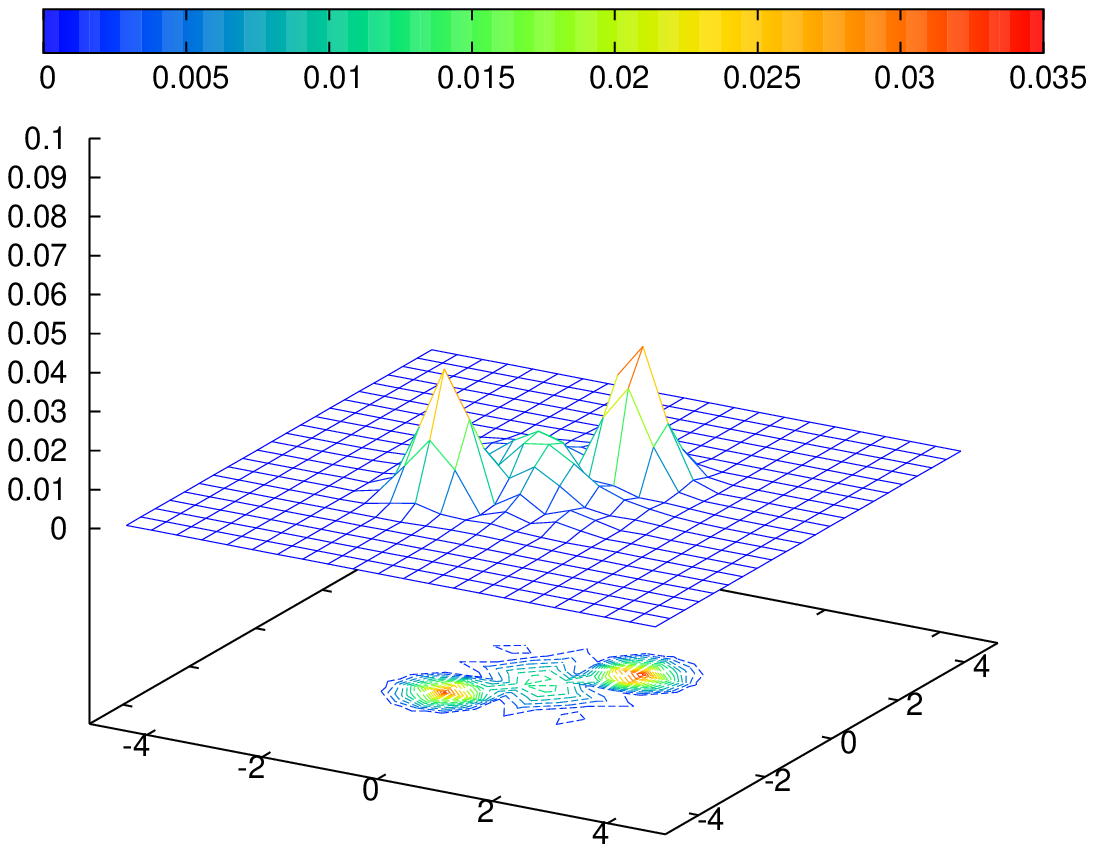}
    \includegraphics[width=0.45\textwidth]{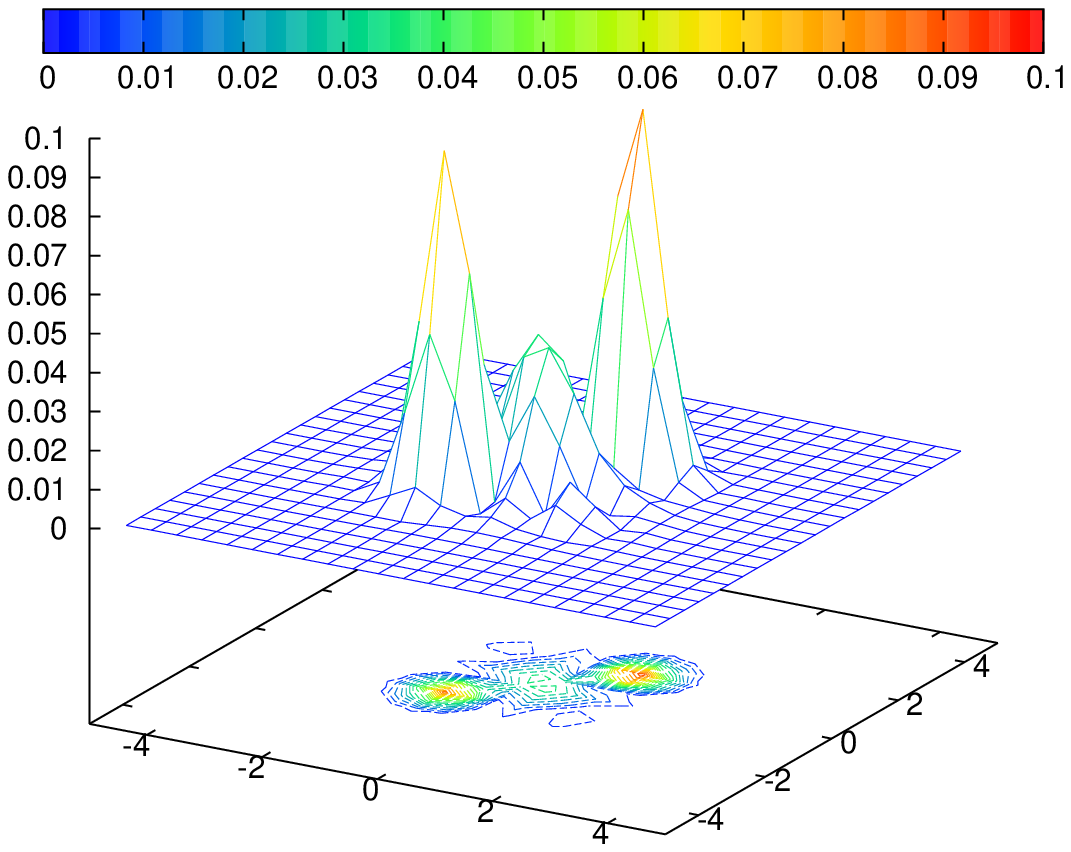}\\
    \includegraphics[width=0.45\textwidth]{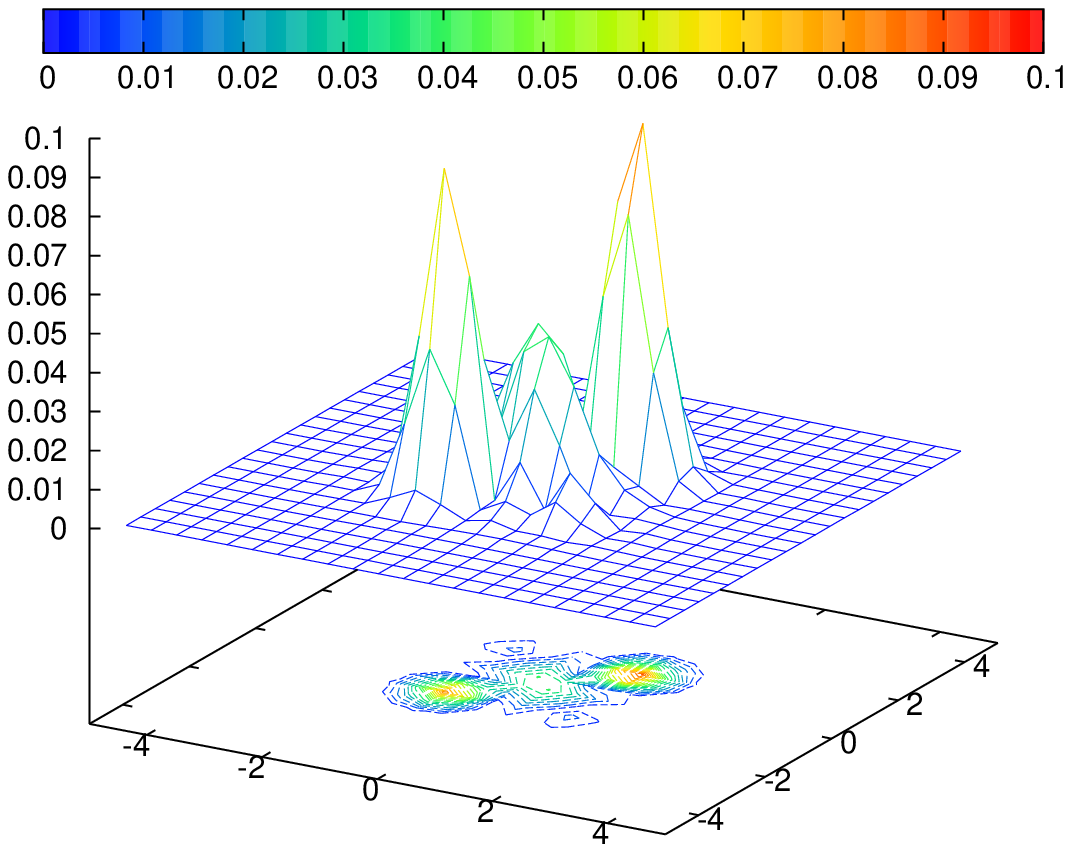}
    \includegraphics[width=0.45\textwidth]{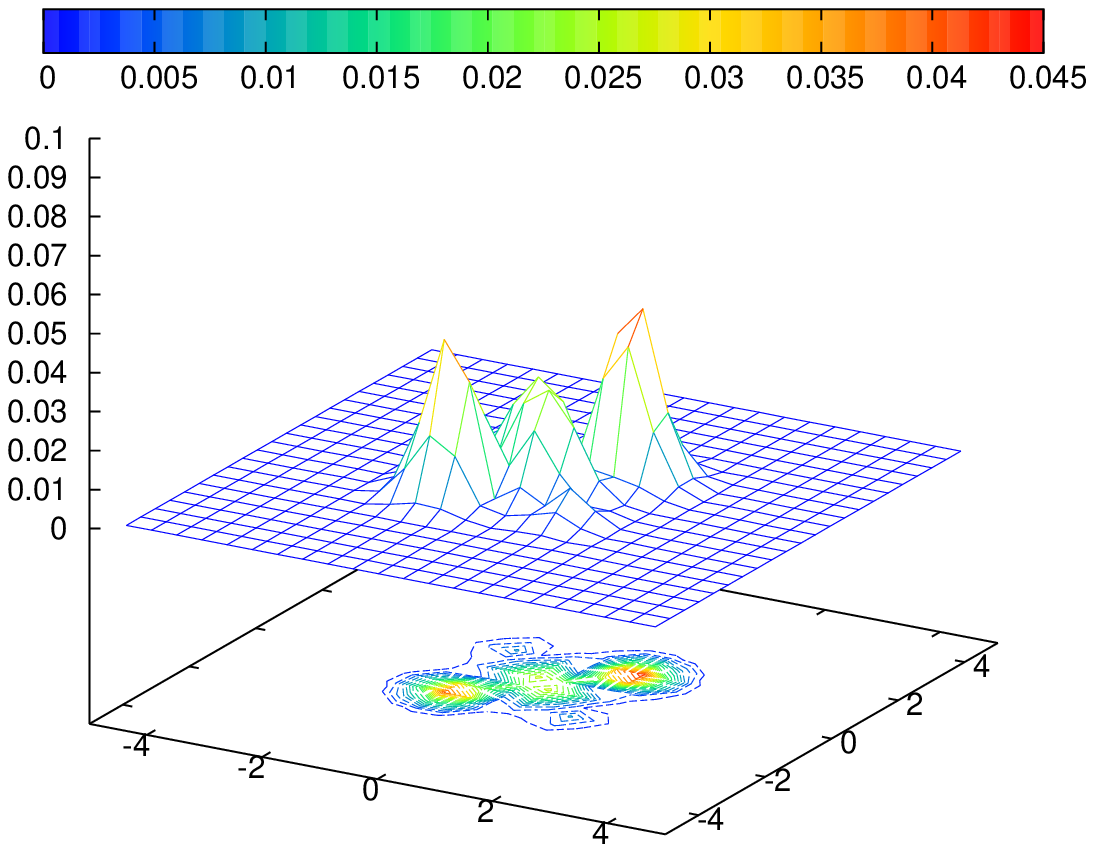}
  \end{center}
  \caption{Test 1.3. Differences in the distribution functions obtained by using the BGK and the Boltzmann models 
    at  $t=0.05$, $t=0.25$, $t=0.5$ and $t=1$ from top-left to bottom-right for $v_z=0$. Hard sphere molecules in three dimensions.}
  \label{fig:test2_comp2}
\end{figure}

The following initial condition is considered \cite{FilbetRusso} 
\begin{gather}
  f(v,t=0) = \frac{1}{2(2\pi\sigma^2)^{3/2}}
 \left[ \exp \left(-\frac{|v-v_1 |^2}{2\sigma^2} \right) + \exp \left(-\frac{|v+v_1 |^2}{2\sigma^2} 
\right) \right],
  \label{eq:Jmgb}
\end{gather}
where $\sigma^2=0.2$ and $v_1$ is $v_1=(v_x,v_y,v_z)=(-1,-1,-0.25)$. The final time is set to $t_{\text{final}} = 2$, the time step is $\Delta t = 0.05$,
while the velocity domain is $[-7;7]^3$ discretized with $N=32^3$ points.
\begin{figure}
\begin{center}
  \hspace{-1cm}
    \includegraphics[width=0.3825\textwidth,height=0.29\textwidth]{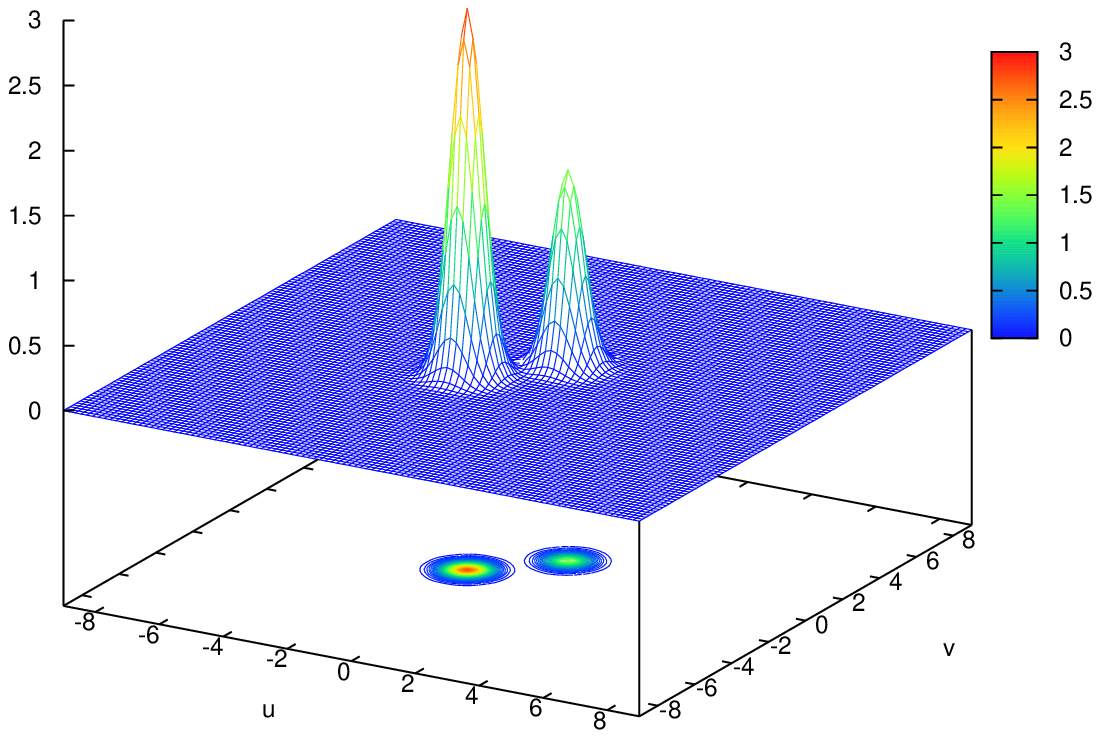}
  \hspace{-1cm} 
    \includegraphics[width=0.3825\textwidth,height=0.29\textwidth]{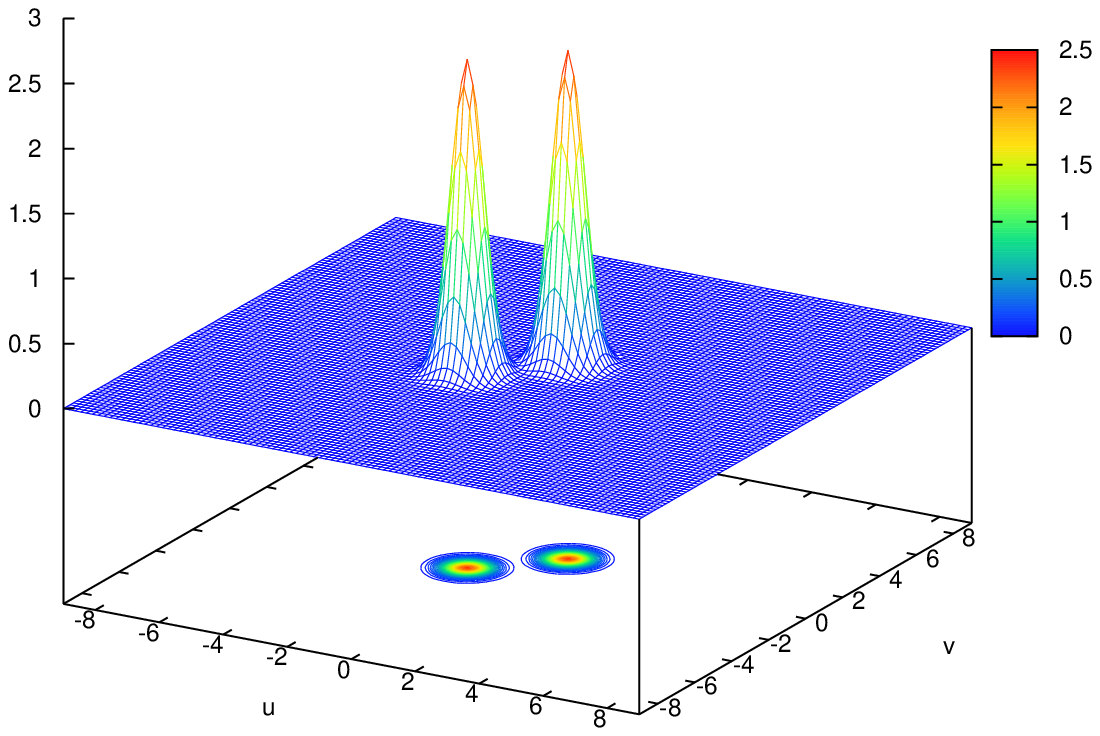} 
  \hspace{-1cm}
    \includegraphics[width=0.3825\textwidth,height=0.29\textwidth]{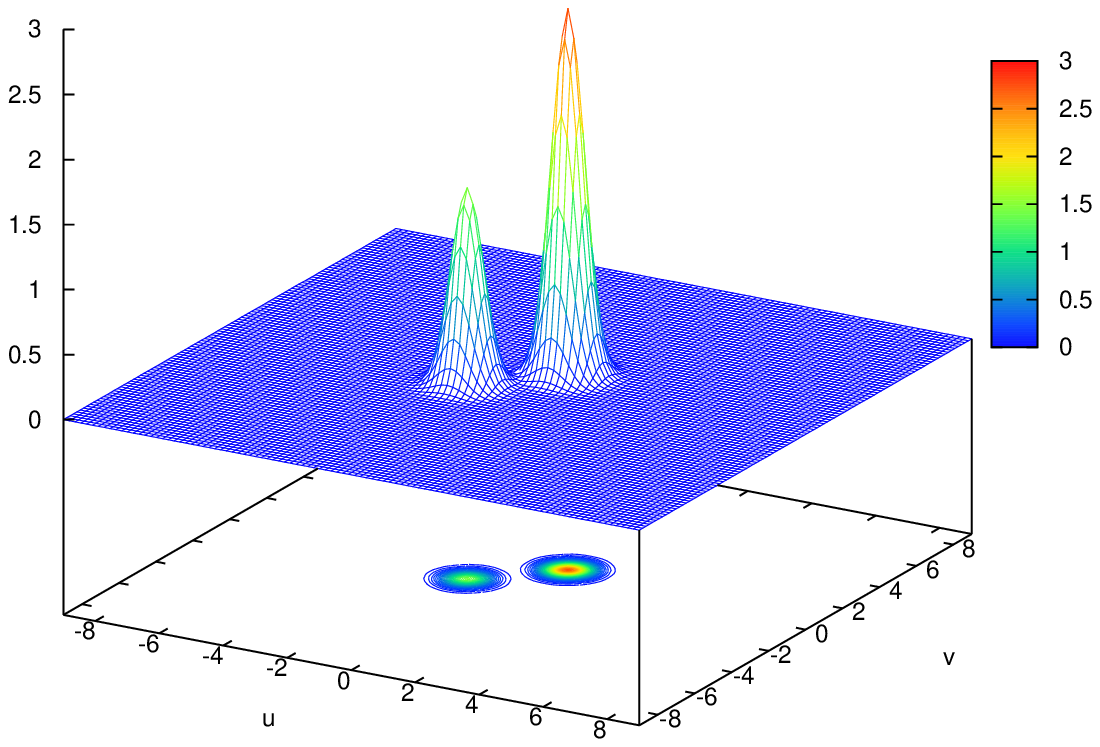} 
    \end{center}    
  \caption{Test 1.3. Sanity check --- 
    Initial distribution function for $f(v_x,v_y,0)$ middle, $f(v_x,v_y,0.25)$ left  and $f(v_x,v_y,-0.25)$ right. Hard sphere molecules.}
  \label{fig:test2_t0}
\end{figure}
In Figure~\ref{fig:test2_t0}, we report the initial data while in Figure \ref{fig:test2_bis}, we report the relaxation to the equilibrium for the hard sphere model on a fixed plane
passing through the points of coordinates $\pm(−1,−1,−0.25)$ and of normal $(−0.1,−0.1,1)$\footnote{This plane passes through the center of the two spheres
defined by (\ref{eq:Jmgb})}. In Figure \ref{fig:test2_comp}, we report as before the comparison between the Boltzmann and the BGK model 
measuring the $L_1$ difference between the two distribution functions in time. Finally, in Figure~\ref{fig:test2_comp2} are presented the details of such differences between the distribution functions 
for the red marked iterations of figure~\ref{fig:test2_comp} at $v_z=0$. As for the Maxwellian molecules in two dimensions, at the beginning of the simulation we observe the larger deviations between the two models while at the end they both converge to the  same limit solution as expected. \\
The results of these homogeneous tests (in two and in three dimensions in velocity) are meant to verify and validate our implementation of the spectral discretization of the Boltzmann collisional operator. Moreover, the differences observed between the two models justify the necessity of the Boltzmann operator for space non homogeneous situations. We can now investigate the coupling of the Boltzmann operator with the Fast semi-Lagrangian method for several space non homogeneous cases.

%
%

%
%
\subsection{Part 2. Numerical results for the one dimensional case in space.} \label{sec:act1}
In this part, we focus on solving the one dimensional in space Boltzmann and BGK equations considering a two or three dimensional
dimensional velocity space setting. The purposes are twofold. First, we numerically demonstrate that the FKS and the spectral accurate Boltzmann kernel solver may be 
appropriately coupled and that they provide a valid kinetic solver for Boltzmann equations. Second, we show that BGK and Boltzmann models provide different results,
justifying the use of a more sophisticated model. For these tests, a classical Riemann problem with Sod like initial data is considered
\bea
\nonumber
 \rho=1,    \quad u=0, \quad T=2.5, & \;\;  \; \text{if} \;\;& x\leq L/2,  \\
\nonumber
 \rho=0.125, \quad u=0,  \quad T=0.25,& \;\;  \; \text{if} \;\;& x> L/2,  
\eea
with $\Omega=[0,2]$. 
Initial local thermodynamics equilibrium is considered for all tests: $f(x,v,t=0)=M[f](x,v,t=0)$. The velocity space is set to $L_v=[-15,15]$ for all cases.
Dirichlet boundary conditions are set on the left/right boundaries of $\Omega$. The two kinetic models are solved by employing a time rescaling factor in order to put in evidence the role
of the collisions in the solutions. The rescaled equations reads\be \partial_t f+v\cdot\nabla_x f=\frac{1}{\tau}Q(f),\ee
where $\tau$ is the rescaled parameter (the frequency of relaxation) which plays the role of the non dimensional Knudsen number. Smaller is the relaxation frequency, faster is the relaxation
of the distribution function towards the equilibrium state. However, the exact rate of convergence is dependent on the type of collision considered either BGK or Boltzmann (Maxwellian molecules or hard spheres).

\subsubsection{Test 2.1. Numerical convergence of the Boltzmann equation. The two dimensional in velocity Maxwellian molecules case.}
\begin{figure}[ht]
  \begin{center}
    \includegraphics[width=0.35\textwidth]{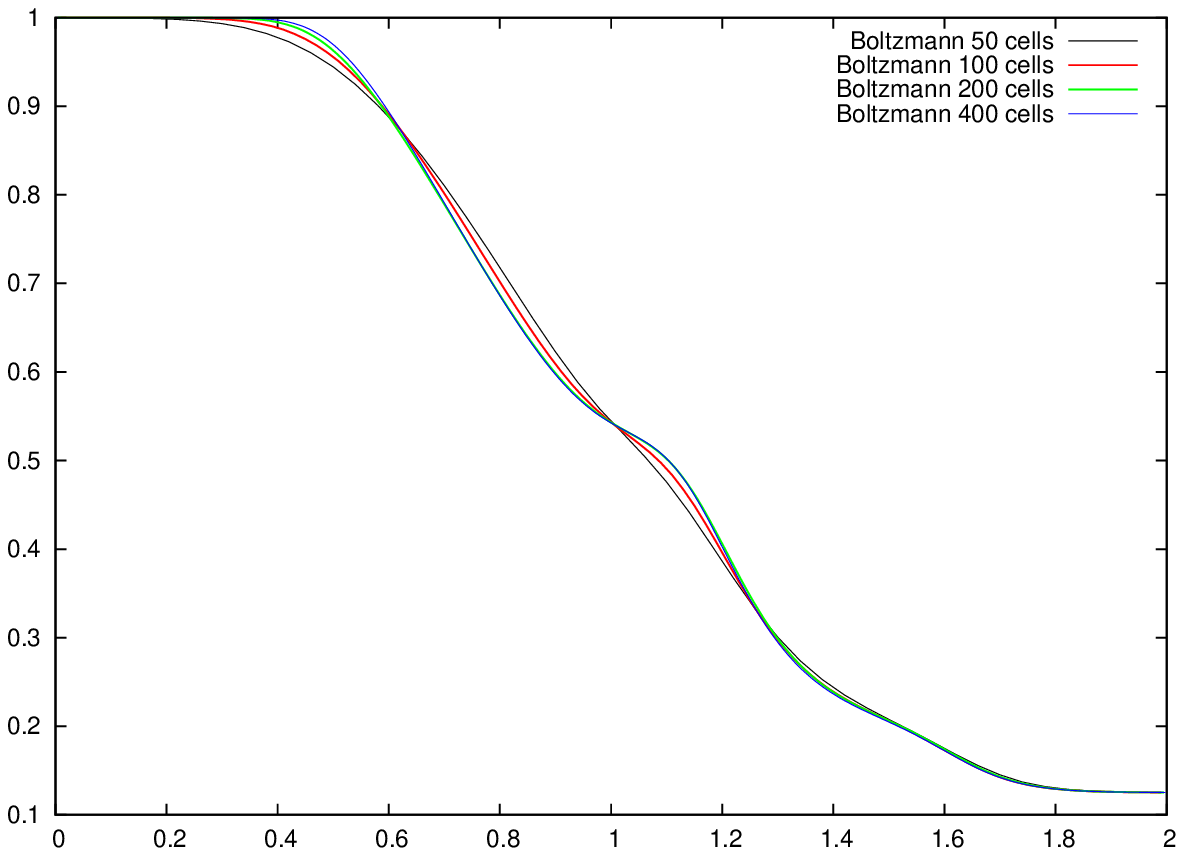} 
    \includegraphics[width=0.35\textwidth]{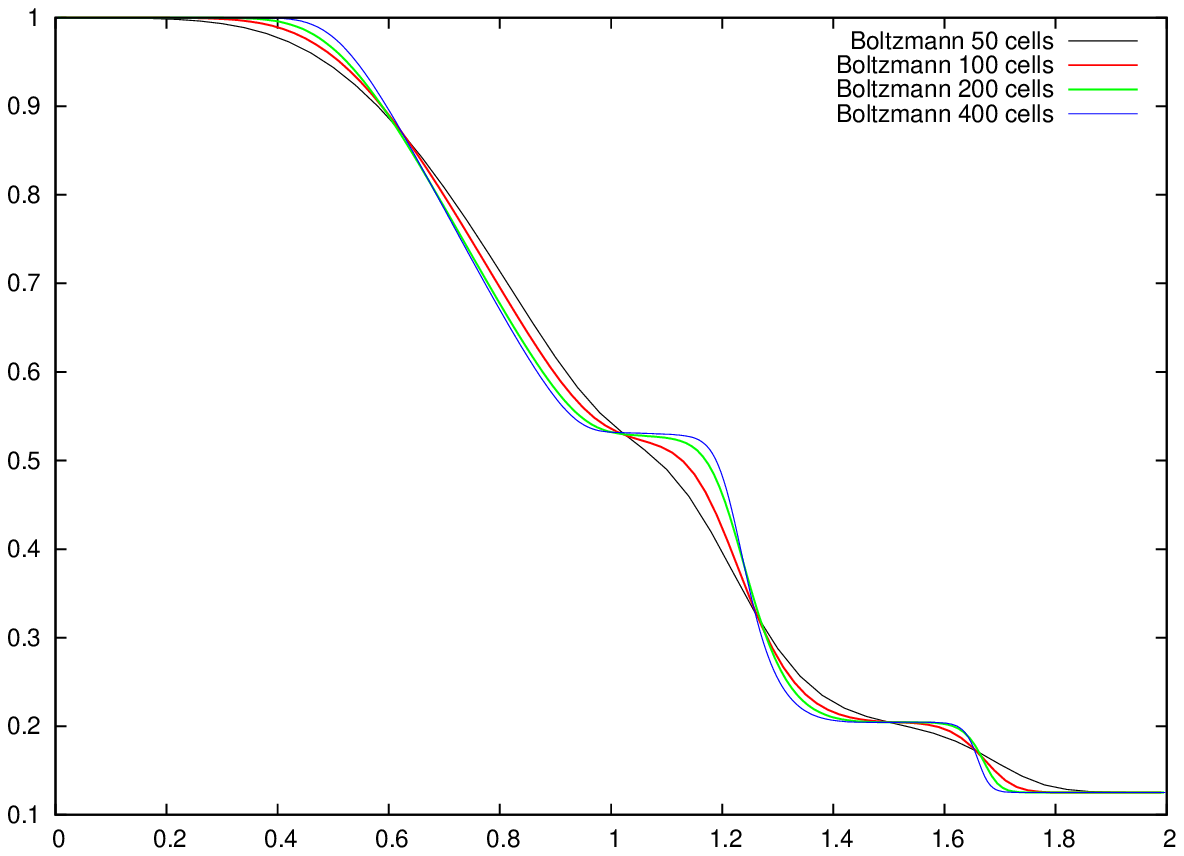}\\ 
    \includegraphics[width=0.35\textwidth]{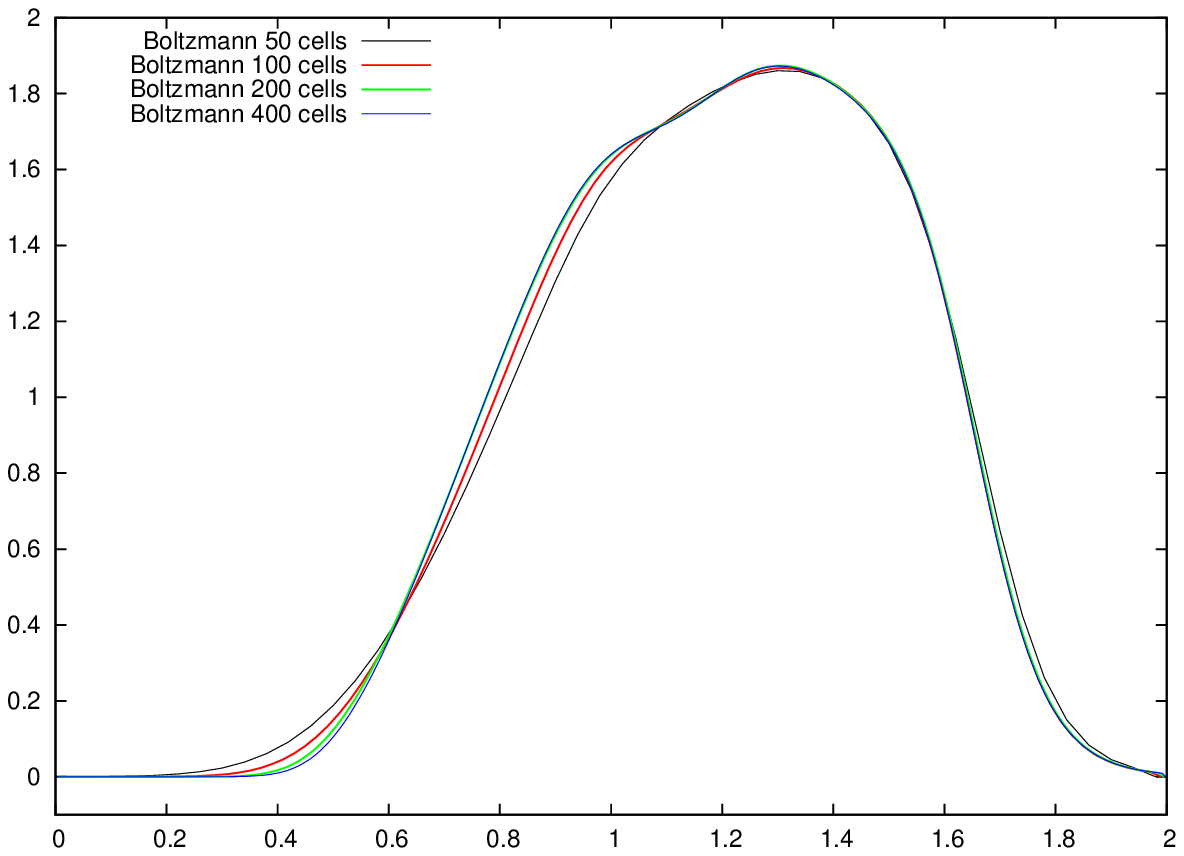} 
    \includegraphics[width=0.35\textwidth]{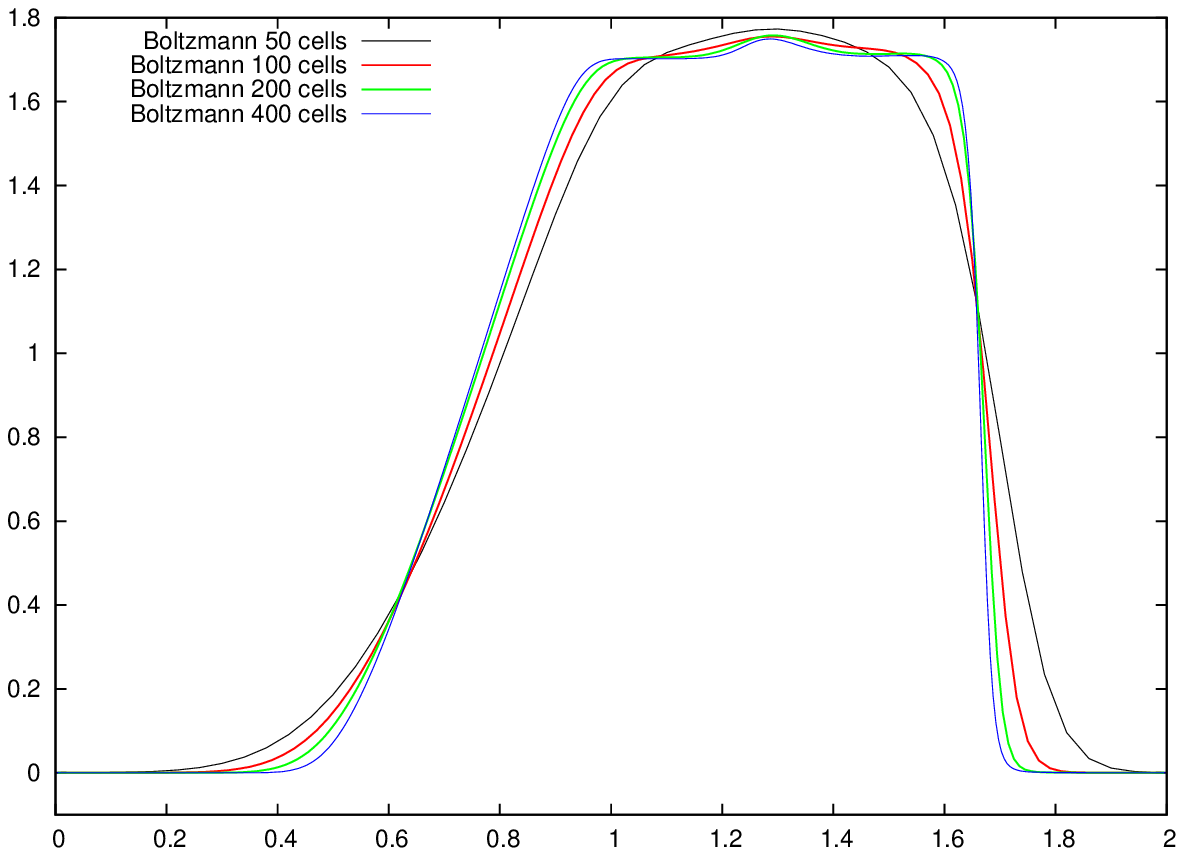}\\
    \includegraphics[width=0.35\textwidth]{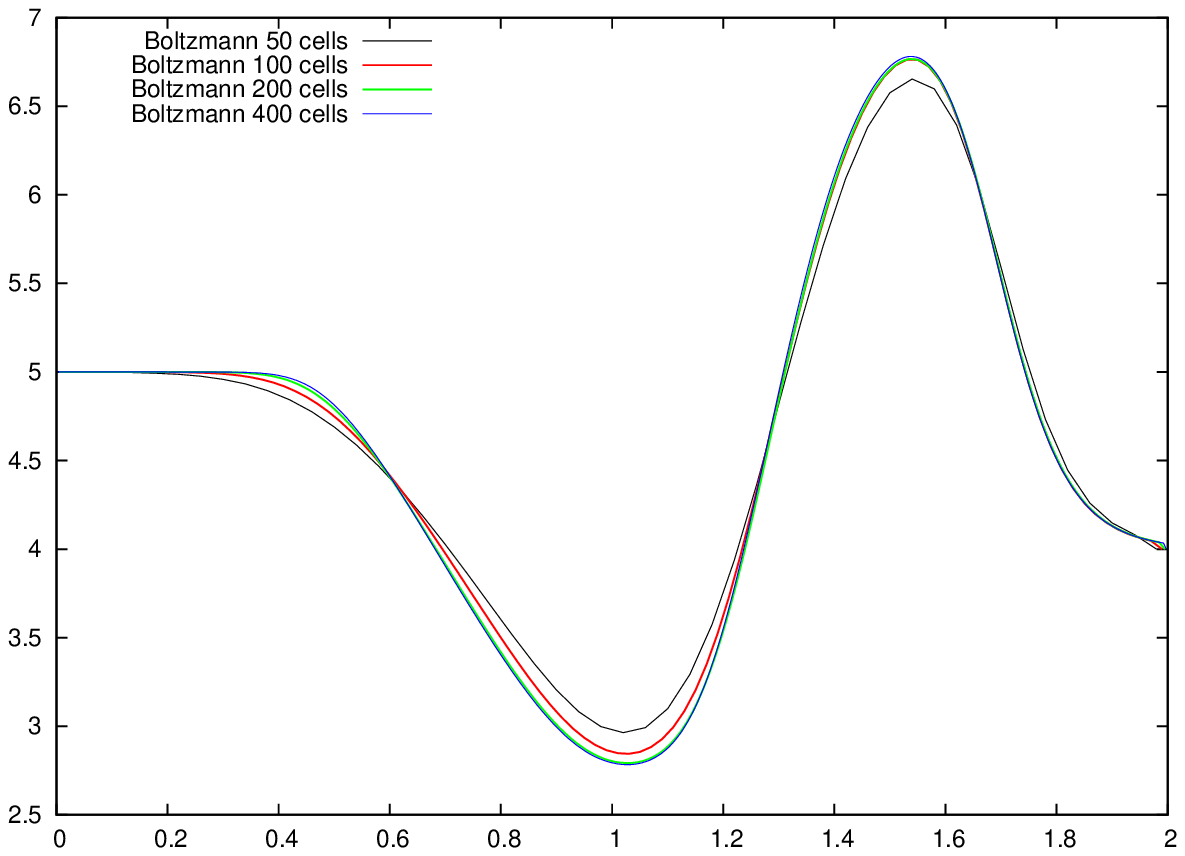} 
    \includegraphics[width=0.35\textwidth]{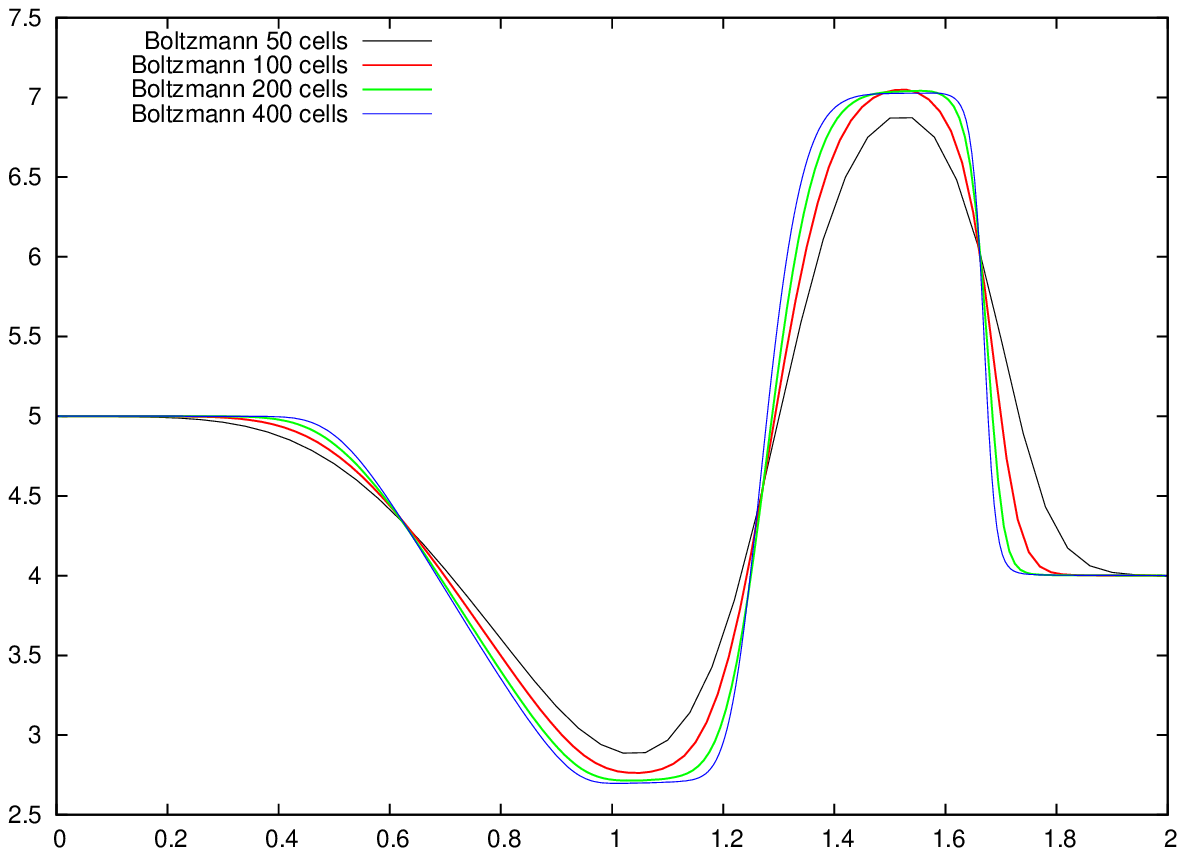} 
  \end{center}
  \caption{Test 2.1. One dimension in space and two dimension in velocity Boltzmann model with Maxwellian molecules for a Sod like test case at $t_{\text{final}} =0.15$.
    Mesh convergence results for $\tau=10^{-3}$ (left) and $\tau=10^{-4}$ (right).
    Density (top), velocity (middle) and temperature (bottom) are shown for $M=50, 100, 200$, and $400$ cells and $N=64^2$ velocity cells.
   }
  \label{fig:test3_0}
\end{figure}

The two dimensional velocity space is first considered leading to a space/velocity mesh of the form $M\times N$
with $N=64^2$ and a varying number of space cells $M$. In Figure~\ref{fig:test3_0}, we present the space convergence results for the density, the velocity and the temperature
when the Boltzmann operator is solved. Successively refined (doubled) spatial meshes are employed, from $50$ to $400$ up to final time $t_{\text{final}} =0.15$. From these data we can observe that the 
simulation results seem to converge towards a given numerical solution in both cases $\tau=10^{-3}$ (left panels) and $\tau=10^{-4}$ (right panels). We can also observe that
the increase in mesh resolution is profitable especially for smaller $\tau$. This is the same behavior observed in \cite{FKS}. In fact, the scheme precision decreases as the equilibrium
state is approached, being virtually exact in non collisional or almost non collisional regimes. The loss of precision observed in fluid dynamic regimes can be recovered with
a similar technique as the one proposed in \cite{FKS_HO}. Here, however we do not consider this possibility. The CFL condition employed is the following
\be \Delta t\leq \min \left(\frac{\Delta x}{|v_{max}|},\frac{\tau}{\rho}\right),\ee
where the first term on the right hand side of the above equation comes from the will of keeping the error small enough in the splitting scheme, while the second term is due to the stability
restriction in the solution of the space homogeneous problem when Maxwellian molecules are employed. In fact, in this case the loss part of the collision integral $Q^-(f)$ can be estimated, giving $Q^-(f)=\rho(f)f$.

\subsubsection{Test 2.2. Comparisons between the BGK model and the Boltzmann model. The two dimensional in velocity Maxwellian molecules case.}
\begin{figure}
  \begin{center}
    \includegraphics[width=0.38\textwidth]{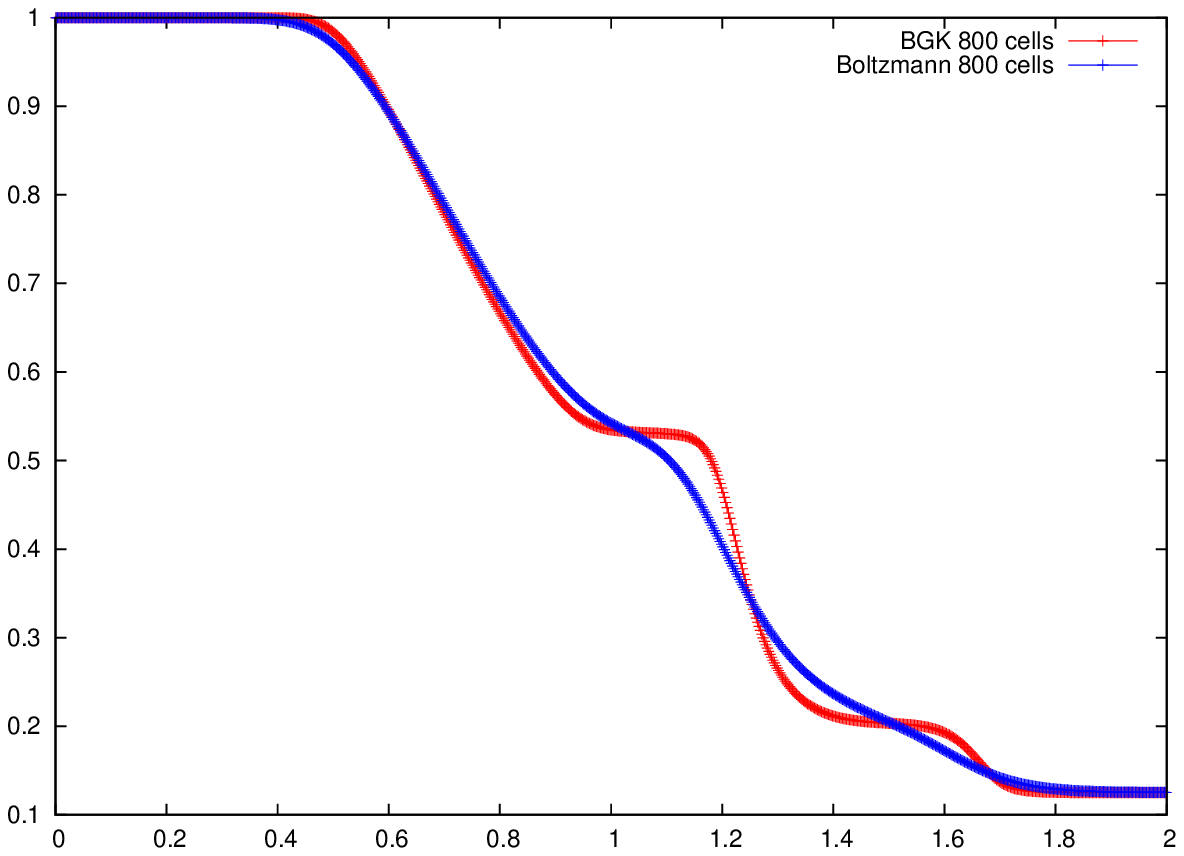} 
    \includegraphics[width=0.38\textwidth]{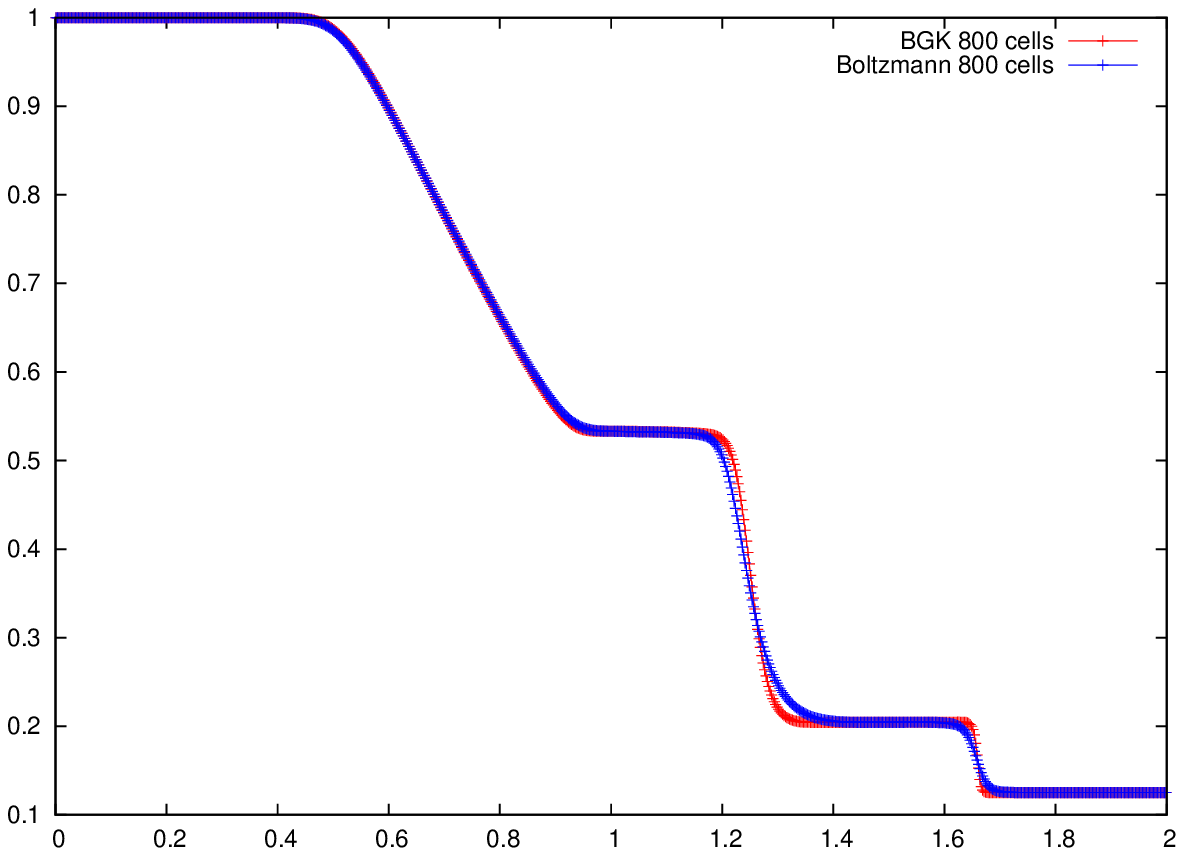}\\ 
    \includegraphics[width=0.38\textwidth]{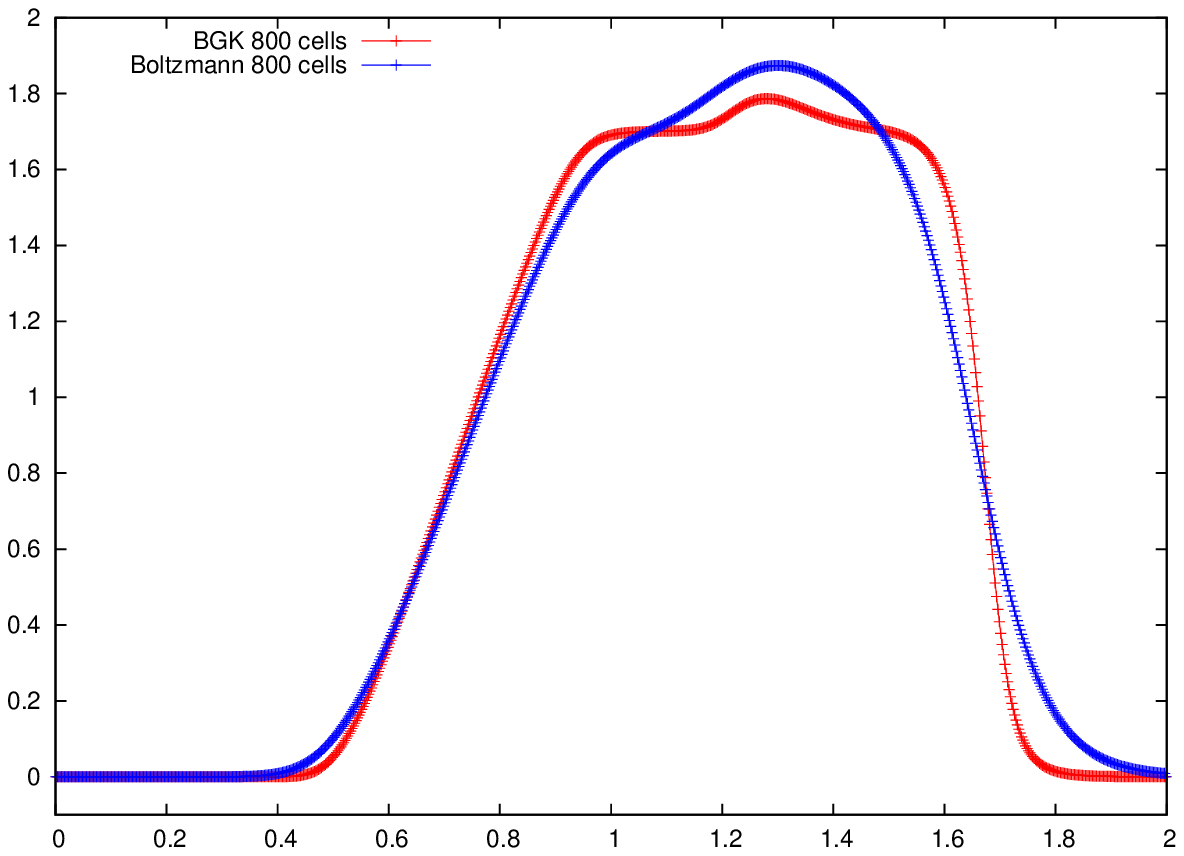} 
    \includegraphics[width=0.38\textwidth]{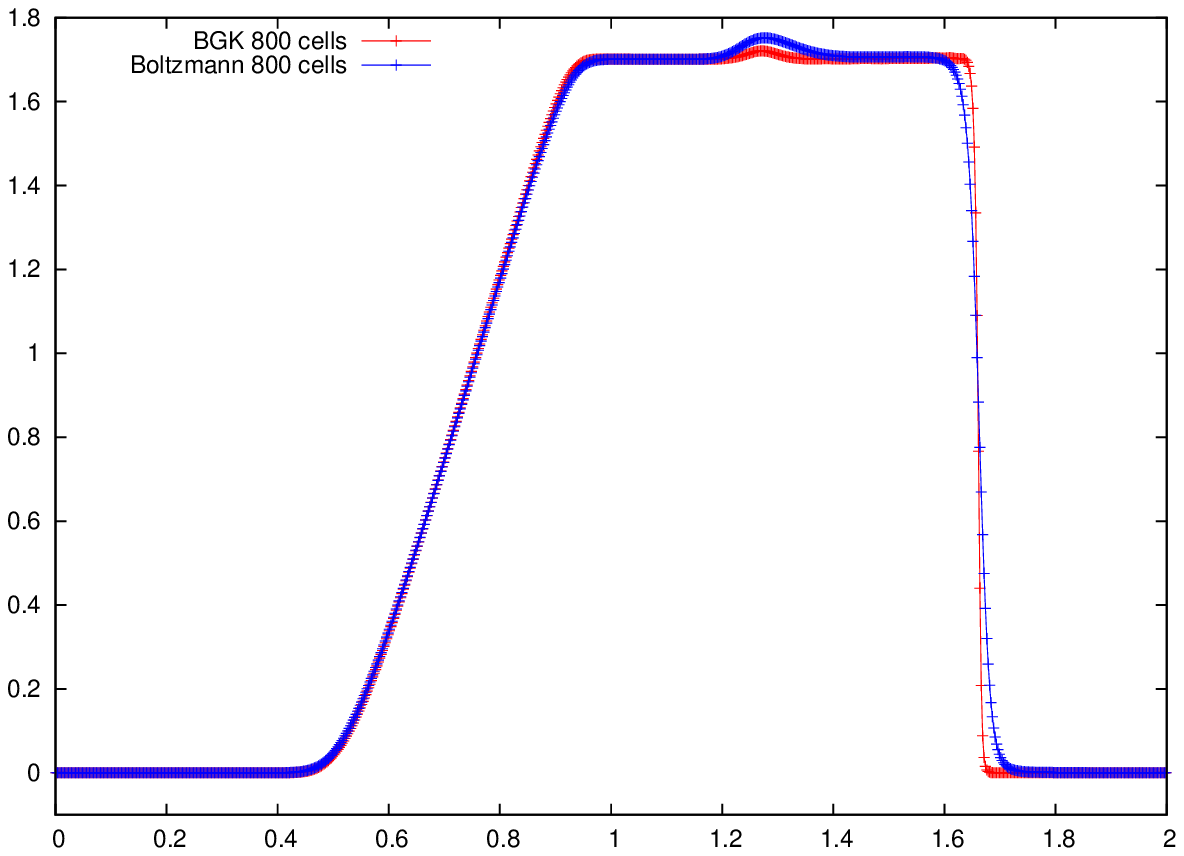}\\ 
    \includegraphics[width=0.38\textwidth]{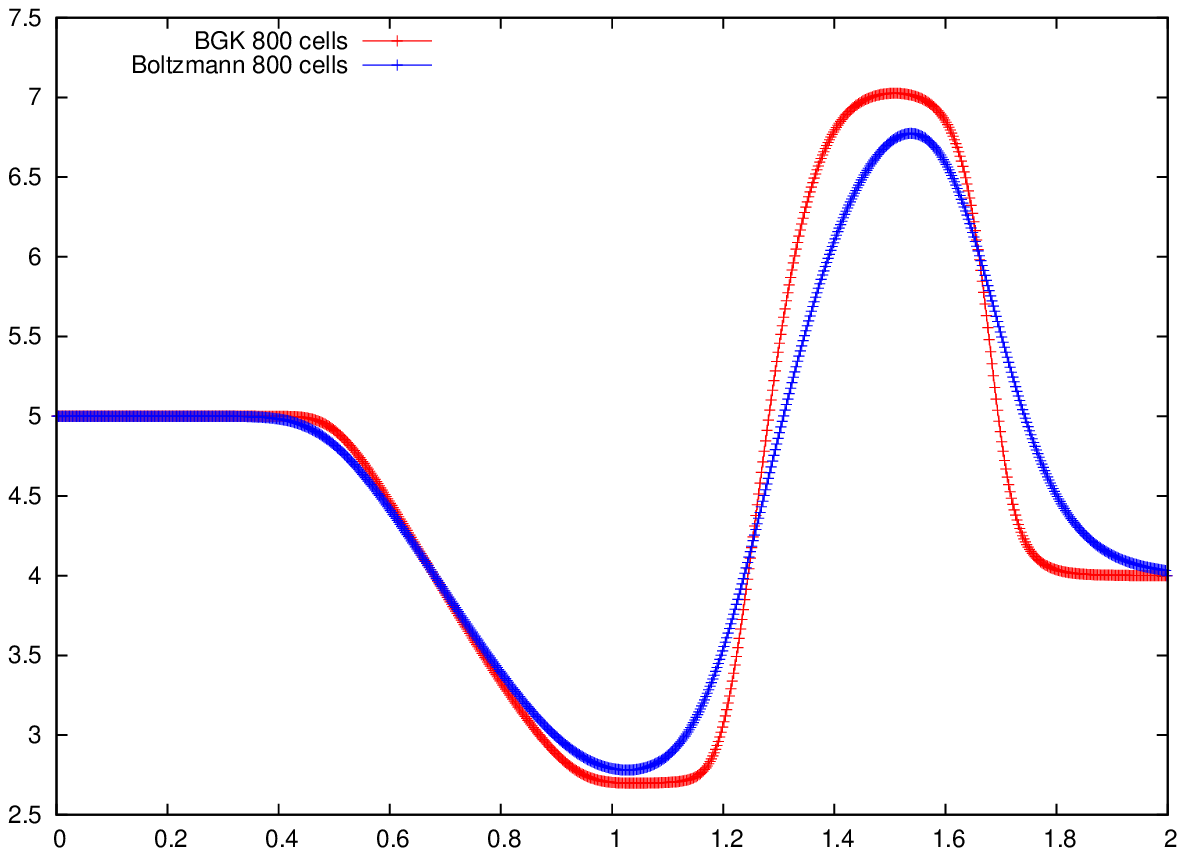} 
    \includegraphics[width=0.38\textwidth]{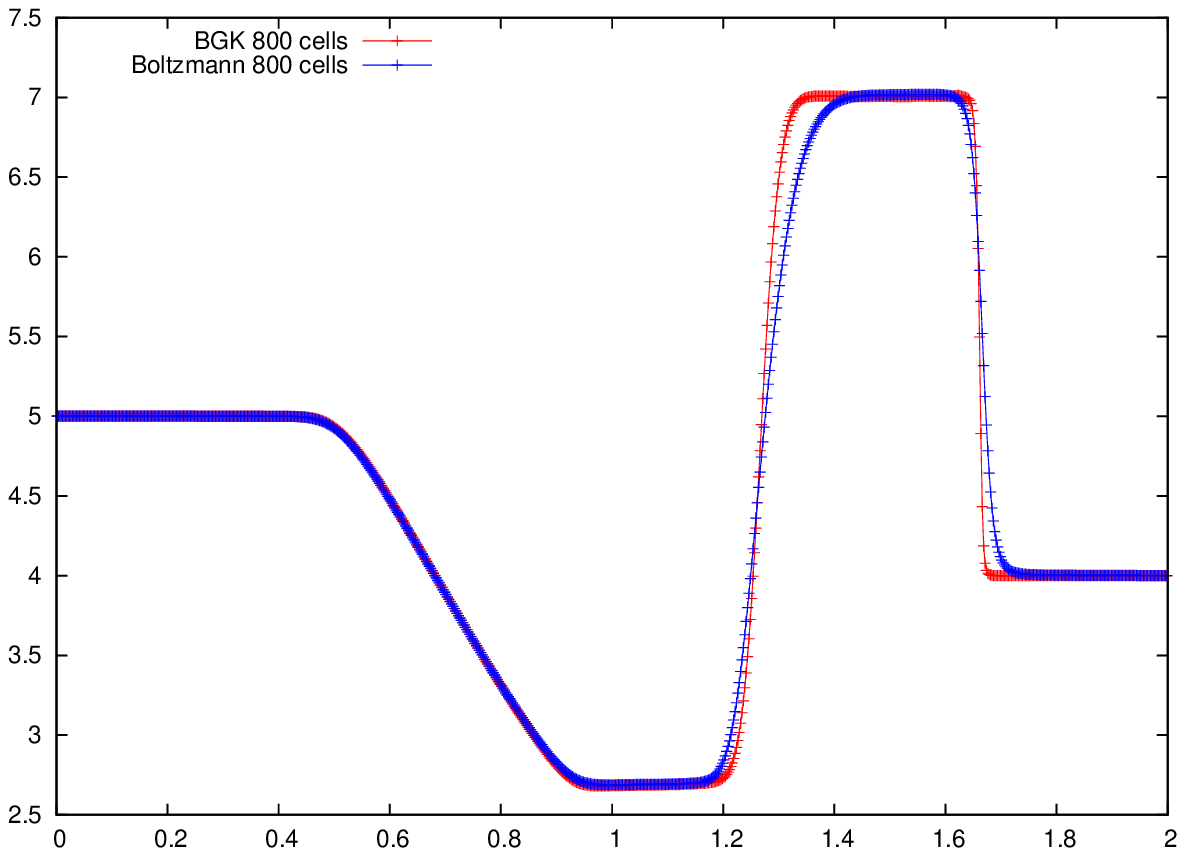} \\
    \includegraphics[width=0.38\textwidth]{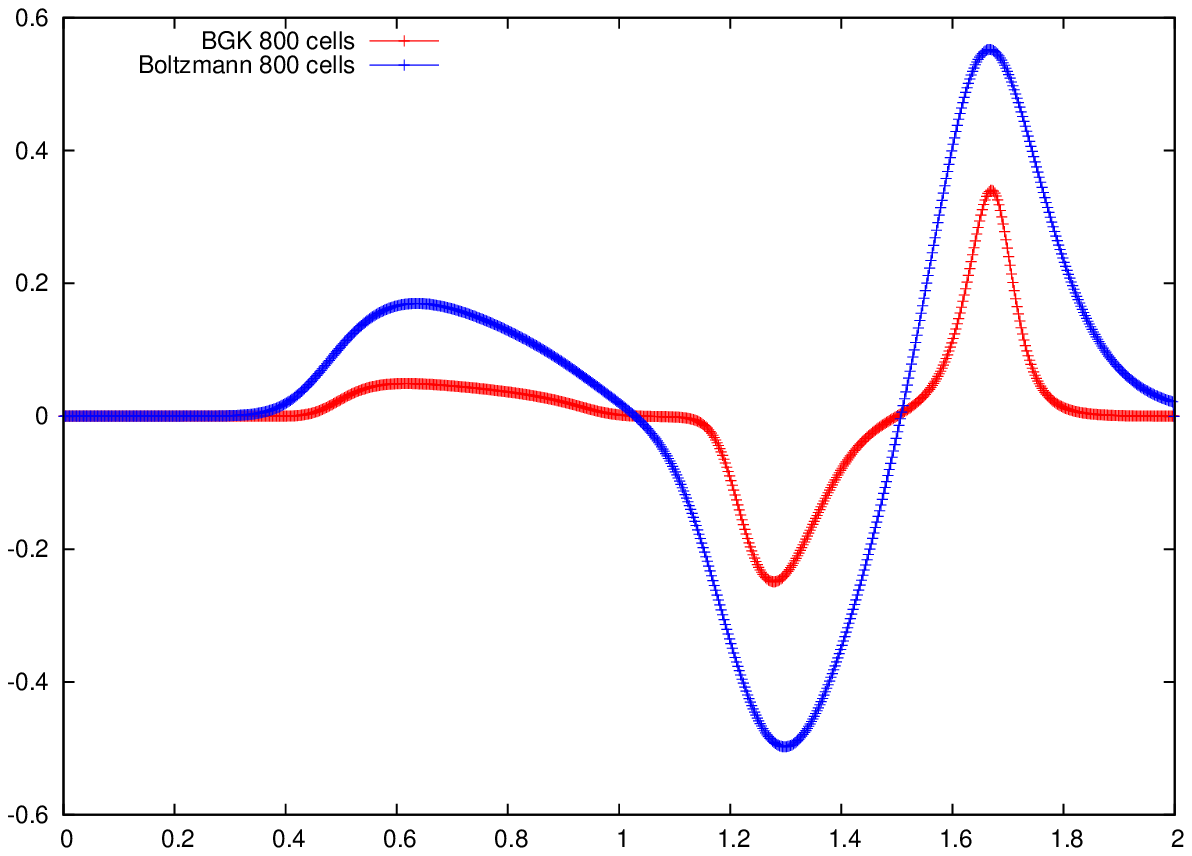} 
    \includegraphics[width=0.38\textwidth]{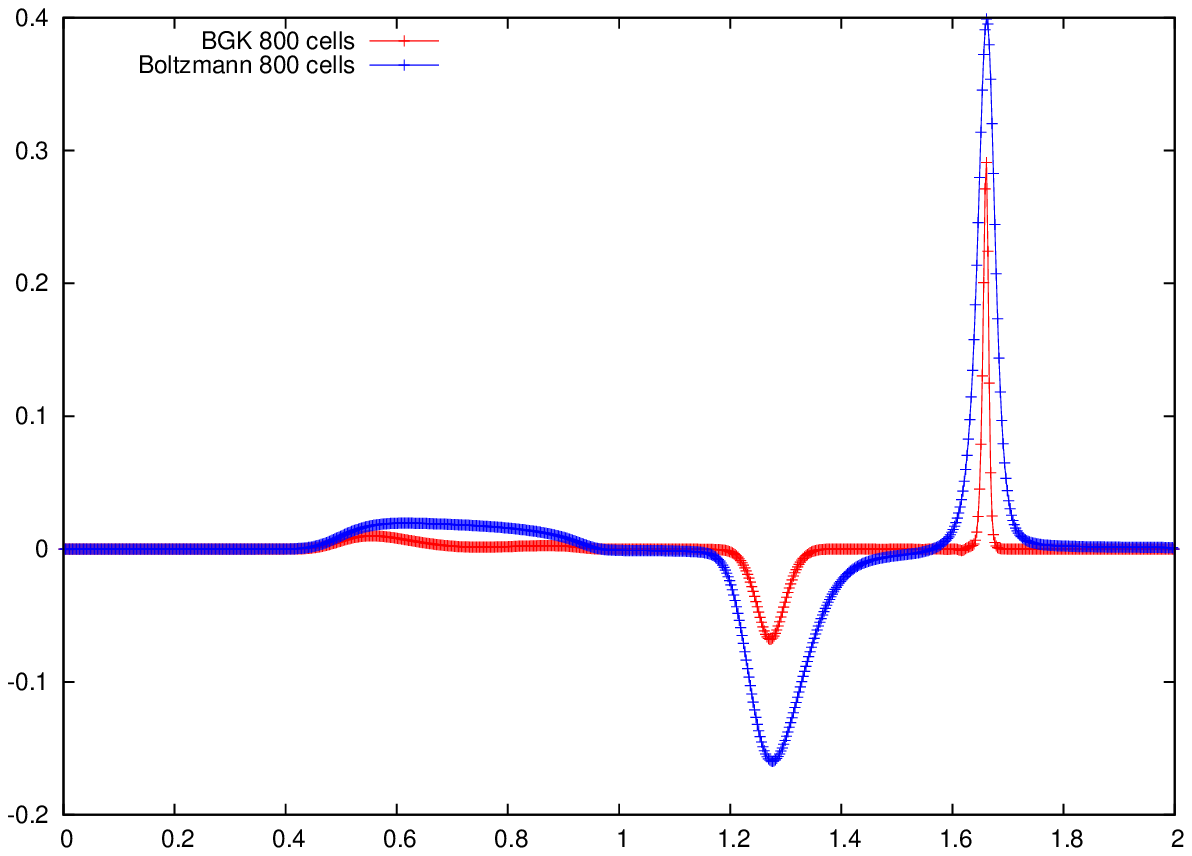} 
  \end{center}
  \caption{Test 2.2. One dimension in space and two dimension in velocity Boltzmann model with Maxwellian molecules (blue)
    and a BGK model (red) for a Sod like test case at $t_{\text{final}} =0.15$ for $\tau=10^{-3}$ (left) and $\tau=10^{-4}$ (right).
    Density (top), velocity (middle top), temperature (middle bottom) and heat flux (bottom) are shown for $M=800$ cells and $N=64^2$ velocity cells.
   }
  \label{fig:test4_1}
\end{figure}

\begin{figure}[ht]
    \includegraphics[width=0.32\textwidth]{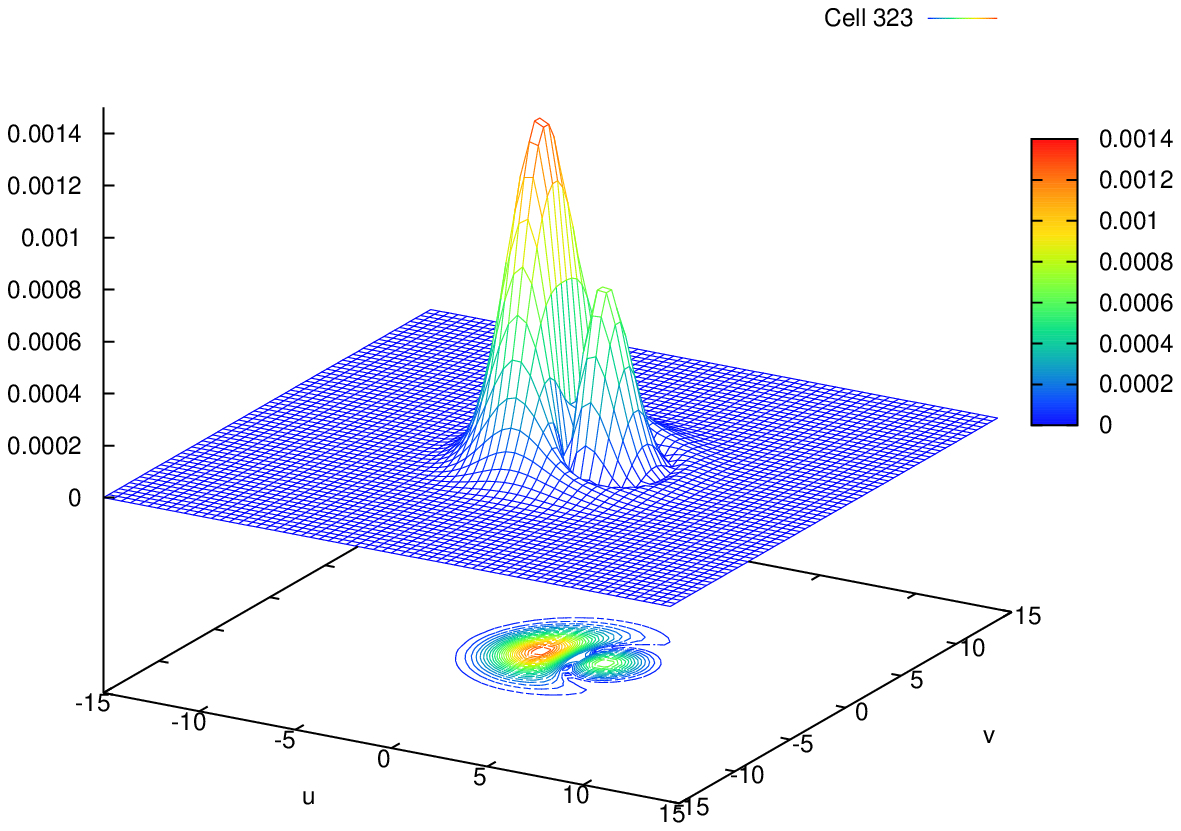}  
\hskip5.5cm 
    \includegraphics[width=0.32\textwidth]{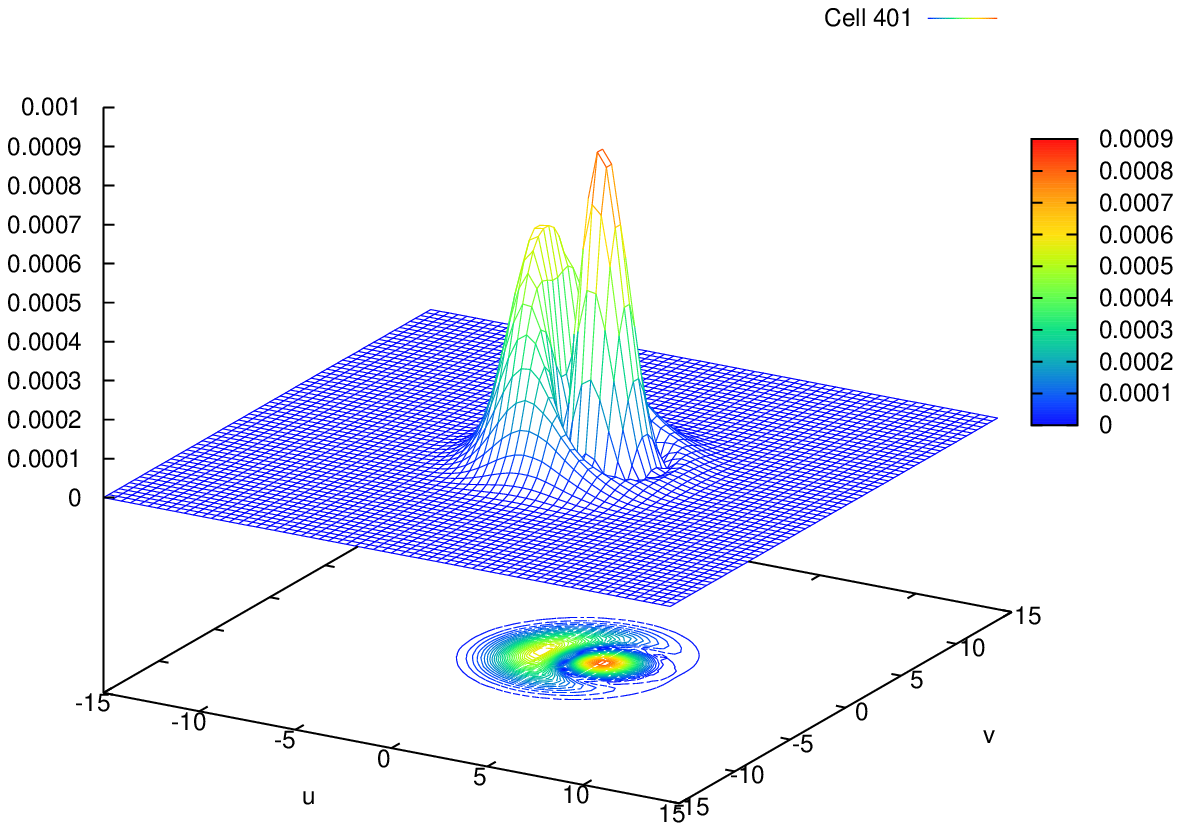} \\ 
    \includegraphics[width=0.32\textwidth]{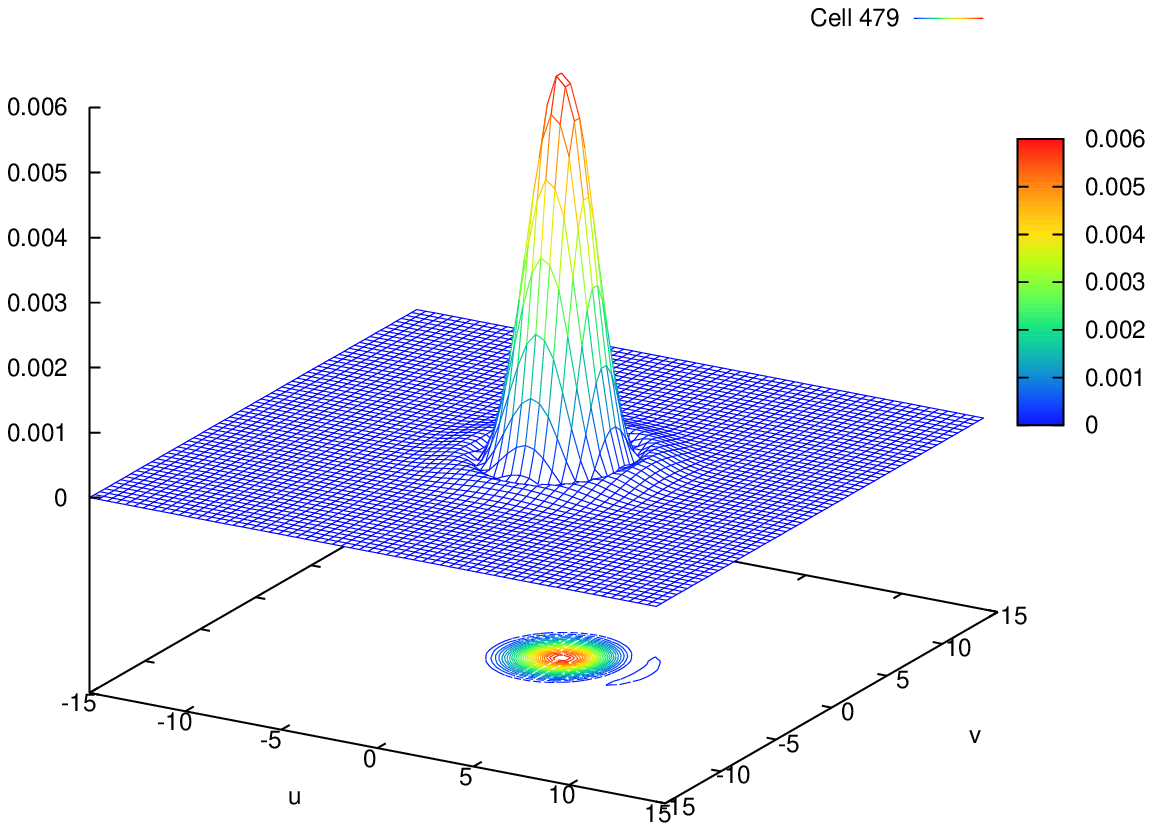} 
    {\includegraphics[width=0.32\textwidth]{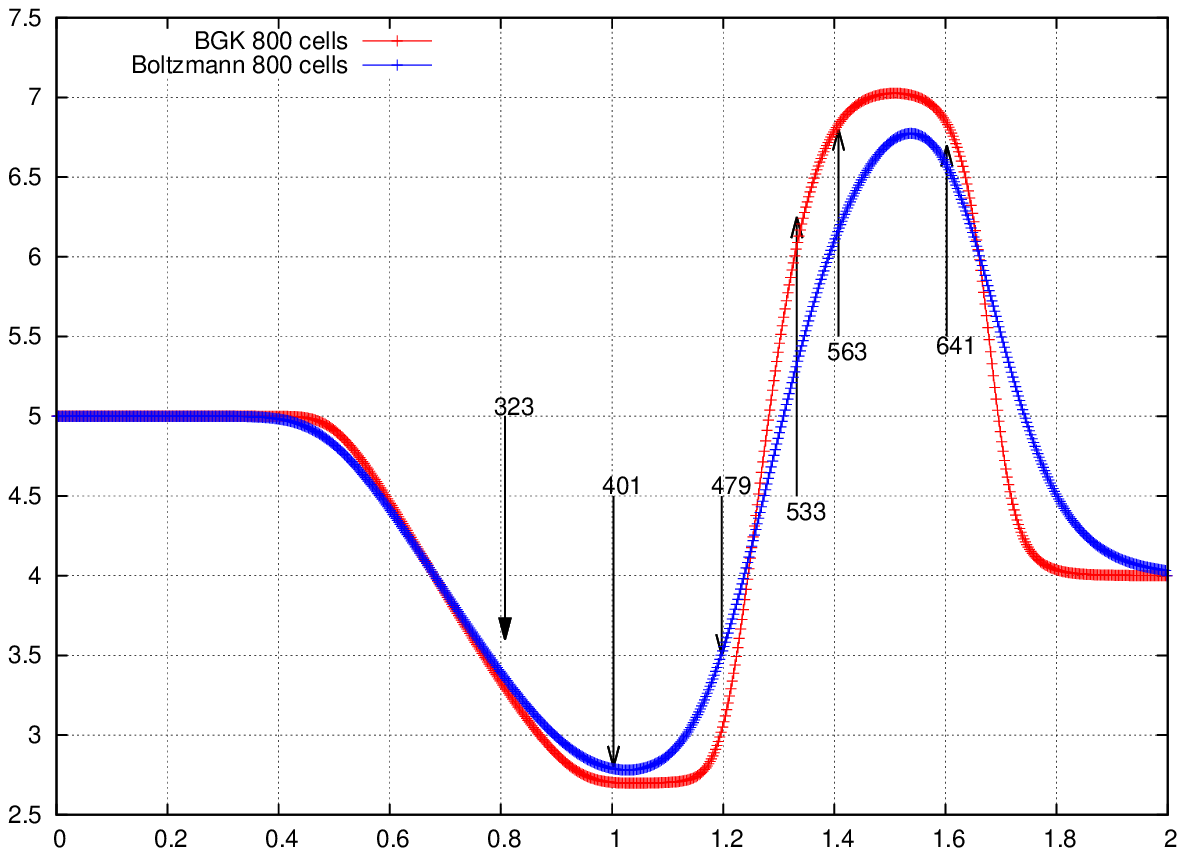}}
   \includegraphics[width=0.32\textwidth]{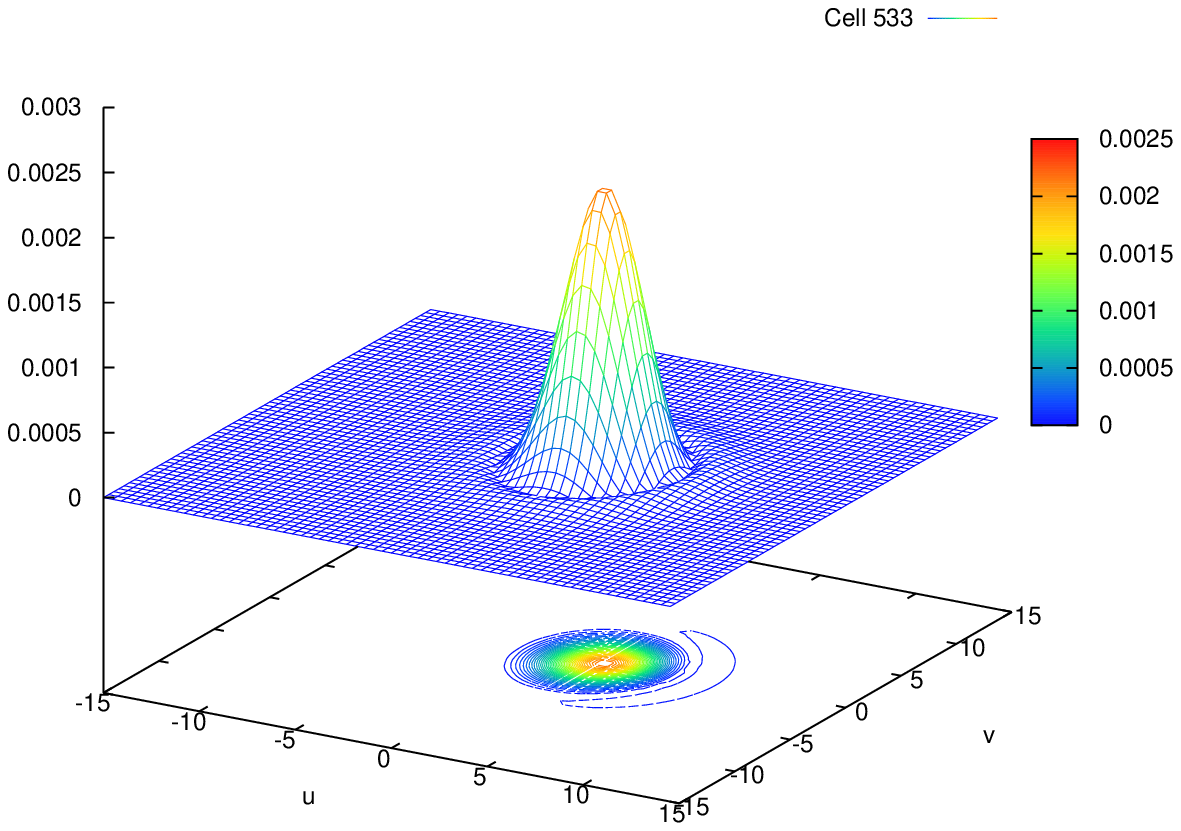}\\ 
      \includegraphics[width=0.32\textwidth]{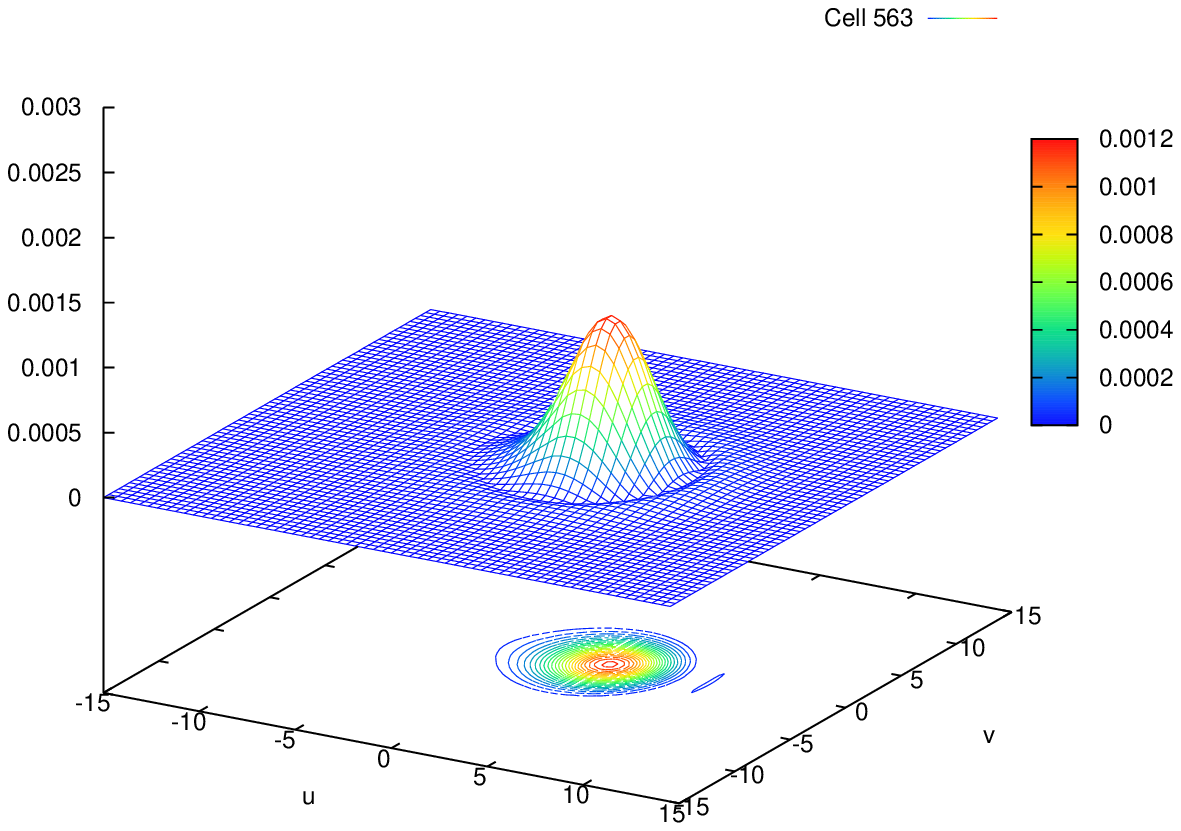} 
      \hskip5.5cm
    \includegraphics[width=0.32\textwidth]{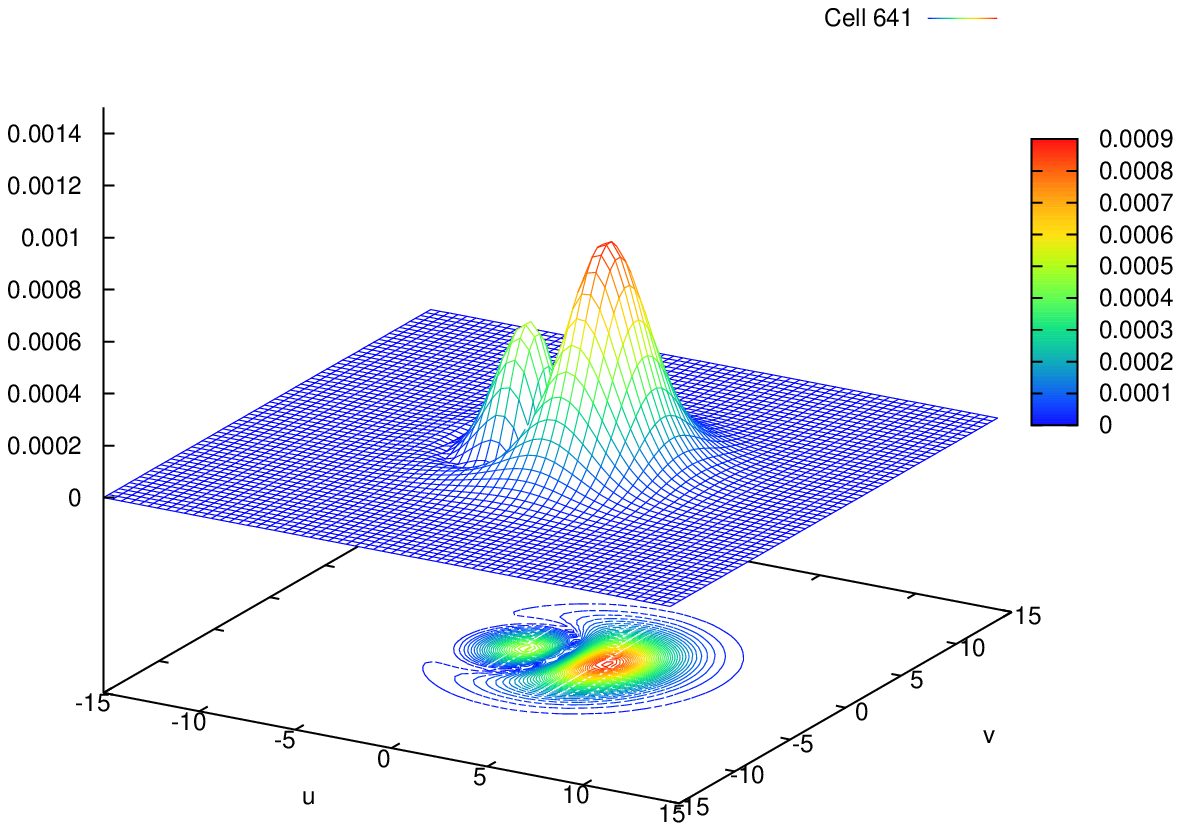} 
  \caption{Test 2.2. One dimension in space and two dimension in velocity at $t_{\text{final}} =0.15$ for BGK and Boltzmann models
    for $\tau=10^{-3}$ with $M=800$ and $N=64^2$ velocity cells.
    Middle panel: temperature for the two models. The arrows indicate the regions for which the difference of the 
    two distribution functions $|f^{\text{BGK}}(x_i,v)-f^{\text{Boltz}}(x_i,v)|$ is reported.
    The vertical scale is kept constant, only the color scale is adapted to the values.
   }
  \label{fig:test4_2}
\end{figure}

Here, the BGK and Boltzmann models are simulated for the same Sod-like problem. In order to have fairest as possible comparisons between the two models, we choose
$\nu=\rho$ for the BGK model. This choice permits to have the same loss part for the two models since for Maxwellian molecules the loss part is close to $\rho f$ as stated in the previous paragraph.
We fix $M=800$ spatial cells and $N=64^2$ cells in velocity space. This permits to consider almost converged results.
In Figure~\ref{fig:test4_1}, we present the results when $\tau=10^{-3}$ (left) and $\tau=10^{-4}$ (right) for the  density (top), the velocity (middle top), the temperature (middle bottom) and
the heat flux (bottom). The very first observation is the relative large differences between the two solutions. The Boltzmann solution systematically presents more dissipated waves than BGK solution. 
However, the main waves are located 
in the same positions for the two models. The difference in term of the macroscopic quantities can be of the order of $10\%-25\%$ of the solution. 
Systematically, the BGK model underestimates the values of the heat flux.
In Figure~\ref{fig:test4_2} we report the absolute value of the difference between the two distribution functions
$f_{\text{BGK}}(x,v,t)$ and $f_{\text{Boltz}}(x,v,t)$ at different locations $x_j$ and at final time $t_{\text{final}}$ as a function of the velocity variables $v$. 
The vertical scale is kept constant for all panels, only the color scale is adapted to the values reported.

\subsubsection{Test 2.3. Numerical convergence of the Boltzmann equation. The three dimensional in velocity hard sphere molecules case.}
\begin{figure}[ht]
  \begin{center}
    \includegraphics[width=0.36\textwidth]{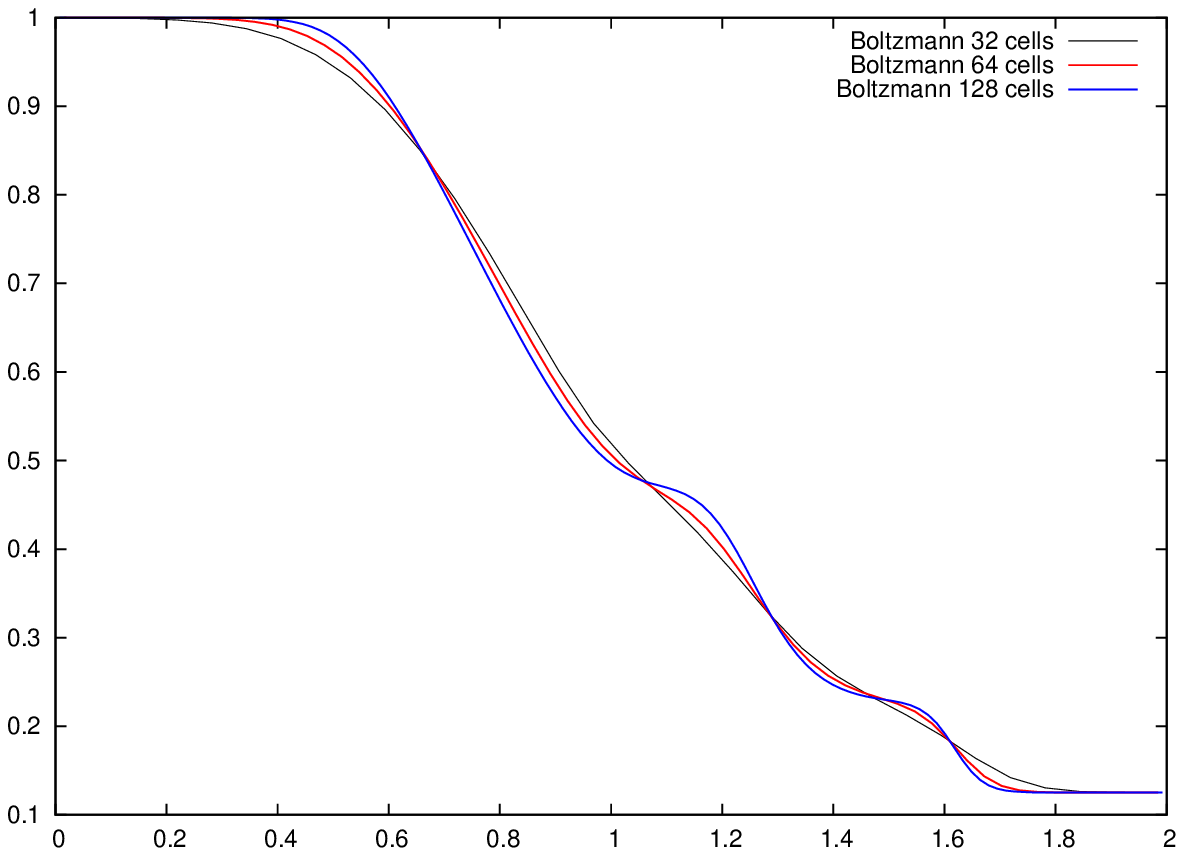} 
    \includegraphics[width=0.36\textwidth]{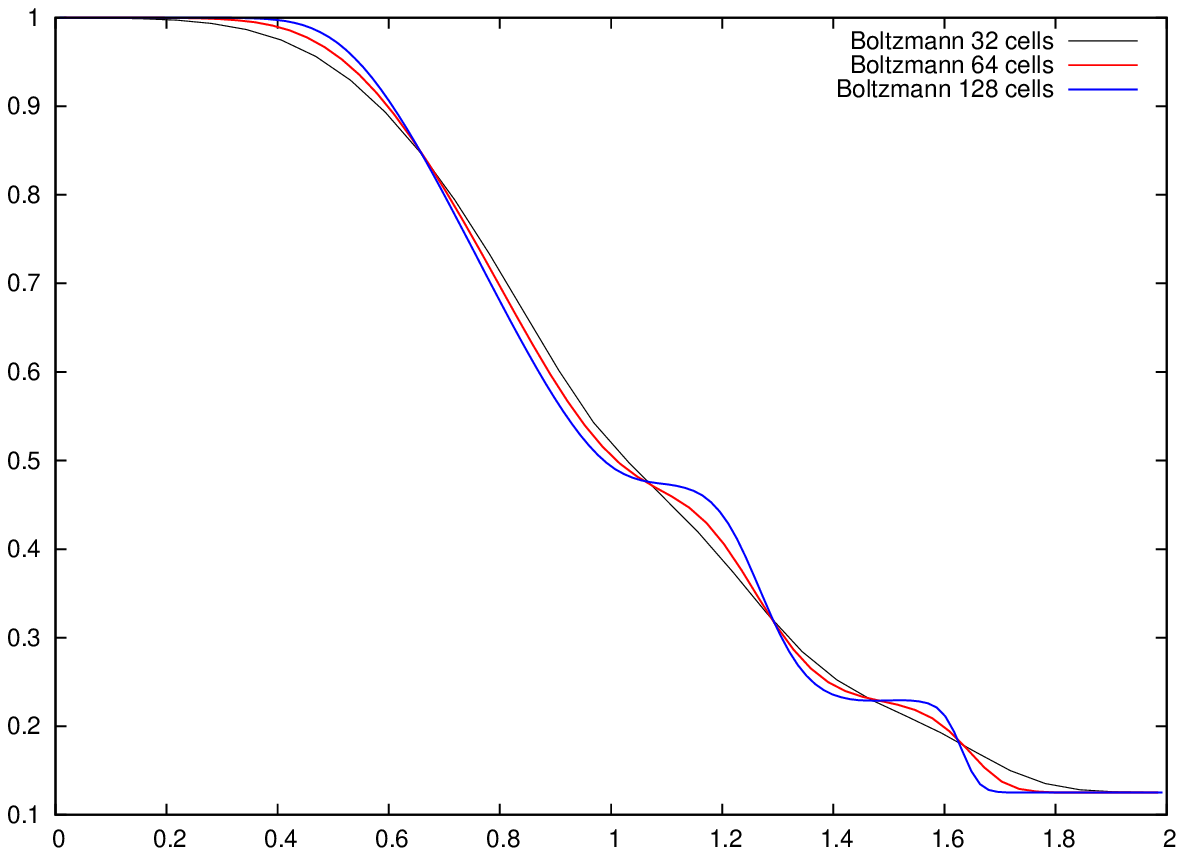} \\
    \includegraphics[width=0.36\textwidth]{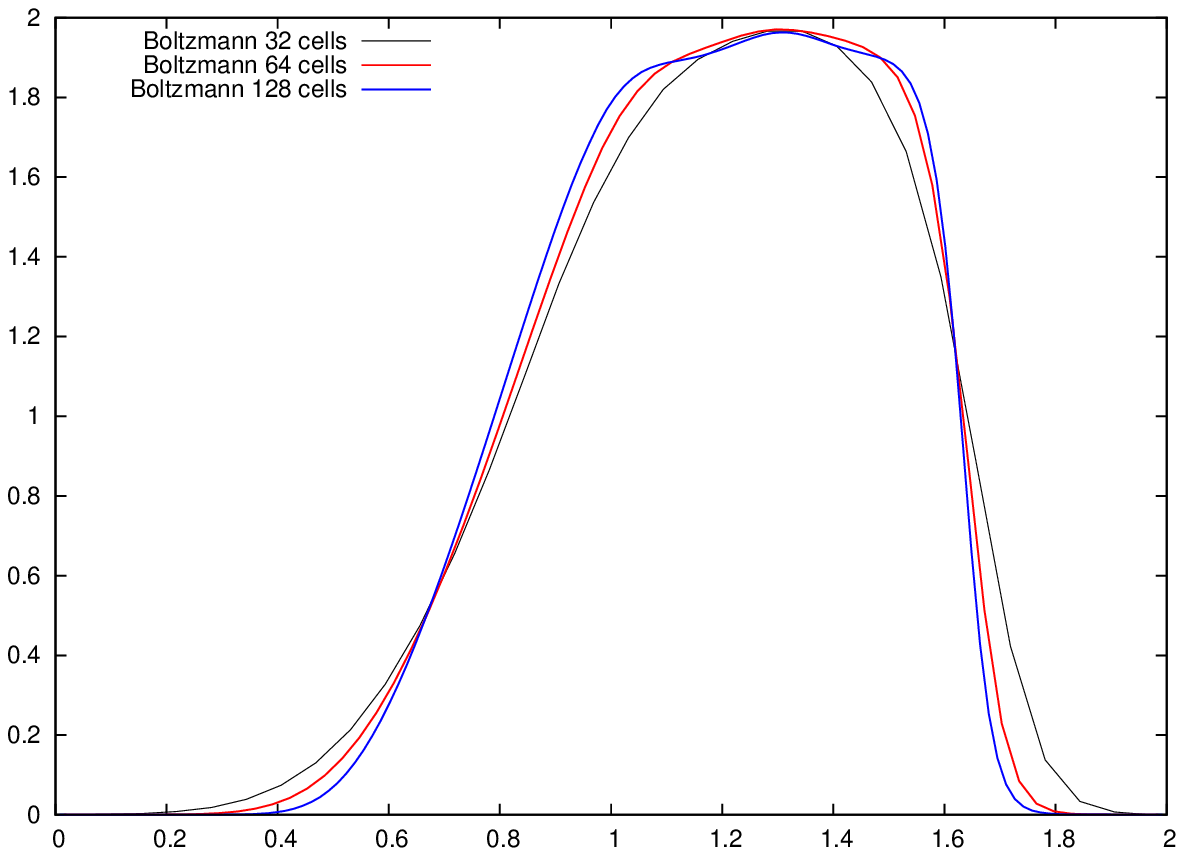} 
    \includegraphics[width=0.36\textwidth]{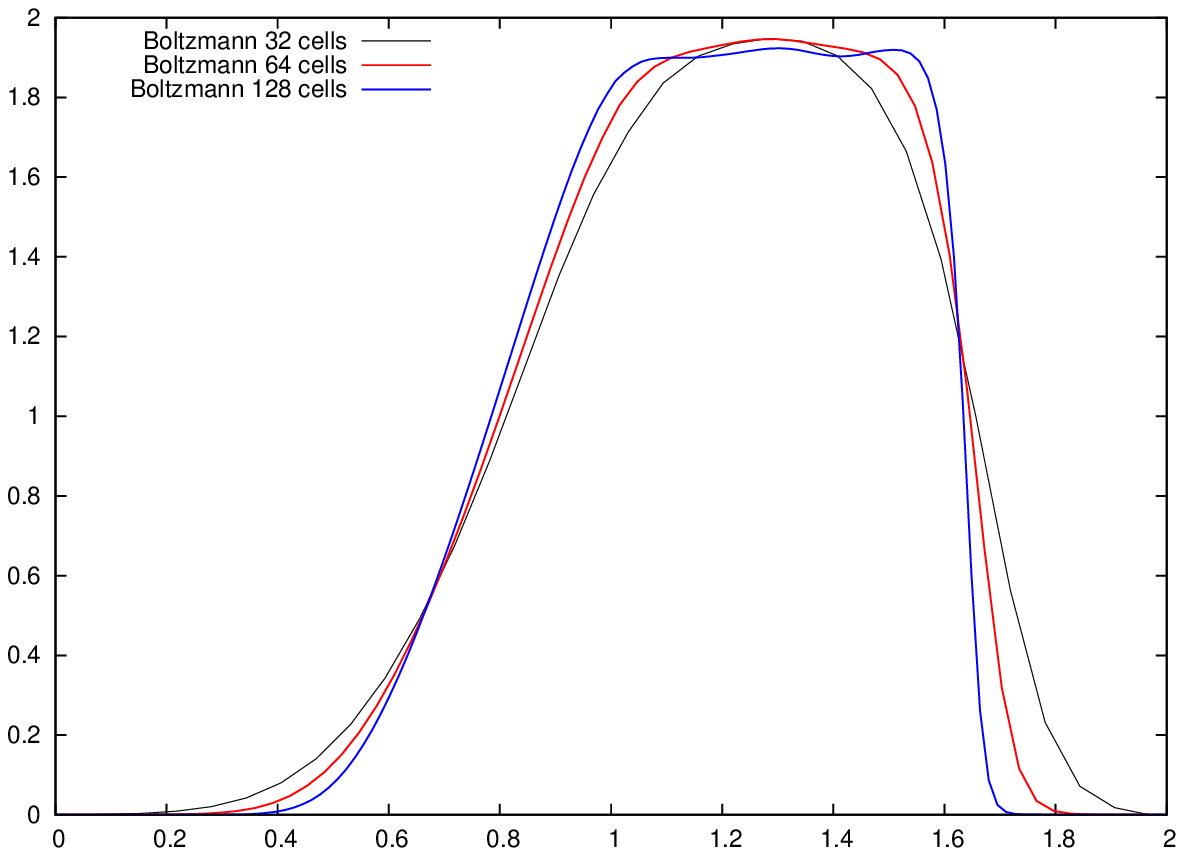} \\
    \includegraphics[width=0.36\textwidth]{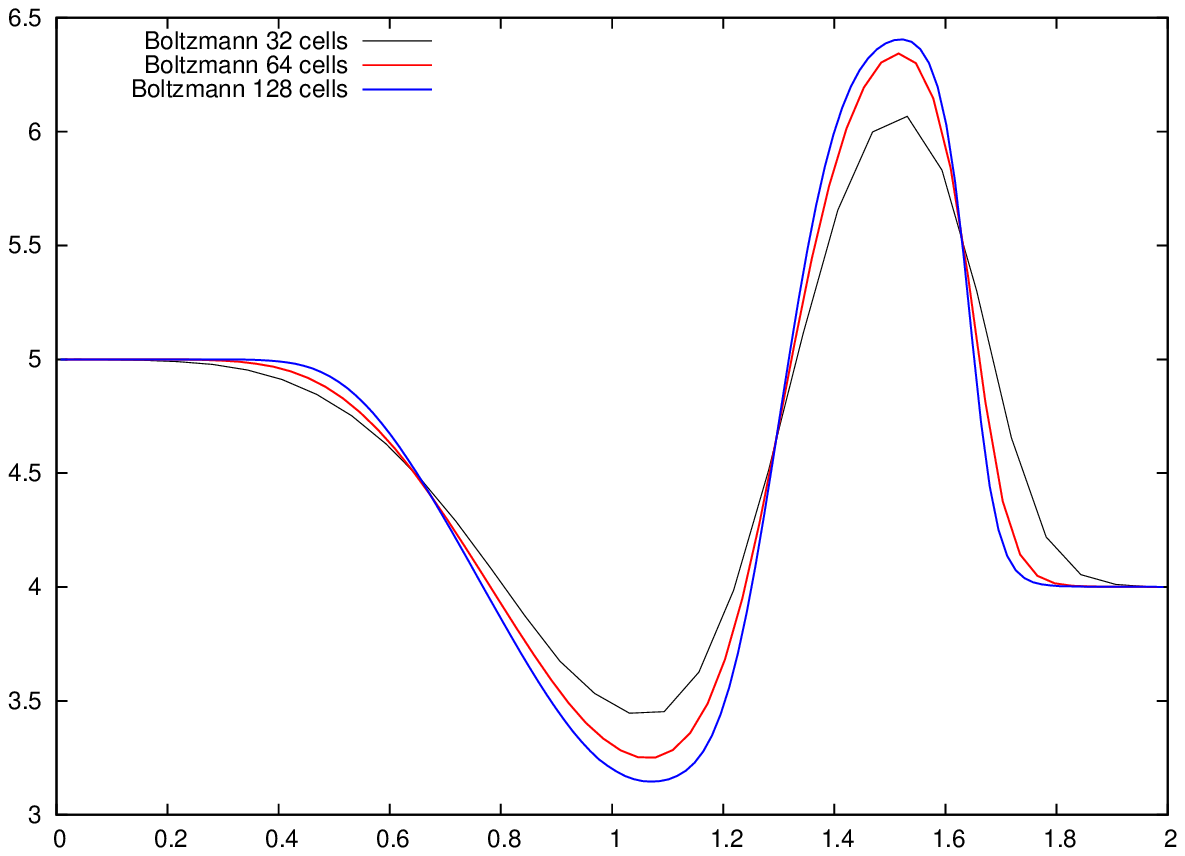} 
    \includegraphics[width=0.36\textwidth]{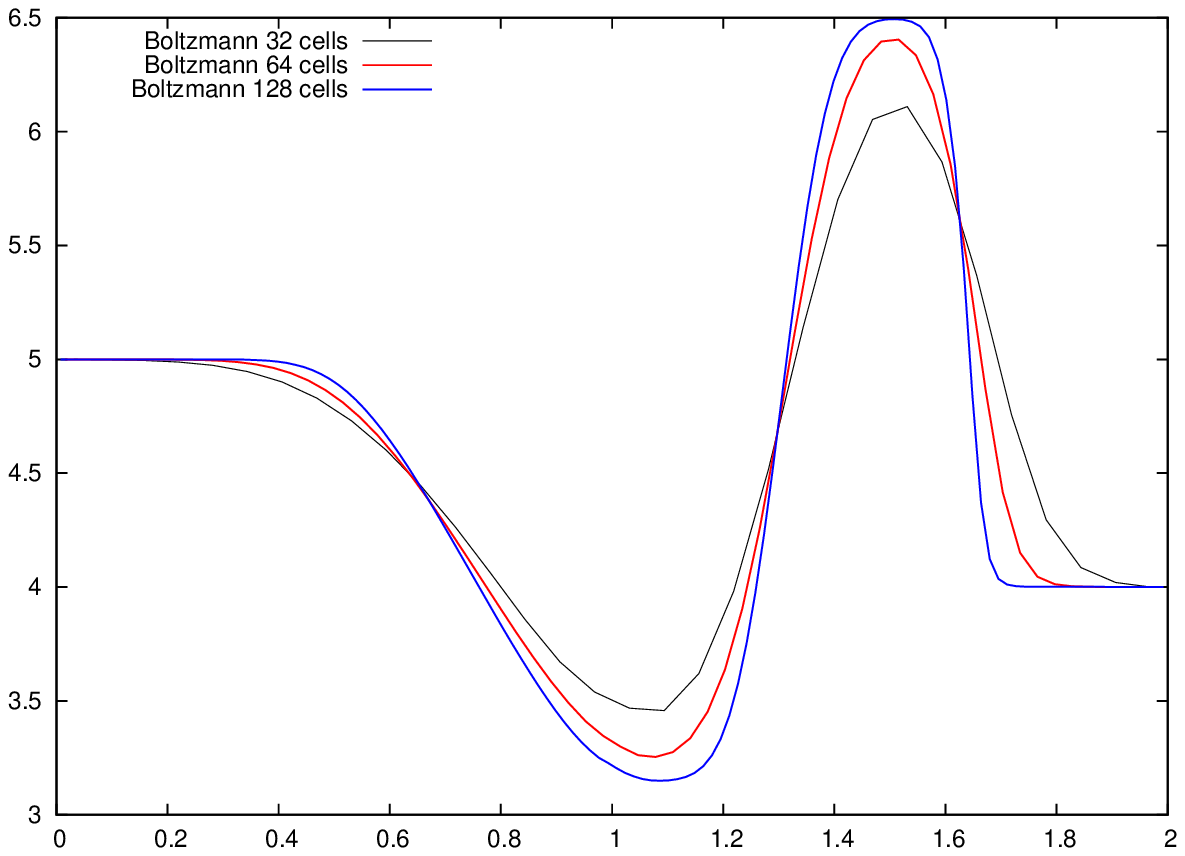} 
  \end{center}
  \caption{Test 2.3:  One dimension in space and three dimension in velocity for Sod like test case at $t_{\text{final}} =0.5$ 
    for  the hard sphere molecules case.
    Mesh convergence results for $\tau=10^{-2}$ (left) and $\tau=2 \ 10^{-3}$ (right).
    Density (top), velocity (middle) and temperature (bottom) are shown for $M=32, 64$ and $128$ cells and $N=32^3$ velocity cells.  }
  \label{fig:test5}
\end{figure}

Let us focus on the one dimension in space and three dimensions in velocity case. We consider the same Sod-like problem up to final time $t_{\text{final}}=0.5$
and two different 
collision frequencies $\tau=10^{-2}$ and $\tau=2 \ 10^{-3}$. The space/velocity mesh chosen is of the form $M\times N$
with $N=32^3$ uniformly spread on a velocity domain $[-16;16]^3$ and varying number of space cells $M$. In Figure~\ref{fig:test5}, we present the space convergence results for the density, the velocity and the temperature for successively refined spatial meshes from $32$ to $128$. From these data we can observe that the 
simulation results seem to converge towards the same numerical solution in both cases $\tau=10^{-2}$ (left panels) and $\tau=2 \ 10^{-3}$ (right panels). The CFL condition employed in this case is the following
\be \Delta t\leq \min \left(\frac{\Delta x}{|v_{max}|},\frac{\tau}{\mu}\right),\ee
where the second term is due to the stability restriction in the solution of the space homogeneous problem when hard sphere molecules are employed. 
In this case, the loss part of the collision integral $Q^-(f)$ can be only estimated, giving $Q^-(f)=L(f)f\leq C_\alpha 4\pi(2\lambda \pi)^\alpha$ and thus $\mu\geq C_\alpha 4\pi(2\lambda \pi)^\alpha$ in order to ensure that the gain part $Q^+(f)$ is positive and monotone.
We have also performed comparisons between the BGK and the Boltzmann models in this one dimensional setting. The results
are close to the ones obtained for the two dimensional Maxwell molecules, i.e. the BGK model tends to over-relax the distribution function to the equilibrium state and for this reason we do not
report them. However, since the computational costs involved in the approximation of the Boltzmann integral in the three dimensional setting are much larger than those obtained with
the simpler two dimensional model, in the next paragraph we analyze the performances of the scheme proposed in both cases.

\subsubsection{Performances}
In this part, we analyze the performances of our scheme in the one dimensional in space setting. In Table~\ref{tab:perf1DxND} we monitor the CPU time for the Boltzmann and BGK operators when a successively refined spatial mesh is considered for a fixed two and three dimensional velocity mesh. The data reported are relative to the same Sod-like problem considered in the previous paragraphs.
We employ $64^2$ velocity cells in the two dimensional case and $32^3$ in the three dimensional case.  
The number of spatial cells is $M=100$, $200$ and $400$ in the two dimensional case and $M=32$, $64$, and $128$ in the three dimensional one.
More precisely, we measure the CPU time $\text{T}$ in seconds, the time per cycle as $T_{\text{cycle}} = \text{T}/N_{\text{cycle}}$, the
time per cycle per spatial cell by $T_{\text{cell}} = \text{T}/N_{\text{cycle}}/M$
and the time per cycle per degree of freedom (d.o.f) $T_{\text{dof}} = \text{T}/N_{\text{cycle}}/(M \times N)$.
The OpenMP parallel version of the scheme is used on a laptop having 8 threads (HP ZBook with Intel 8 Core i7-4940MX CPU @ 3.10GHz 
on Ubuntu 15.10 (64 bits)).
These simulations have been run in parallel on non dedicated computer, as such the 
results are to be understood as rough estimates. By no mean we pretend that they cannot be improved.
In the two dimensional case in velocity, the Boltzmann model results are about $15$ times more expensive than the BGK model, while in three dimensions, 
it can reach about $250$ times. Notice that, as expected, these ratio do increase when the collision frequency $\tau$ becomes smaller, for instance when $\tau=2 \times 10^{-3}$ in three dimensions
then  the ratio in CPU time between the two models is of the order of $500$. 
Let us notice that the number of cycles for hard sphere molecules is fixed for the reported simulations. This is due to the fact that
the hard sphere model considered has a stability condition which is more restrictive than the CFL condition chosen for the transport part as opposite to the Maxwellian molecule model
for the two dimensional case. For fair comparison, the number of cycles of the two models has been kept constant in the computation of the costs, even if, 
the space homogeneous BGK model does not present stability requirements.
\begin{table}
\begin{center}
  \begin{tabular}{|p{0.35cm}|cc|c||c|c|c|c|c|}
    \hline
    \multicolumn{9}{|c|}{\textbf{Sod like Riemann problem in 1D$\times$2D} $\tau=10^{-3}$} \\
    \hline
    \hline
    \multirow{3}{*}{\begin{sideways}\textbf{Model}\end{sideways}}
    & \multicolumn{2}{|c|}{\textbf{Velocity}}  
      &\textbf{Cell \#}    & \textbf{Cycle} & \textbf{Time} & T\textbf{/cycle} & T\textbf{/cell} & T/\textbf{d.o.f} \\
    & $N$ & {\begin{sideways} Vel.\end{sideways}}&$M \times N$ &  $N_{\text{cycle}}$  &  $\text{T}$ (s)  & $T_{\text{cycle}}$ (s) & $T_{\text{cell}}$ (s) & $T_{\text{dof}}$ (s)\\
    \hline
     \hline
     \multirow{6}{*}{\begin{sideways}  \textbf{BGK} \end{sideways}}  
     & \multirow{6}{*}{ $64^2$} &
     \multirow{6}{*}{\begin{sideways} $[-15,15]$  \end{sideways}}
     &$100\times 64^2$  & \multirow{2}{*}{$111$}  & \multirow{2}{*}{$0.935$} & \multirow{2}{*}{$8.42\times 10^{-3}$}  & \multirow{2}{*}{$8.42\times 10^{-5}$}& \multirow{2}{*}{$2.06\times 10^{-8}$} \\
     &&&$\simeq 4.1 \times 10^5$ & &  & & &  \\
     \cline{4-9}
     &&&$200 \times 64^2$  & \multirow{2}{*}{$222$}& \multirow{2}{*}{$2.738$} & \multirow{2}{*}{$1.23\times 10^{-2}$} & \multirow{2}{*}{$6.17\times 10^{-5}$}& \multirow{2}{*}{$1.51\times 10^{-8}$}\\
     &&&$\simeq 8.2\times 10^{5}$ & & & & &  \\
     \cline{4-9}
     &&&$400 \times 64^2$ & \multirow{2}{*}{$436$} & \multirow{2}{*}{$7.822$} & \multirow{2}{*}{$1.77\times 10^{-2}$} & \multirow{2}{*}{$4.41\times 10^{-5}$}& \multirow{2}{*}{$1.08\times 10^{-8}$} \\
     &&&$\simeq 16.4\times 10^5$ & &  &  & &  \\
     \hline
     \hline
     \multirow{6}{*}{\begin{sideways}  \textbf{Boltzmann} \end{sideways}}  
     & \multirow{6}{*}{ $64^2$} &
     \multirow{6}{*}{\begin{sideways} $[-15,15]$  \end{sideways}}
     &$100\times 64^2$  & \multirow{2}{*}{$111$}  & \multirow{2}{*}{$8.174$}  & \multirow{2}{*}{$7.36\times 10^{-2}$} & \multirow{2}{*}{$7.36\times 10^{-4}$}& \multirow{2}{*}{$5.64\times 10^{-8}$} \\ 
     &&&$\simeq 4.1 \times 10^5$ & &  & & &  \\
     \cline{4-9}
     &&&$200 \times 64^2$  & \multirow{2}{*}{$218$}& \multirow{2}{*}{$28.095$} & \multirow{2}{*}{$1.27\times 10^{-1}$} & \multirow{2}{*}{$6.33\times 10^{-4}$}& \multirow{2}{*}{$4.42\times 10^{-8}$}\\
     &&&$\simeq  8.2\times 10^{5}$ & & & & & \\
     \cline{4-9}
    &&&$400 \times 64^2$ & \multirow{2}{*}{$443$} & \multirow{2}{*}{$113.495$} & \multirow{2}{*}{$2.56\times 10^{-1}$}& \multirow{2}{*}{$6.40\times 10^{-4}$}& \multirow{2}{*}{$1.56\times 10^{-7}$} \\
     &&&$\simeq  16.4\times 10^5$ & &  &  & & \\
     \hline   
     \hline
     \hline
     \multicolumn{9}{|c|}{\textbf{Sod like Riemann problem in 1D$\times$3D} $\tau=10^{-2}$} \\
     \hline
     \hline
    \multirow{3}{*}{\begin{sideways}\textbf{Model}\end{sideways}}
    & \multicolumn{2}{|c|}{\textbf{Velocity}}  
      &\textbf{Cell \#}    & \textbf{Cycle} & \textbf{Time} & T\textbf{/cycle} & T\textbf{/cell} & T/\textbf{d.o.f} \\
    & $N$ & {\begin{sideways} Vel.\end{sideways}}&$M \times N$ &  $N_{\text{cycle}}$  &  $\text{T}$ (s)  & $T_{\text{cycle}}$ (s) & $T_{\text{cell}}$ (s) & $T_{\text{dof}}$ (s)\\
    \hline
     \hline
     \multirow{6}{*}{\begin{sideways}  \textbf{BGK} \end{sideways}}  
     & \multirow{6}{*}{ $32$} &
     \multirow{6}{*}{\begin{sideways} $[-16,16]$  \end{sideways}}
     &$32\times 32^2$  & \multirow{2}{*}{$395$}   & \multirow{2}{*}{$12.64$}   & \multirow{2}{*}{$3.20\times 10^{-2}$} & \multirow{2}{*}{$1.00\times 10^{-3}$} & \multirow{2}{*}{$3.05\times 10^{-8}$}\\
     &&&$\simeq 1.05 \times 10^6$ & &  & & &  \\
     \cline{4-9}
     &&&$64 \times 32^2$  & \multirow{2}{*}{$395$} & \multirow{2}{*}{$19.55$} & \multirow{2}{*}{$4.95\times 10^{-2}$} & \multirow{2}{*}{$7.73\times 10^{-4}$} & \multirow{2}{*}{$2.36\times 10^{-8}$}\\
     &&&$\simeq 2.10\times 10^{6}$ & & & & &  \\
     \cline{4-9}
     &&&$128 \times 32^2$ & \multirow{2}{*}{$395$} & \multirow{2}{*}{$40.75$} & \multirow{2}{*}{$1.03\times 10^{-1}$} & \multirow{2}{*}{$8.06\times 10^{-4}$}& \multirow{2}{*}{$2.46\times 10^{-8}$} \\
     &&&$\simeq 4.19\times 10^6$ & &  &  & &  \\
     \hline
     \hline
     \multirow{6}{*}{\begin{sideways}  \textbf{Boltzmann} \end{sideways}}  
     & \multirow{6}{*}{ $32$} &
     \multirow{6}{*}{\begin{sideways} $[-16,16]$  \end{sideways}}
     &$32\times 32^2$  & \multirow{2}{*}{$395$}    & $2104.17$ & \multirow{2}{*}{$5.33$}& \multirow{2}{*}{$1.66\times 10^{-1}$}& \multirow{2}{*}{$5.08\times 10^{-6}$} \\ 
     &&&$\simeq 1.05 \times 10^6$ & & $\sim 35$mn & & &  \\
     \cline{4-9}
     &&&$64 \times 32^2$  & \multirow{2}{*}{$395$} & $4799.19$ & \multirow{2}{*}{$12.10$} & \multirow{2}{*}{$1.90\times 10^{-1}$}& \multirow{2}{*}{$5.79\times 10^{-6}$} \\
     &&&$\simeq 2.10\times 10^{6}$ & & $\sim 1.3$h & & & \\
     \cline{4-9}
     &&&$128 \times 32^2$ & \multirow{2}{*}{$395$} & $9114.56$ & \multirow{2}{*}{$23.10$} & \multirow{2}{*}{$1.80\times 10^{-1}$}& \multirow{2}{*}{$5.50\times 10^{-6}$} \\
     &&&$\simeq 4.19\times 10^6$ & & $\sim 2.5$h &  & & \\
    \hline
    \hline
  \end{tabular}
  \end{center} 
  \caption{ \label{tab:perf1DxND}
    One dimensional in space, two dimensional in velocity Maxwellian molecules and three dimensional in velocity hard sphere molecules simulations. 
    Comparisons between the BGK and Boltzmann models for spatial mesh variation.
    Monitoring of CPU time.
    Time per cycle is obtained by $T_{\text{cycle}} = \text{T}/N_{\text{cycle}}$,
    time per cycle per cell by $T_{\text{cell}} = \text{T}/N_{\text{cycle}}/M$ and 
    time per cycle per degree of freedom $T_{\text{dof}} = \text{T}/N_{\text{cycle}}/(M \times N)$.
  }
\end{table}
Next in table~\ref{tab:perf1DxND_vel} we monitor the CPU time of the simulation when the number of cells is fixed in space to $M=100$.
The number of velocity cells increases from $N=8^2$ to $64^2$ in two and from $N=8^3$ to $64^3$ in three dimensions. The same Sod test case as in table~\ref{tab:perf1DxND} is simulated.
Because the CPU times have been obtained by parallel simulations, the analysis of such table must take into account the fact that the
dimension of the mesh plays a role in the performance. For instance, in two dimensions, the smaller mesh ($100\times 8^2$)
usually presents a larger CPU time per degree of freedom. This is due to the fact that a decent amount of cells is needed to observe the benefit of using a parallel machine.
For the Boltzmann simulations, in three dimensions, $8$ velocity cells in each direction is not large enough to have a stable simulation. This is probably due to the large loss of energy
caused by the spectral method which can be only partly cured by the $L_2$ projection technique detailed in Section $3.2$. In fact, even if such renormalization permits to keep the correct
energy after the collision step, it may transform the distribution function in an unphysical manner. This operation repeated multiple times give rise to instability in the spectral scheme. 
However, for completeness using the first completed iterations we report in the table an approximation of the CPU time that would be needed using this small amount of points.
The Hard Sphere molecule simulations, as can be seen on the Figures reported, are extremely time consuming already in the one dimensional case in space.
In two dimensions, the cost of the Boltzmann model is about $3$ to $8$ times (from the smallest to largest mesh) more expensive than the BGK model.
In three dimensions, this ratio ranges between $40$ and $700$ times.
\begin{table}
\begin{center}
  \begin{tabular}{|p{0.35cm}|cc|c||c|c|c|c|c|}
    \hline
    \multicolumn{9}{|c|}{\textbf{Sod like Riemann problem in 1D$\times$2D} $\tau=10^{-3}$} \\
    \hline
    \hline
    \multirow{3}{*}{\begin{sideways}\textbf{Model}\end{sideways}}
    & \multicolumn{2}{|c|}{\textbf{Velocity}}  
      &\textbf{Cell \#}    & \textbf{Cycle} & \textbf{Time} & T\textbf{/cycle} & T\textbf{/cell} & T/\textbf{d.o.f} \\
    & $N$ & {\begin{sideways} Vel.\end{sideways}}&$M \times N$ &  $N_{\text{cycle}}$  &  $\text{T}$ (s)  & $T_{\text{cycle}}$ (s) & $T_{\text{cell}}$ (s) & $T_{\text{dof}}$ (s)\\
    \hline
     \hline
     \multirow{8}{*}{\begin{sideways}  \textbf{BGK} \end{sideways}}  
     & \multirow{2}{*}{ $8^2$} &
     \multirow{8}{*}{\begin{sideways} $[-15,15]$  \end{sideways}}
                               &$100\times 8^2$  & \multirow{2}{*}{$105$}  & \multirow{2}{*}{$0.126$} & \multirow{2}{*}{$1.20\times 10^{-3}$}  & \multirow{2}{*}{$1.20\times 10^{-5}$}& \multirow{2}{*}{$1.87\times 10^{-7}$} \\
     & &&$\simeq 6.4 \times 10^4$ & &  & & &  \\
     \cline{4-9}
     &\multirow{2}{*}{ $16^2$} &&$100 \times 16^2$  & \multirow{2}{*}{$106$}& \multirow{2}{*}{$0.265$} & \multirow{2}{*}{$2.50\times 10^{-3}$} & \multirow{2}{*}{$2.50\times 10^{-5}$}& \multirow{2}{*}{$9.78\times 10^{-8}$}\\
     &&&$\simeq 2.56\times 10^{4}$ & & & & &  \\
     \cline{4-9}
     &\multirow{2}{*}{ $32^2$} &&$100 \times 32^2$  & \multirow{2}{*}{$109$}& \multirow{2}{*}{$0.291$} & \multirow{2}{*}{$2.67\times 10^{-3}$} & \multirow{2}{*}{$2.67\times 10^{-5}$}& \multirow{2}{*}{$2.61\times 10^{-8}$}\\
     &&&$\simeq 1.024\times 10^{5}$ & & & & &  \\
     \cline{4-9}
     &\multirow{2}{*}{ $64^2$} &&$100 \times 64^2$ & \multirow{2}{*}{$111$} & \multirow{2}{*}{$0.935$} & \multirow{2}{*}{$8.42\times 10^{-3}$} & \multirow{2}{*}{$8.42\times 10^{-5}$}& \multirow{2}{*}{$2.06\times 10^{-8}$} \\
     &&&$\simeq 4.096\times 10^5$ & &  &  & &  \\
     \hline
     \hline
     \multirow{8}{*}{\begin{sideways}  \textbf{Boltzmann} \end{sideways}}  
     & \multirow{2}{*}{ $8^2$} &
     \multirow{8}{*}{\begin{sideways} $[-15,15]$  \end{sideways}}
                               &$100\times 8^2$  & \multirow{2}{*}{$105$}  & \multirow{2}{*}{$0.387$} & \multirow{2}{*}{$3.69\times 10^{-3}$}  & \multirow{2}{*}{$3.69\times 10^{-5}$}& \multirow{2}{*}{$5.76\times 10^{-7}$} \\
     & &&$\simeq 6.4 \times 10^4$ & &  & & &  \\
     \cline{4-9}
     &\multirow{2}{*}{ $16^2$} &&$100 \times 16^2$  & \multirow{2}{*}{$106$}& \multirow{2}{*}{$0.733$} & \multirow{2}{*}{$3.69\times 10^{-3}$} & \multirow{2}{*}{$3.69\times 10^{-5}$}& \multirow{2}{*}{$5.76\times 10^{-7}$}\\
     &&&$\simeq 2.56\times 10^{4}$ & & & & &  \\
     \cline{4-9}
     &\multirow{2}{*}{ $32^2$} &&$100 \times 32^2$  & \multirow{2}{*}{$109$}& \multirow{2}{*}{$2.302$} & \multirow{2}{*}{$2.11\times 10^{-2}$} & \multirow{2}{*}{$2.11\times 10^{-4}$}& \multirow{2}{*}{$2.06\times 10^{-7}$}\\
     &&&$\simeq 1.024\times 10^{5}$ & & & & &  \\
     \cline{4-9}
     &\multirow{2}{*}{ $64^2$} &&$100 \times 64^2$ & \multirow{2}{*}{$111$} & \multirow{2}{*}{$7.834$} & \multirow{2}{*}{$7.06\times 10^{-2}$} & \multirow{2}{*}{$7.06\times 10^{-4}$}& \multirow{2}{*}{$1.72\times 10^{-7}$} \\
     &&&$\simeq 4.096\times 10^5$ & &  &  & &  \\
     \hline   
     \hline
     \hline
     \multicolumn{9}{|c|}{\textbf{Sod like Riemann problem in 1D$\times$3D} $\tau=10^{-2}$} \\
     \hline
     \hline
    \multirow{3}{*}{\begin{sideways}\textbf{Model}\end{sideways}}
    & \multicolumn{2}{|c|}{\textbf{Velocity}}  
      &\textbf{Cell \#}    & \textbf{Cycle} & \textbf{Time} & T\textbf{/cycle} & T\textbf{/cell} & T/\textbf{d.o.f} \\
    & $N$ & {\begin{sideways} Vel.\end{sideways}}&$M \times N$ &  $N_{\text{cycle}}$  &  $\text{T}$ (s)  & $T_{\text{cycle}}$ (s) & $T_{\text{cell}}$ (s) & $T_{\text{dof}}$ (s)\\
    \hline
     \hline
     \multirow{8}{*}{\begin{sideways}  \textbf{BGK} \end{sideways}}  
     & \multirow{2}{*}{ $8^3$} &
     \multirow{8}{*}{\begin{sideways} $[-15,15]$  \end{sideways}}
                               &$100\times 8^3$  & \multirow{2}{*}{$395$}  & \multirow{2}{*}{$0.126$} & \multirow{2}{*}{$1.20\times 10^{-3}$}  & \multirow{2}{*}{$1.20\times 10^{-5}$}& \multirow{2}{*}{$1.87\times 10^{-7}$} \\
     & &&$\simeq 5.12 \times 10^4$ & &  & & &  \\
     \cline{4-9}
     &\multirow{2}{*}{ $16^3$} &&$100 \times 16^3$  & \multirow{2}{*}{$395$}& \multirow{2}{*}{$0.265$} & \multirow{2}{*}{$2.50\times 10^{-3}$} & \multirow{2}{*}{$2.50\times 10^{-5}$}& \multirow{2}{*}{$9.78\times 10^{-8}$}\\
     &&&$\simeq 4.10\times 10^{5}$ & & & & &  \\
     \cline{4-9}
     &\multirow{2}{*}{ $32^3$} &&$100 \times 32^3$  & \multirow{2}{*}{$395$}& \multirow{2}{*}{$0.291$} & \multirow{2}{*}{$2.67\times 10^{-3}$} & \multirow{2}{*}{$2.67\times 10^{-5}$}& \multirow{2}{*}{$2.61\times 10^{-8}$}\\
     &&&$\simeq 3.28\times 10^{6}$ & & & & &  \\
     \cline{4-9}
     &\multirow{2}{*}{ $64^3$} &&$100 \times 64^3$ & \multirow{2}{*}{$395$} & \multirow{2}{*}{$0.935$} & \multirow{2}{*}{$8.42\times 10^{-3}$} & \multirow{2}{*}{$8.42\times 10^{-5}$}& \multirow{2}{*}{$2.06\times 10^{-8}$} \\
     &&&$\simeq 2.62\times 10^7$ & &  &  & &  \\
     \hline
     \hline
     \multirow{8}{*}{\begin{sideways}  \textbf{Boltzmann} \end{sideways}}  
     & \multirow{2}{*}{ $8^3$} &
     \multirow{8}{*}{\begin{sideways} $[-15,15]$  \end{sideways}}
                               &$100\times 8^3$  & \multirow{2}{*}{$395$}  & FAIL & FAIL & FAIL & FAIL \\
     & &&$\simeq 5.12 \times 10^4$ & & $(\sim 25)$ & $(\sim 6.50\times 10^{-2})$ & $(\sim 6.50\times 10^{-4})$ &  $(\sim 1.27\times 10^{-6})$ \\
     \cline{4-9}
     &\multirow{2}{*}{ $16^3$} &&$100 \times 16^3$  & \multirow{2}{*}{$395$}& $241$ & \multirow{2}{*}{$6.11\times 10^{-1}$} & \multirow{2}{*}{$6.11\times 10^{-3}$}& \multirow{2}{*}{$1.49\times 10^{-6}$}\\
     &&&$\simeq 4.10\times 10^{5}$ & & $\sim 4$mn & & &  \\
     \cline{4-9}
     &\multirow{2}{*}{ $32^3$} &&$100 \times 32^3$  & \multirow{2}{*}{$395$}& $6438$ & \multirow{2}{*}{$1.63\times 10^{1}$} & \multirow{2}{*}{$1.63\times 10^{-1}$}& \multirow{2}{*}{$4.97\times 10^{-6}$}\\
     &&&$\simeq 3.28\times 10^{6}$ & & $\sim 1.8$h& & &  \\
     \cline{4-9}
     &\multirow{2}{*}{ $64^3$} &&$100 \times 64^3$ & \multirow{2}{*}{$395$} & $79099$ & \multirow{2}{*}{$2.00\times 10^{2}$} & \multirow{2}{*}{$2.00\times 10^{0}$}& \multirow{2}{*}{$7.64\times 10^{-6}$} \\
     &&&$\simeq 2.62\times 10^7$ & & $\sim 22$h  &  & &  \\
    \hline
    \hline
  \end{tabular}
  \end{center} 
  \caption{ \label{tab:perf1DxND_vel}
    One dimensional in space, two dimensional in velocity Maxwellian molecules and three dimensional in velocity hard sphere molecules simulations. 
    Simulations are performed using the OpenMP version of the scheme run on 8 hreads.
    Comparisons between the BGK and Boltzmann models for velocity mesh variation.
    Monitoring of CPU time.
    Time per cycle is obtained by $T_{\text{cycle}} = \text{T}/N_{\text{cycle}}$,
    time per cycle per cell by $T_{\text{cell}} = \text{T}/N_{\text{cycle}}/M$ and 
    time per cycle per degree of freedom $T_{\text{dof}} = \text{T}/N_{\text{cycle}}/(M \times N)$.
  }
\end{table}

%
%
\subsection{Part 3. Numerical results for the space two dimensional case.
} \label{sec:act2}
In this part, we focus on solving the two dimensional in space and velocity Boltzmann and BGK equations. The purposes are twofold. 
First, we want to show the differences which arises between the two models. Second, we want to analyze the performances of the method
in the two dimensional setting by monitoring the cost of such simulations and by performing a profiling of the scheme in terms of the main routines in order
to highlight the eventual bottlenecks. This permits to understand in which part of the scheme one should concentrate to improve the efficiency in the future.
For all reported simulations Maxwellian molecules are considered for the Boltzmann model. 
The CFL condition employed is as for the one dimensional case and for all tests the following
\be \Delta t\leq \min \left(\frac{\Delta x}{|v_{max}|},\frac{\tau}{\rho}\right).\ee

\subsubsection{Test 3.1. Two dimensional vortex in motion.}
\begin{figure}[h]
  \begin{center}
    \hspace{-1cm}
    \includegraphics[width=0.33\textwidth]{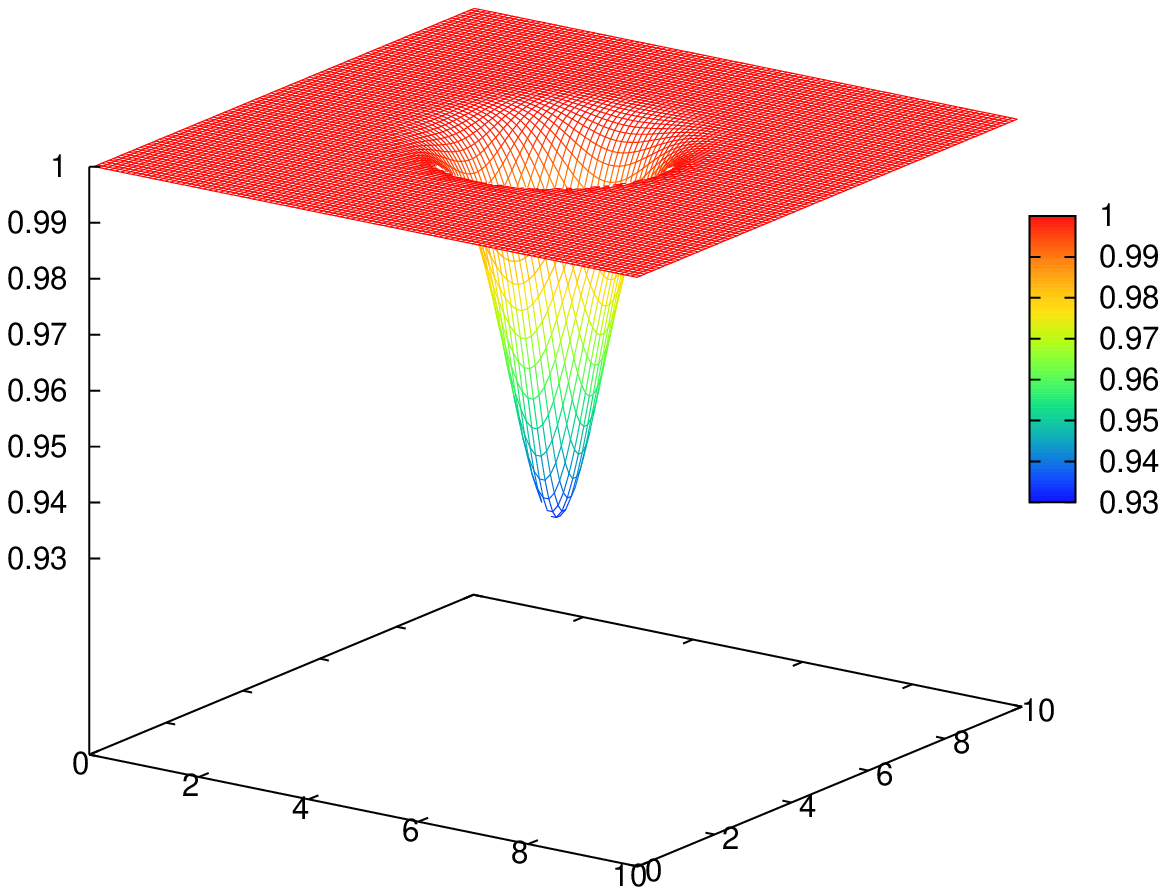} 
    \includegraphics[width=0.33\textwidth]{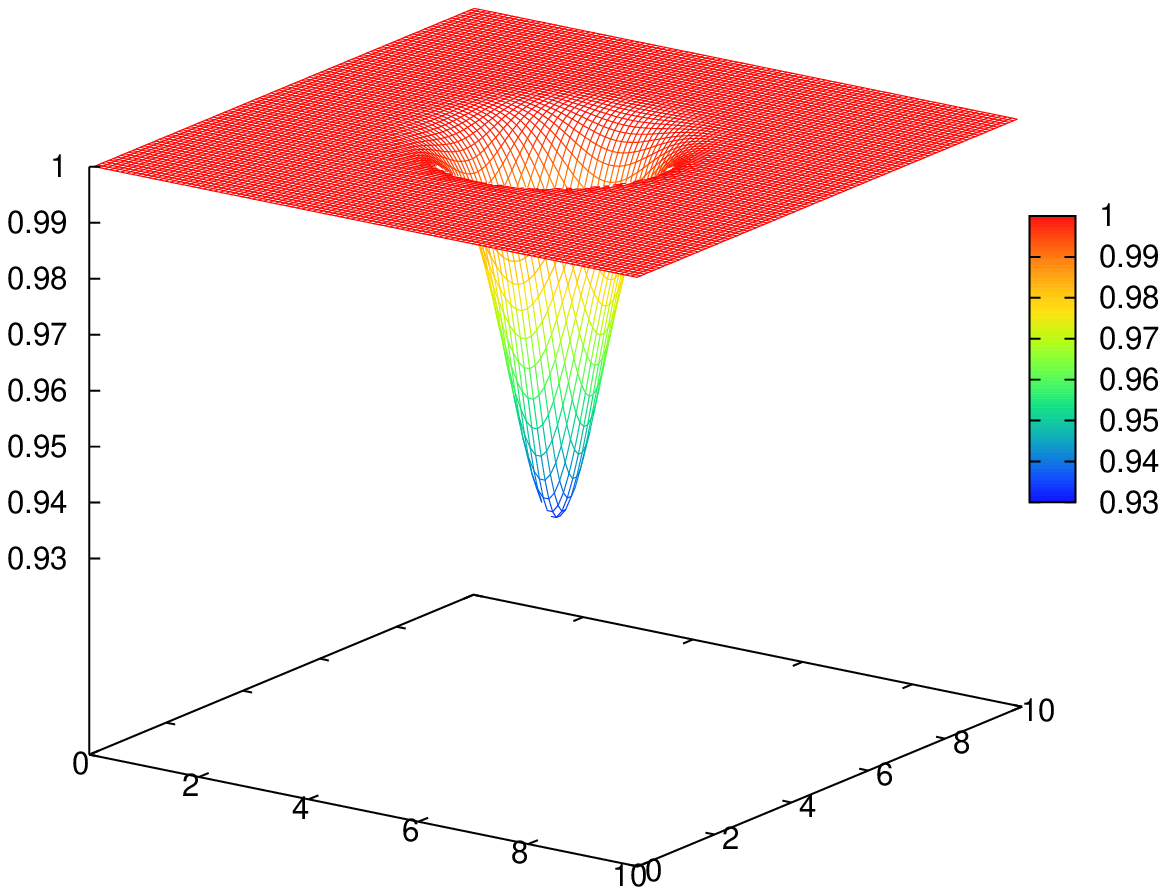} 
    \includegraphics[width=0.36\textwidth]{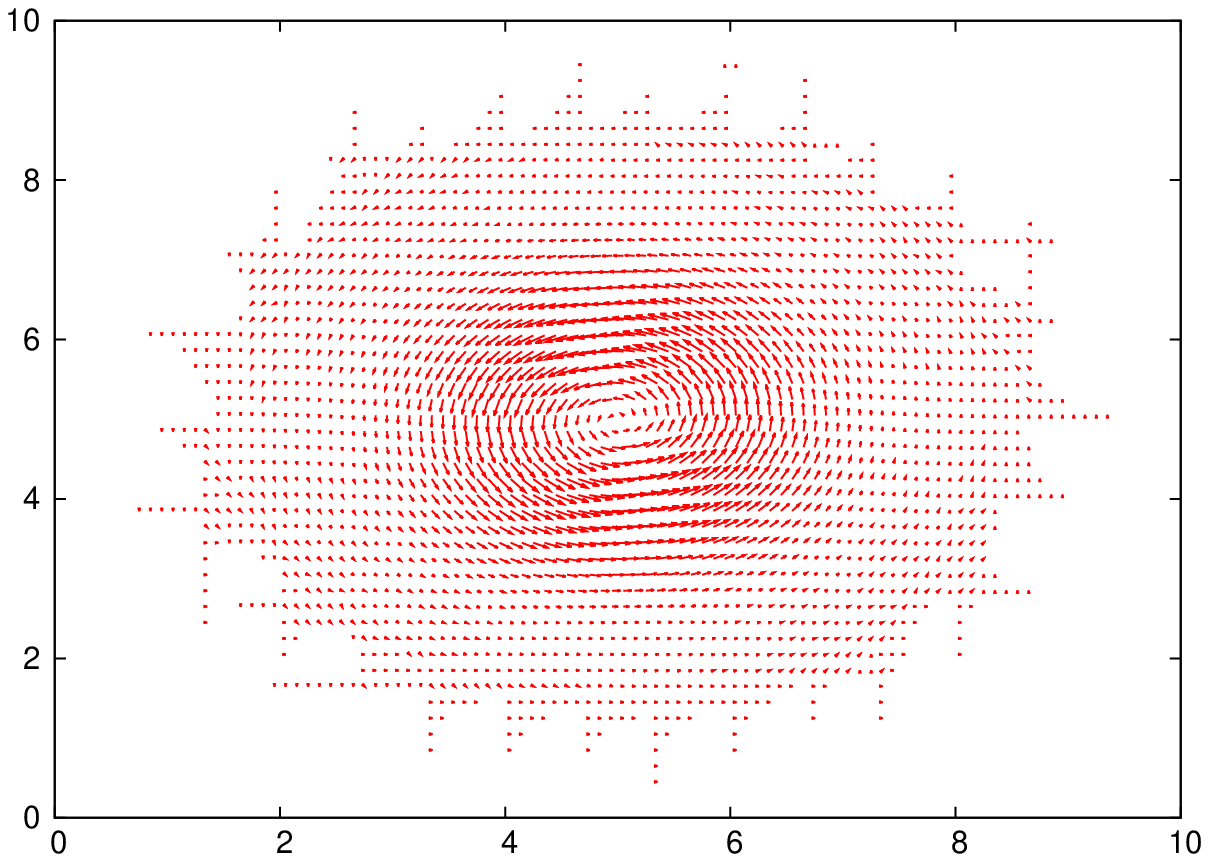} 
  \end{center}
  \caption{Test 3.1. Two dimensional vortex test case for $\tau=10^{-1}$ 
    with $M=100\times 100$ spatial cells and $N=32^2$ velocity cells.    
    Initial density, temperature 
    and velocity vector $(\delta u_x, \delta u_y)$.
   }
  \label{fig:test6_0}
\end{figure}

The test case consists of an isentropic vortex in motion initially introduced for the collisional regime, i.e. the compressible Euler equations, in two dimensions in \cite{Shu1}. 
This problem has an exact smooth solution expressed analytically in the fluid regime. The computational domain is $\Omega = [0,10]^2$. The ambient flow is characterized by a gas with density, 
mean velocity and temperature respectively of $\rho_\infty=1.0$, $u_{x,\infty}=1.0$, $u_{y,\infty}=1.0$, $T_\infty=1.0$. A vortex is centered at $(x_{\text{v}}, y_{\text{v}})=(5,5)$ and supplemented at the 
initial time $t=0$ with conditions
$u_x(t=0) = u_{x,\infty} + \delta u$, $ u_y(t=0) = u_{y,\infty} + \delta v$,  $ T(t=0) = T_\infty + \delta T$ with
\begin{eqnarray}
  \nonumber
  \delta u_x   = -y' {\frac {\beta} {2 \pi}} \exp \left( {\frac {1-r^2} {2}} \right),  \quad
  \delta u_y   = x' {\frac {\beta} {2 \pi}} \exp \left( {\frac {1-r^2} {2}} \right), \quad
  \delta T = - { \frac {(\gamma - 1 ) \beta} {8 \gamma \pi^2}} \exp \left( {1-r^2} \right),
\end{eqnarray}
where $r$ is the Euclidean distance in the two dimensional space, i.e. $ r = \sqrt{{x'}^2 + {y'}^2}$, and $x' = x - x_{\text{v}}, y' = y - y_{\text{v}}$.
The vortex strength depends on the value $\beta$ fixed here to $5.0$.
The initial density is given by
\begin{eqnarray}
\rho(t=0)
       =   \rho_\infty \left( {\frac{T(t=0)}{T_\infty} } \right)^{\frac{1}{\gamma-1} }
       = \left(1 - { \frac {(\gamma - 1 ) \beta} {8 \gamma \pi^2}} \exp \left( {1-r^2} \right)
            \right)^{\frac{1}{\gamma-1} }.
\end{eqnarray}
Periodic boundary conditions are prescribed everywhere. At the final time chosen: $t_{\text{final}}=10$, the vortex is back to its original position and, in the collisional
regime, the initial and final conditions are alike. In a rarefied regime, the exact solution is not known, but at least the cylindrical symmetry of the problem must be retrieved. 
$M=100\times 100$ uniform spatial cells are considered on domain $[0;10]^2$  with $N=32\times 32$ uniform velocity cells on a velocity domain $[-7.5;7.5]^2$. 
The relaxation frequency is fixed to $\tau=10^{-1}$. The initial data are presented in Figure~\ref{fig:test6_0}.

\begin{figure}[ht]
  \begin{tabular}{cccc}
    \hspace{-1.cm}
    & \textbf{Density} & \textbf{Temperature} & \textbf{Velocity} 
    \\
    \hspace{0.cm}
    \rotatebox{90}{\hspace{0.5cm} \textbf{BGK results}}
    & 
    \hspace{-0.4cm}
    \includegraphics[width=0.32\textwidth]{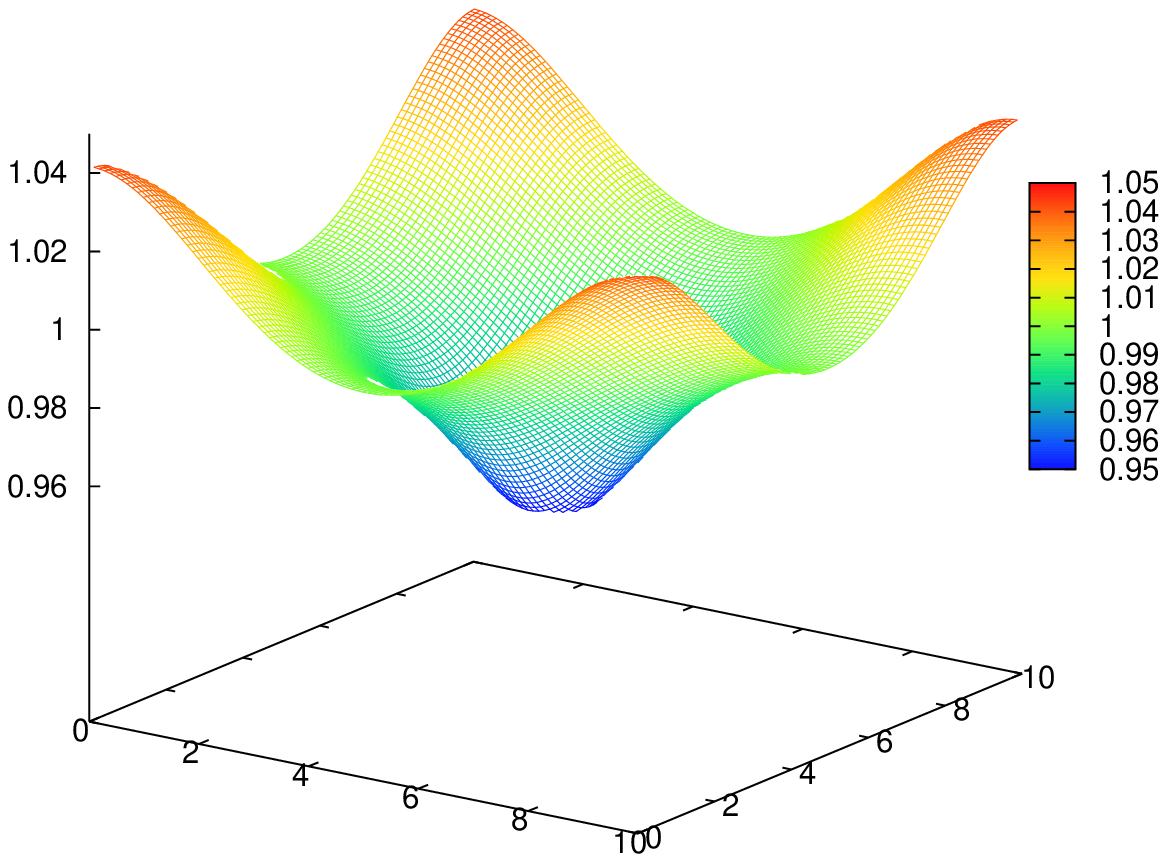} &
    \hspace{-0.5cm}
    \includegraphics[width=0.32\textwidth]{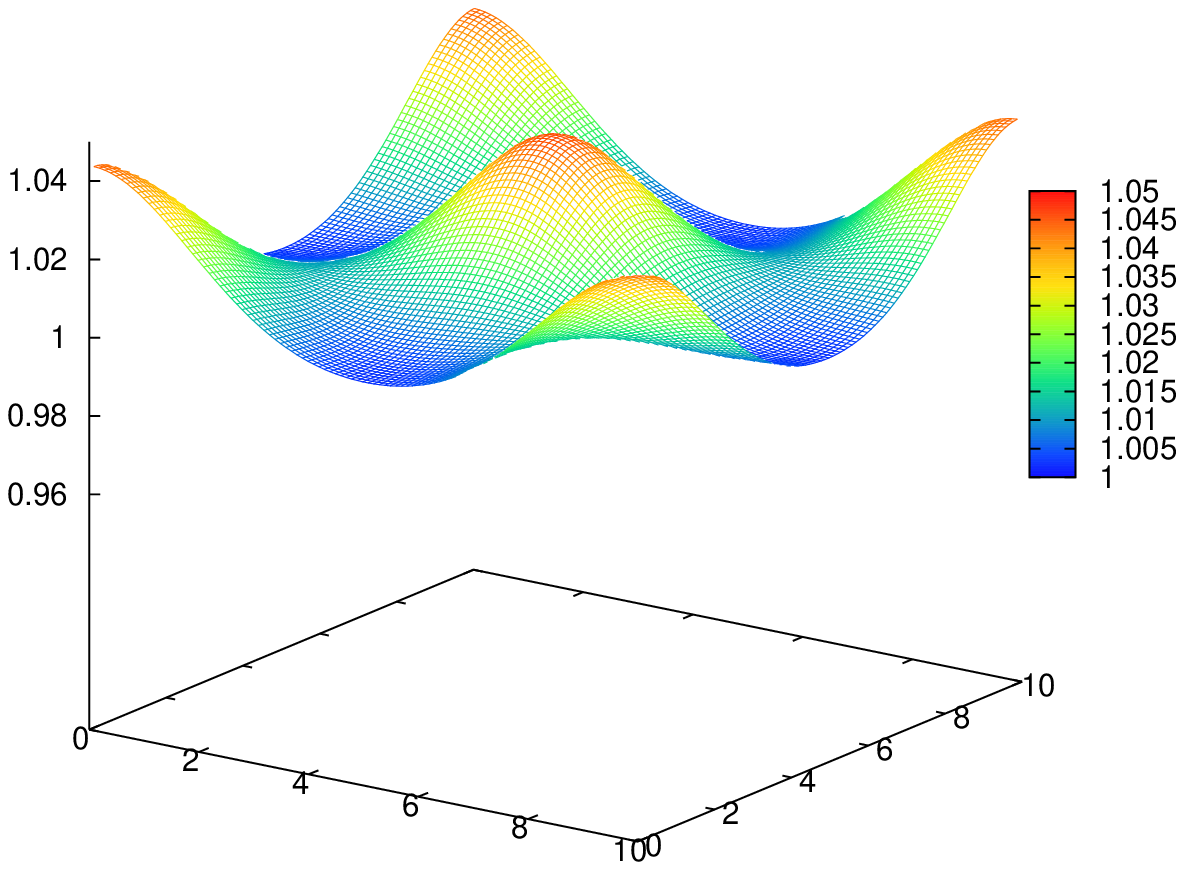} &
    \hspace{-0.5cm}
    \includegraphics[width=0.32\textwidth]{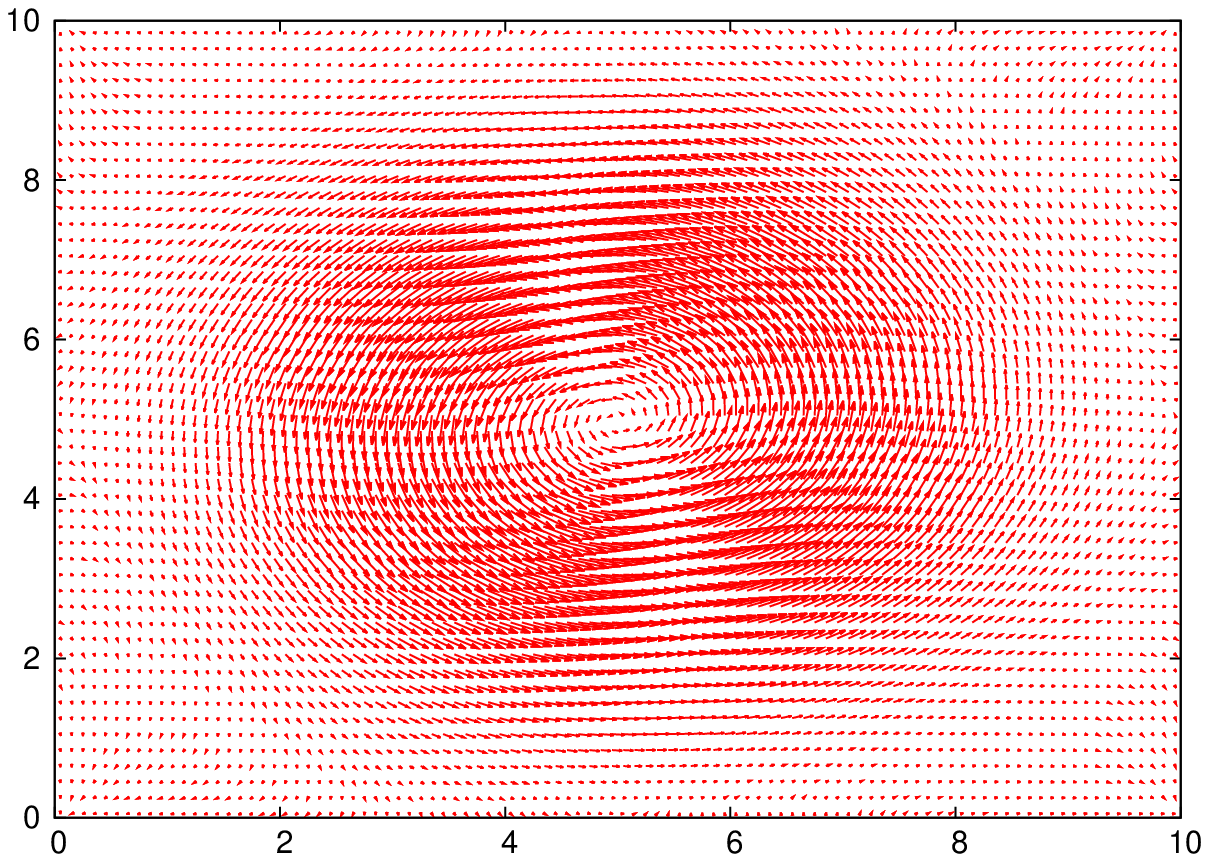} \\    
    \hspace{0.cm}
    \rotatebox{90}{\hspace{0.32cm} \textbf{Boltzmann results}}
    & 
    \hspace{-0.4cm}    
    \includegraphics[width=0.32\textwidth]{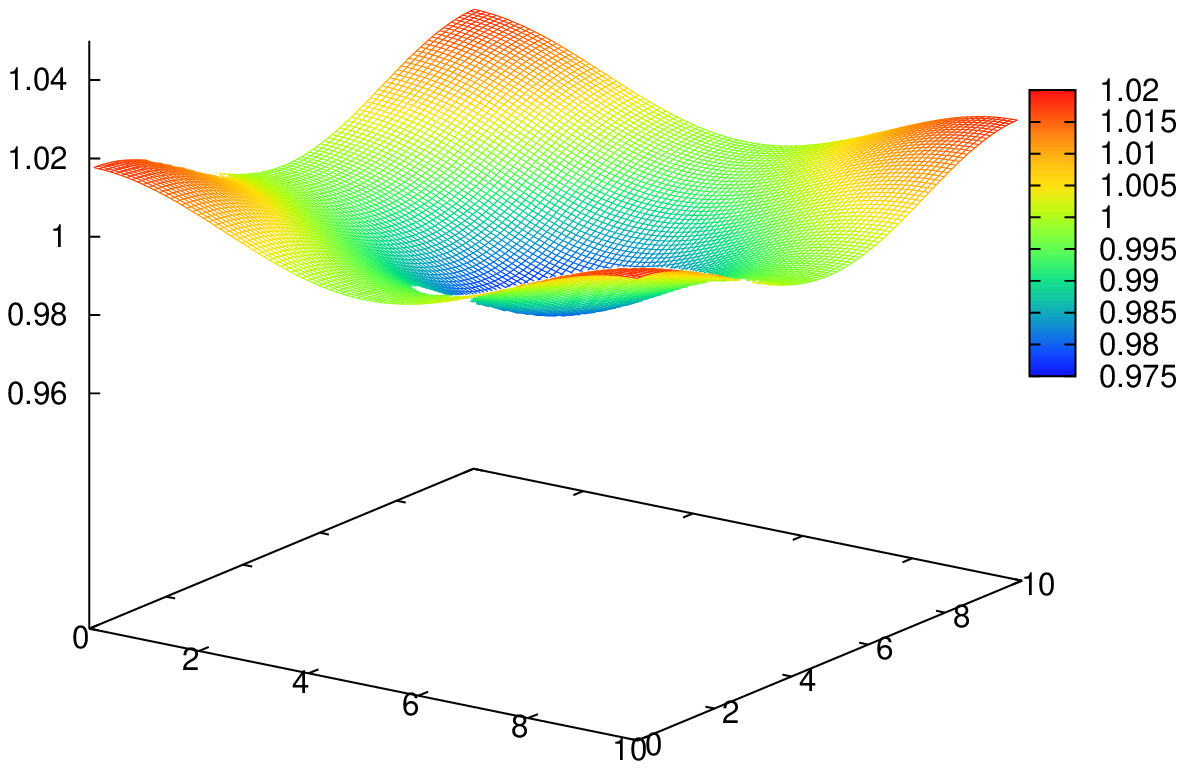} &
    \hspace{-0.5cm}
    \includegraphics[width=0.32\textwidth]{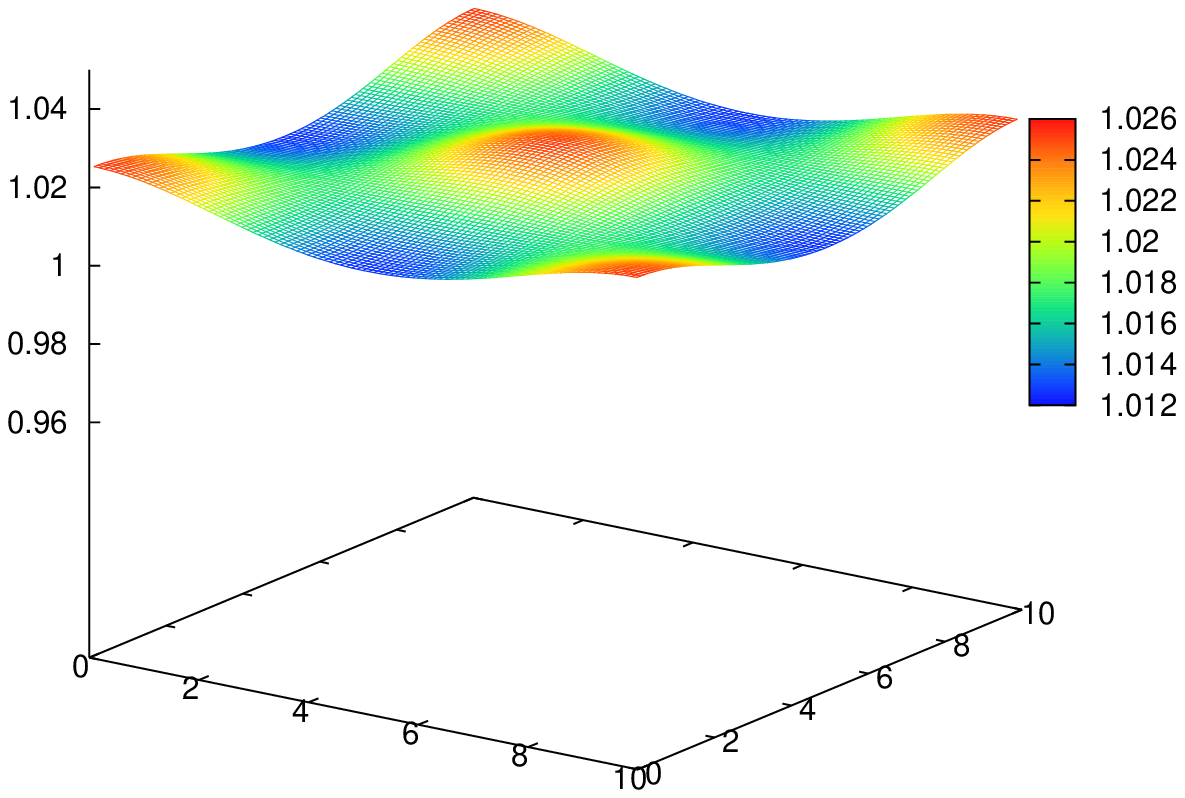} &
    \hspace{-0.5cm}
    \includegraphics[width=0.32\textwidth]{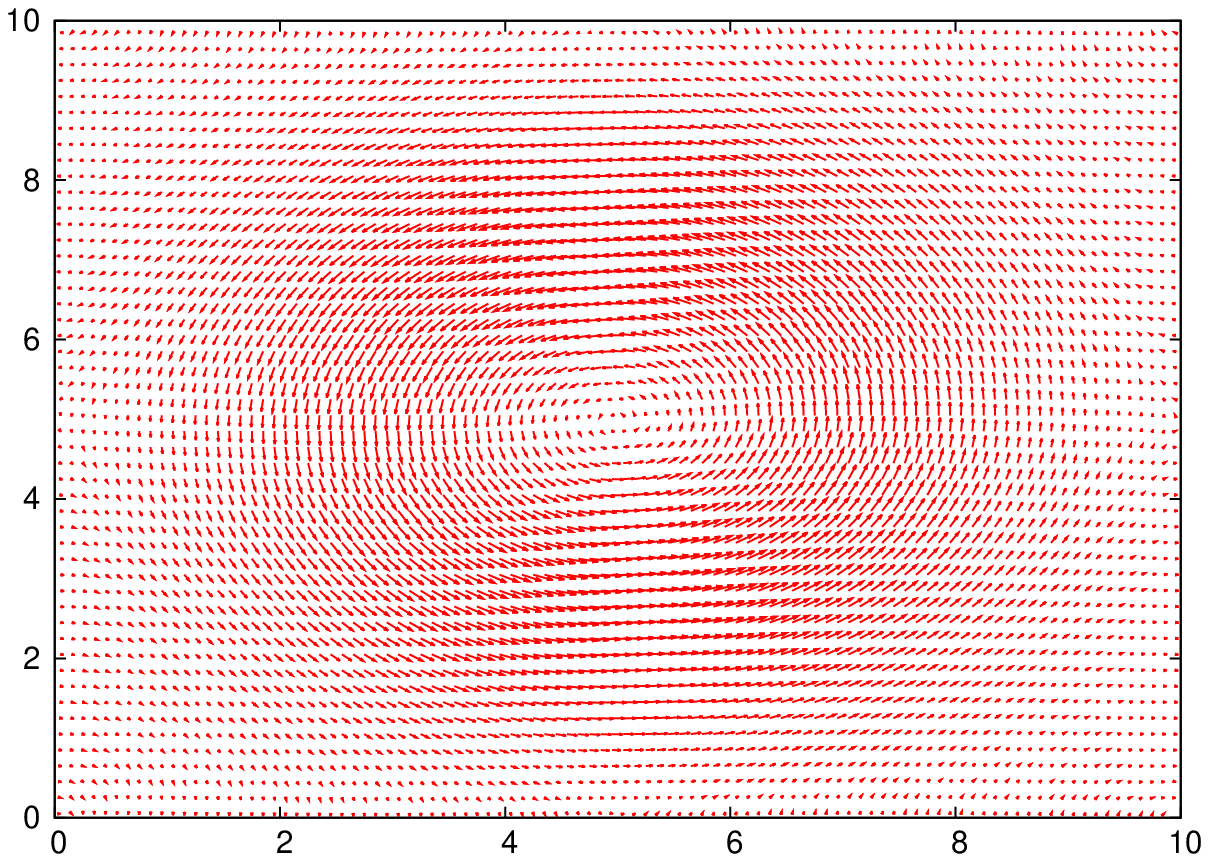}
  \end{tabular}
  \caption{Test 3.1. Two dimensional vortex test case for $\tau=10^{-1}$ 
    with $M=100\times 100$ spatial cells and $N=32^2$ velocity cells.   
    Top/bottom line: results for BGK/Boltzmann models at $t_{\text{final}} =10.0$.
    Left to right: density, temperature, velocity vector $(\delta u_x, \delta u_y)$.
   }
  \label{fig:test6_1}
\end{figure}

\begin{figure}[ht]
  \begin{center}
    \begin{tabular}{ccc} 
      \includegraphics[width=0.27\textwidth]{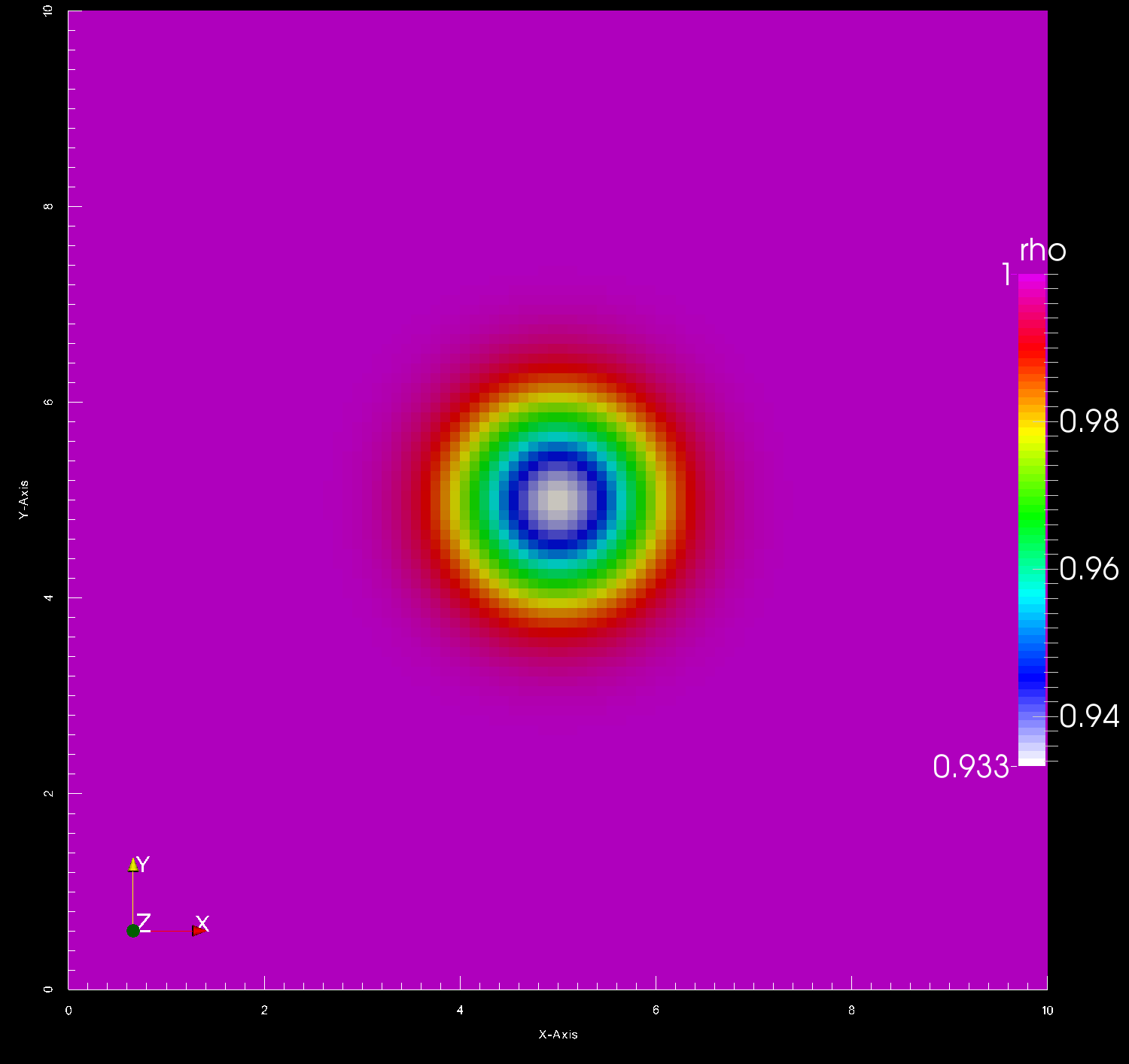} &
      \includegraphics[width=0.27\textwidth]{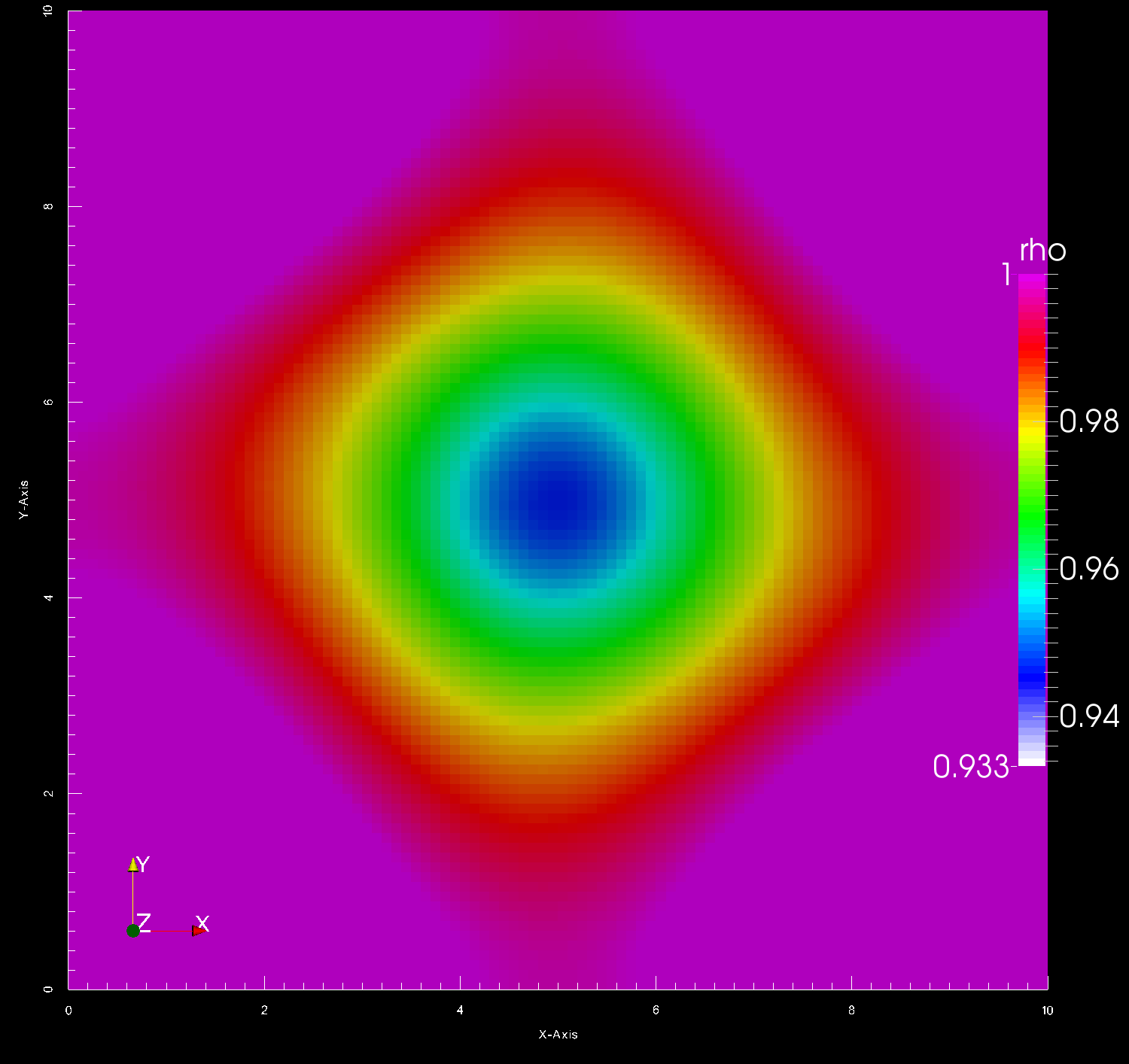} &
      \includegraphics[width=0.27\textwidth]{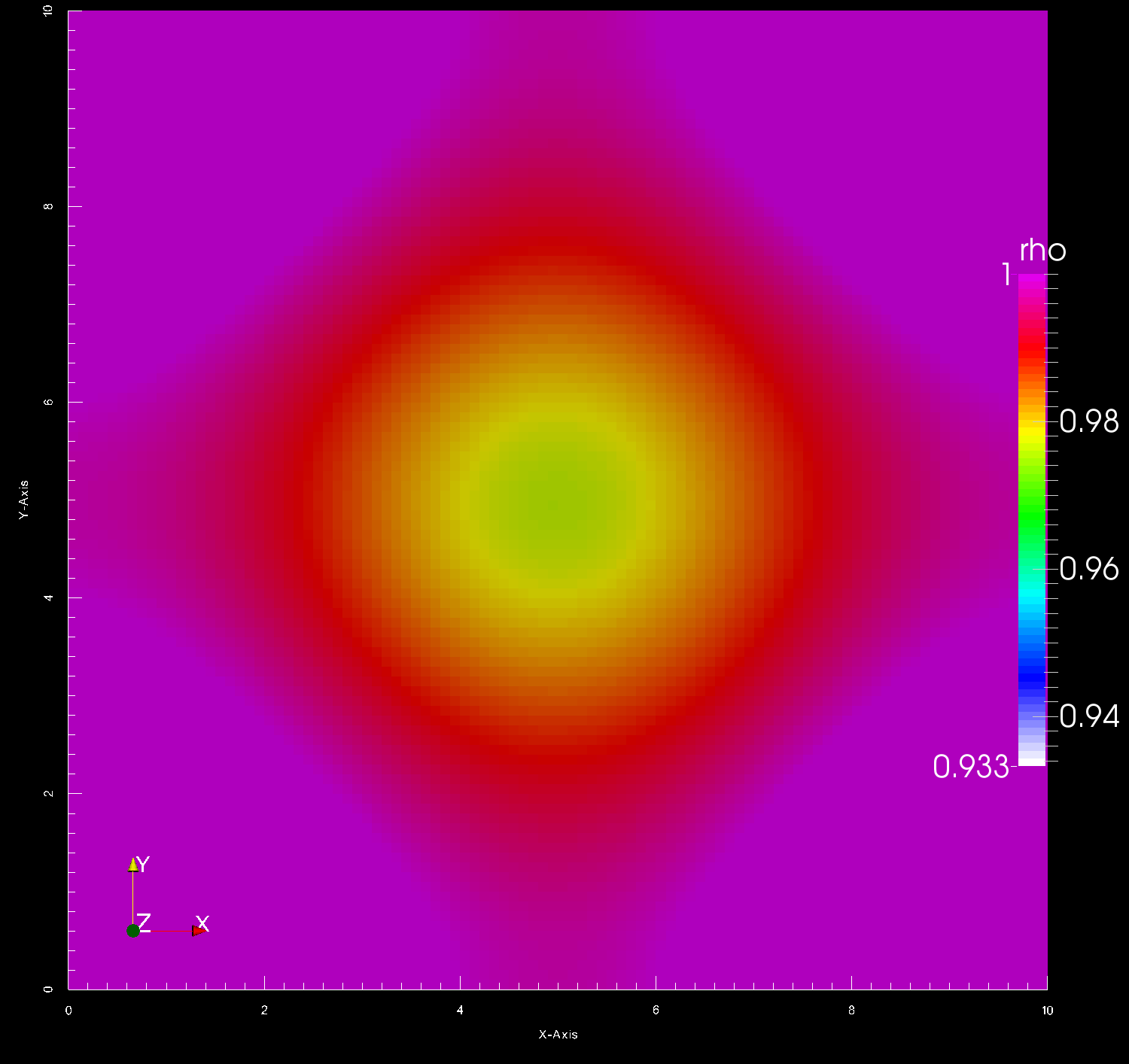} \\
      Initial density & BGK model results & Boltzmann model results 
    \end{tabular}
    \caption{Test 3.1. Two dimensional vortex test case for $\tau=10^{-1}$ 
    with $M=100\times 100$ spatial cells and $N=32^2$ velocity cells.   
     Left: initial density.
      Middle/right: results for BGK/Boltzmann model at $t_{\text{final}} =10.0$ using the 
      same color scale of the initial density.
    }
  \label{fig:test6_2}
  \end{center}
\end{figure}

This problem is simulated using both the BGK and Boltzmann models with $\nu=\rho$ for the BGK case. 
The results are presented in Figure~\ref{fig:test6_1} where density, temperature and 
the velocity are plotted. The top line presents BGK model results while
the bottom line presents Boltzmann ones. The same scale is used to ease the comparison.
The velocity fluctuation $(\delta u_x, \delta u_y)$ is represented with the same scaling as to observe that the vortex
is more dissipated by the Boltzmann model than by BGK one.
In Figure~\ref{fig:test6_2} we replot the initial density along with the final 
BGK and Boltzmann results using the same color scale used for showing the initial data.
From the Figures we can observe that both models reproduce a vortex at the correct
final location but with different rates of dissipation. The BGK model furnishes results closer to the compressible Euler solution results which means over-relaxation.

We analyze now the performances of the scheme. In order to obtain the previous results both models used $1400$ time steps to reach $t_{\text{final}}$.
This result is due to the choice of the time integrator: the same for both models. The Open MP version of the code is run on $8$ cores and the total CPU times 
on a laptop HP ZBook Intel Core i7-4940MX CPU @ 3.10GHz$\times$8 operated by Ubuntu 15.10 64bits,
are of the order $\text{CPU}_{\text{BGK}}\simeq 102\text{s}$, and $\text{CPU}_{\text{Boltz}}\simeq 2008\text{s}$. We split the code in several distinct and conceptually important parts which
are denoted as \texttt{ Transport}, \texttt{ ToConservative}, \texttt{ ToPrimitive}, 
and \texttt{ Collision}. The first routine implements the transport phase, the last routine implement the collision phase, while the two routines in the middle reconstruct the conserved macroscopic
and primitive variables from the kinetic distribution. These two routines are necessary for defining the Maxwellian distribution for performing the collision step in the BGK model, while
they only serve to show the results for the Boltzmann model. Moreover, since the macroscopic quantities change only due to the transport phase, one can think to associate the cost of these
routines to the transport part. All routines are monitored during a simulation of the vortex problem using the same mesh reported above. The results are reported in Table~\ref{tab:profiling2Dx2D}.
The Boltzmann model is $20$ times more expensive than BGK one. More in details, the relative cost
of the collision steps jumps from $70\%$ for BGK model to $97\%$ for Boltzmann model.
Indeed, while for the BGK model, the collision routine and the routine recomputing the conservative variables
have some impact on the total CPU time, for the Boltzmann model the collision part is the only one participating to the global cost at least for the mesh used.
\begin{table}
  \begin{center}
  \numerikNine
  \begin{tabular}{|c|c|c||c||cc|}
    \hline
    &  \textbf{Cycle} &  \textbf{CPU}  &  \textbf{Main routines} 
    &  \textbf{Cost CPU }  &  \textbf{Cost}  \\
    & & (s) & & vs total (s) & vs total (\%)  \\
    \hline
    \hline
     \multirow{5}{*}{\begin{sideways} \textbf{BGK}  \end{sideways}}
    & \multirow{5}{*}{ 1400 }
    & \multirow{5}{*}{ 102 }
    &Transport          &   0.36 & 0.35\%  \\
    \cline{4-6}
    & & &ToConservative &  30.26 & 29.77\%  \\
    \cline{4-6}
    & & &ToPrimitive    &  0.19  & 0.19\%  \\
    \cline{4-6}
    & & &Collision     &  70.83 & 69.69\% \\  
    \cline{4-6}
    \cline{4-6}
    & & & = & 101.64 & 100\%  \\
    \hline
    \hline
    \multirow{5}{*}{\begin{sideways} \textbf{Boltzmann}  \end{sideways}}
    & \multirow{5}{*}{ 1400 }
    & \multirow{5}{*}{ 2008 }
    &Transport          &    1.06 &  0.053\%  \\
    \cline{4-6}
    & & &ToConservative &   66.75 &  3.325\%  \\
    \cline{4-6}
    & & &ToPrimitive    &    0.29 &  0.015\%  \\
    \cline{4-6}
    & & &Collision     & 1939.63 & 96.61\%  \\  
    \cline{4-6}
    \cline{4-6}
    & & & = & 2007.73 & 100\%  \\
    \hline
  \end{tabular}
  \caption{ \label{tab:profiling2Dx2D} Profiling of the average cost
  for each routine of the code on the test 3.1. (two dimensional vortex problem)
  simulated on a laptop using $M=100^2$ and $N=32^2$ points.
  }
  \end{center}
\end{table}
We finally gather in Table~\ref{tab:perf2Dx2D} the results in terms of computational costs when
the number of spatial cell in each direction doubles from $M=25^2$ to $M=200^2$ keeping the number 
of velocity cell fixed to $N=32\times 32$ on velocity domain $[-7.5;7.5]^2$. We monitor the number of cycles, the CPU time $\text{T}$ in second(s), and
the time per cycle is obtained by $T_{\text{cycle}} = \text{T}/N_{\text{cycle}}$, 
the time per cycle per cell by $T_{\text{cell}} = \text{T}/N_{\text{cycle}}/M$ and 
the time per cycle per degree of freedom $T_{\text{dof}} = \text{T}/N_{\text{cycle}}/(M^d\times N^{d_v})$.
Doubling the number of spatial cells leads to an increase in CPU time by a factor $7$ to $8$ for 
both models and the number of time cycles is exactly doubled.
Boltzmann model demands on average $20$ times more CPU resources than the BGK one. 
The time per cycle per cell  $T_{\text{cell}}$ is rather constant ($7.5\times 10^{-6}$ for BGK 
and $1.4\times 10^{-4}$ for Boltzmann (thanks to the almost linear complexity of the fast spectral solver), it becomes relatively easy to estimate the cost of further 
refined simulations keeping the velocity mesh fixed.

\begin{table}
\begin{center}
  \begin{tabular}{|c|cc|c||c|c|c|c|c|}
    \hline
    \multicolumn{9}{|c|}{\textbf{Vortex problem in two dimension}} \\
    \hline
    \hline
    \multirow{3}{*}{\begin{sideways}\textbf{Model}\end{sideways}}
    & \multicolumn{2}{|c|}{\textbf{Velocity}}  
      &\textbf{Cell \#}    & \textbf{Cycle} & \textbf{Time} & T\textbf{/cycle} & T\textbf{/cell} &  T\textbf{/d.o.f} \\
    & $N$ & {\begin{sideways} Vel.\end{sideways}}&$M^d \times N^{d_v}$ &  $N_{\text{cycle}}$  &  $\text{T}$ (s)  & $T_{\text{cycle}}$ (s) & $T_{\text{cell}}$ (s) & $T_{\text{dof}}$ (s) \\
    \hline
    \hline
    \hline
    \multirow{9}{*}{\begin{sideways}  \textbf{BGK} \end{sideways}}  
    & \multirow{9}{*}{ $32^2$} &
    \multirow{9}{*}{\begin{sideways} $[-7.5,7.5]$  \end{sideways}}
    &$25^2\times 32^2$  & \multirow{2}{*}{$351$}     & \multirow{2}{*}{$2.61$}  & \multirow{2}{*}{$0.0035$} & \multirow{2}{*}{$5.60\times 10^{-6}$} & \multirow{2}{*}{$5.47\times 10^{-9}$}\\
    &&&$=64 \times 10^4$ & &  & & & \\
    \cline{4-9}
    &&&$50^2 \times 32^2$  & \multirow{2}{*}{$701$}  & \multirow{2}{*}{$13.77$} & \multirow{2}{*}{$0.0196$} & \multirow{2}{*}{$7.84\times 10^{-6}$}  & \multirow{2}{*}{$7.66\times 10^{-9}$} \\
    &&&$=256\times 10^{4}$ & & & & &  \\
    \cline{4-9}
    &&&$100^2 \times 32^2$  & \multirow{2}{*}{$1400$}& \multirow{2}{*}{$102.42$}& \multirow{2}{*}{$0.0732$} & \multirow{2}{*}{$7.32\times 10^{-6}$} & \multirow{2}{*}{$7.15\times 10^{-9}$} \\
    &&&$=1024\times 10^{4}$ & & & & & \\
    \cline{4-9}
    &&&$200^2 \times 32^2$ & \multirow{2}{*}{$2800$} & \multirow{2}{*}{$785.54$}& \multirow{2}{*}{$0.2806$} & \multirow{2}{*}{$7.02\times 10^{-6}$} & \multirow{2}{*}{$6.85\times 10^{-9}$} \\
    &&&$=4096\times 10^4$ & &  &  & & \\
    \hline
    \hline
    \hline
    \multirow{9}{*}{\begin{sideways}  \textbf{Boltzmann} \end{sideways}}  
    & \multirow{9}{*}{ $32^2$} &
    \multirow{9}{*}{\begin{sideways} $[-7.5,7.5]$  \end{sideways}}
    &$25^2\times 32^2$  & \multirow{2}{*}{$351$}     & \multirow{2}{*}{$32.92$}  & \multirow{2}{*}{$0.0938$} & \multirow{2}{*}{$1.50\times 10^{-4}$} & \multirow{2}{*}{$1.47\times 10^{-7}$}\\
    &&&$=64 \times 10^4$ & &  & & & \\
    \cline{4-9}
    &&&$50^2 \times 32^2$  & \multirow{2}{*}{$701$}  & \multirow{2}{*}{$245.03$} & \multirow{2}{*}{$0.350$} & \multirow{2}{*}{$1.40\times 10^{-4}$} & \multirow{2}{*}{$1.37\times 10^{-7}$} \\
    &&&$=256\times 10^{4}$ & & & & & \\
    \cline{4-9}
    &&&$100^2 \times 32^2$  & \multirow{2}{*}{$1400$}& \multirow{2}{*}{$2008.56$}& \multirow{2}{*}{$1.435$} & \multirow{2}{*}{$1.44\times 10^{-4}$}& \multirow{2}{*}{$1.40\times 10^{-7}$} \\
    &&&$=1024\times 10^{4}$ & & & & & \\
    \cline{4-9}
    &&&$200^2 \times 32^2$ & \multirow{2}{*}{$2800$} & \multirow{2}{*}{$15762$}& \multirow{2}{*}{$5.630$} & \multirow{2}{*}{$1.41\times 10^{-4}$} & \multirow{2}{*}{$1.37\times 10^{-7}$}\\
    &&&$=4096\times 10^4$ & &  &  & & \\
    \hline
  \end{tabular}
  \caption{ \label{tab:perf2Dx2D}
    Test 3.1. Two dimensional vortex problem simulations 
    with spatial mesh variation.
    Monitoring of CPU time.
    Time per cycle is obtained by $T_{\text{cycle}} = \text{T}/N_{\text{cycle}}$, 
    time per cycle per cell by $T_{\text{cell}} = \text{T}/N_{\text{cycle}}/M$ and 
    time per cycle per degree of freedom $T_{\text{dof}} = \text{T}/N_{\text{cycle}}/(M^d \times N^{d_v})$.
  }
\end{center}
\end{table}

\subsubsection{Test 3.2. Re-entry test in two dimensions with changing angle of attack in time.}
This test is inspired from re-entry test cases described in \cite{capsule,FKS_DD}. The computational domain is set to $\Omega=[0;4]\times[0;4]$. 
Within this domain we initiate three static objects, two small rectangles upfront ($[x_0;x_1]\times[y_0;y_1]$ and $[x_0;x_1]\times[y'_0;y'_1]$) 
and one larger one behind ($[x'_0;x'_1]\times[(y_0+y'_0)/2;(y_1+y'_1)/2]$) where
$x_0=1.5$, $x_1=1.7$, $x'_0=1.8$, $x'_1=2$, $y_0=1.7$, $y_1=1.95$, $y'_0=2.05$, $y'_1=2.3$.
The computational mesh in physical space is made of $800\times 800$ square cells.
The velocity space is $[-10,10]$ and is discretized with $32^2$ points.
We report the solution obtained in the fluid regime, i.e. $\tau=0$, obtained by projecting the distribution $f$ over the equilibrium state $M$ after each transport phase and
the results obtained with $\tau=10^{-2}$ using both BGK and Boltzmann models.
The initial density is set to $\rho(t=0)=1$, the velocity $(u_x,u_y)(t=0)=(3,0)$ and the temperature to $T(t=0)=1$ everywhere.
The final time is set to $t_{\text{final}}=10$. On the boundaries with the objects reflective boundary conditions are employed.
Inflow boundary conditions are imposed to the west boundary whereas outflow boundary conditions
are set elsewhere.
The inflow boundary conditions are evolving in time and equal to 
\begin{equation} \label{eq:BCs}
(u_x,u_y)_{BC}(t) = \left\{
  \begin{array}{lll}
    (3,0) & \text{ if } &   0 \leq t \leq t_1 \\
    (\sqrt{9-g^2(t)}, g(t) ) & \text{ if } & t_1 \leq t \leq t_2 \\
    (\frac{3\sqrt{2}}{2}, \frac{3\sqrt{2}}{2} ) & \text{ if } & t_2 \leq t \leq t_{\text{final}}\\    
  \end{array}
\right.
\end{equation}
where $t_1 = 3/2$, $g(t) = t -t_1$ and $t_2 = 3\sqrt{2}/2 + t_1$.
Given these initial data, we expect a detached shock wave to occur upfront the objects and
some complex wave pattern behind them.
Moreover, setting the inflow boundary conditions to (\ref{eq:BCs}) splits the simulation into three stages. 
The first stage consists in the inflow boundary conditions facing
the objects up to $t_1$. For this stage the upfront detached shock and the complex flow structure
behind the objects are formed but they are not yet steady.
Next, for the second stage, the inflow boundary condition is changing its direction by smoothly increasing the $y$ component
of $(u_x,u_y)_{BC}(t)$ up to $t=t_2$. Note that this mimics a modification of the angle of attack of the objects with time. This change modifies the entire
flow structure. Last, for the third stage the inflow boundary condition is fixed to $(u_x,u_y)_{BC}=(u_0,v_0)$ up to the final time.
As such the flow reaches an almost steady state.

\begin{figure}
  \begin{center}
    \begin{tabular}{cccc} 
      \hspace{-0.5cm}
      \includegraphics[width=0.32\textwidth]{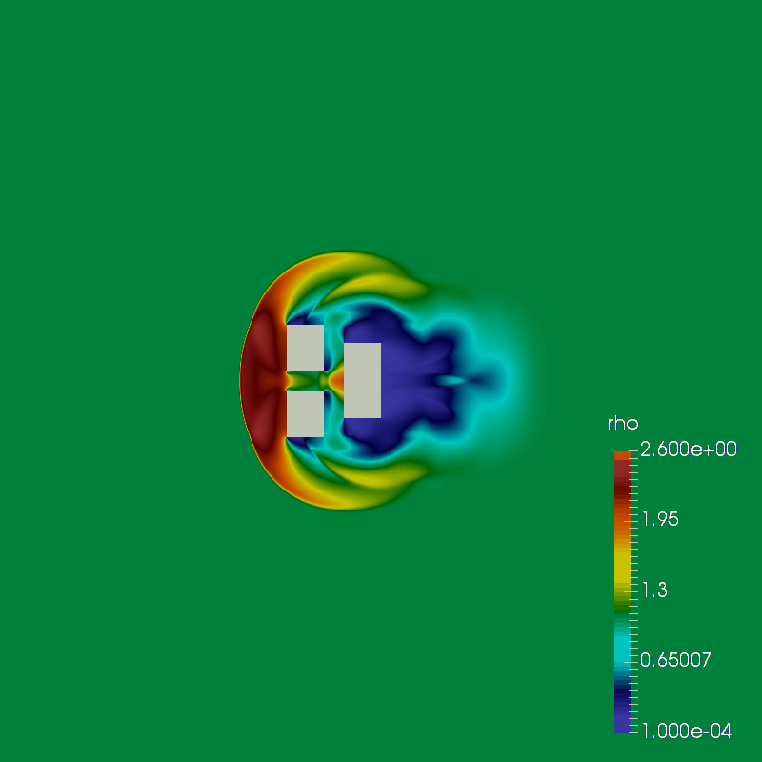} &
      \hspace{-0.5cm}
      \includegraphics[width=0.32\textwidth]{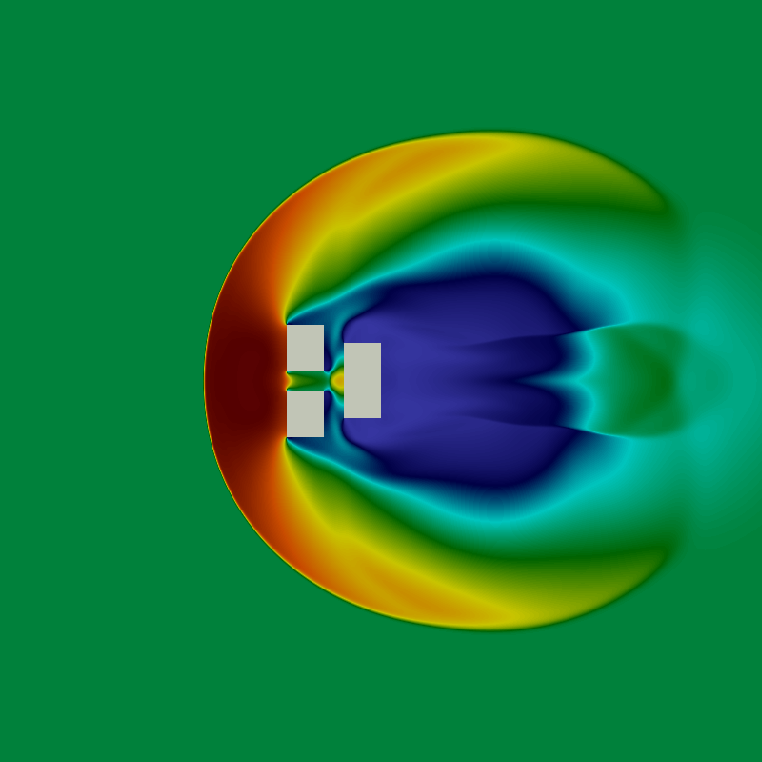} &
      \hspace{-0.5cm}
      \includegraphics[width=0.32\textwidth]{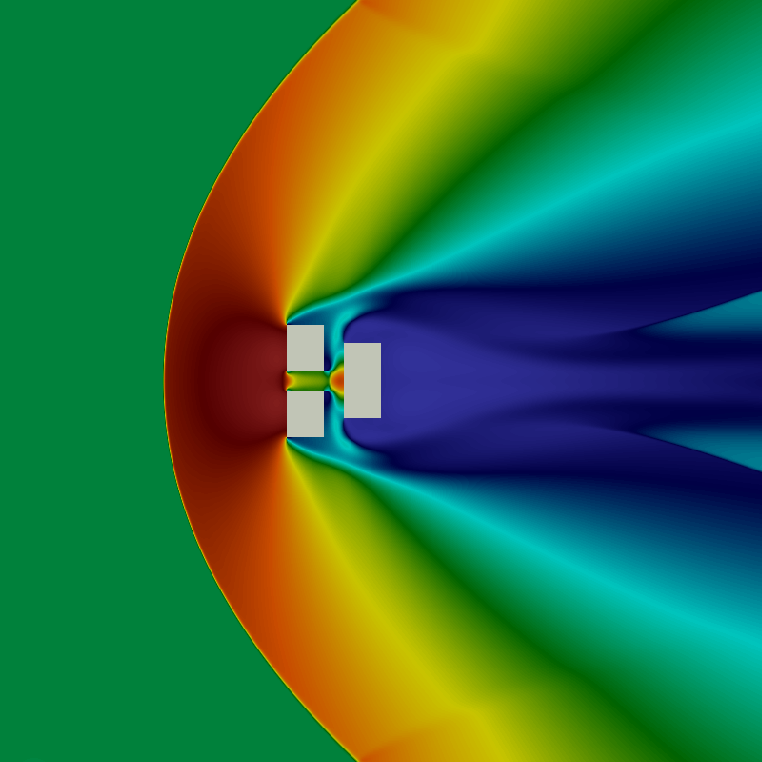} & 
      \\   
      \hspace{-0.5cm}
      \includegraphics[width=0.32\textwidth]{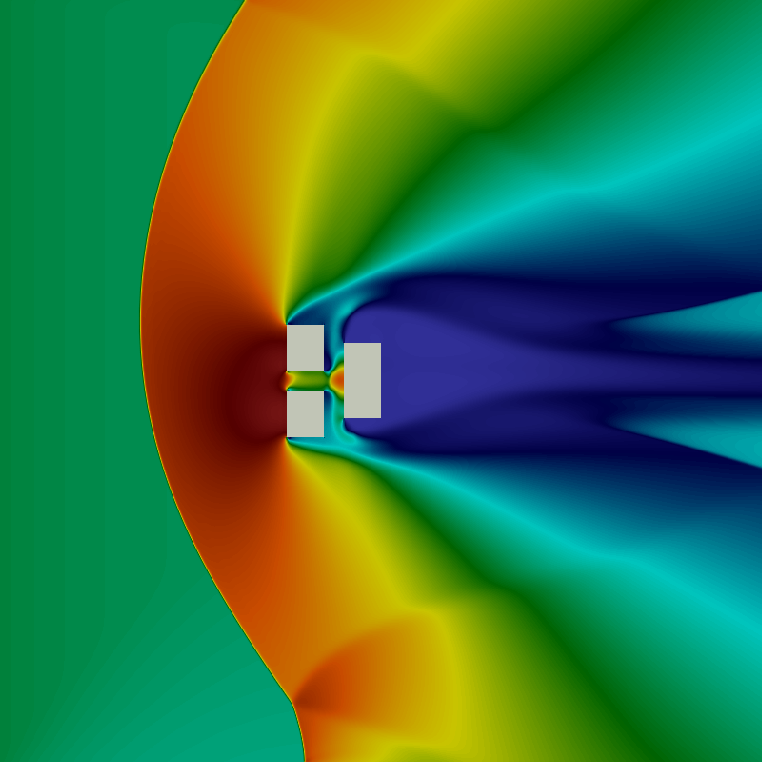} &
      \hspace{-0.5cm}
      \includegraphics[width=0.32\textwidth]{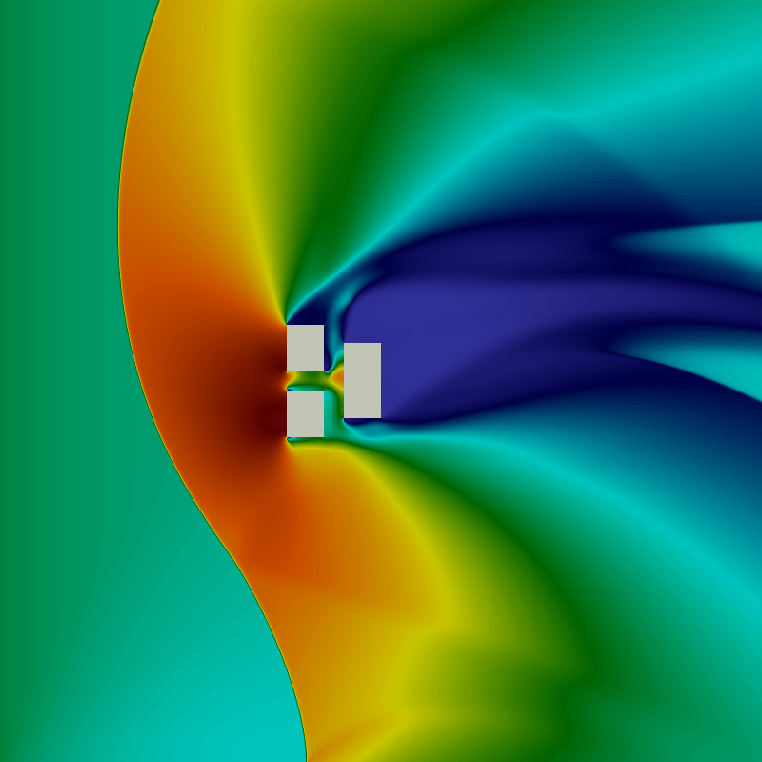} &
      \hspace{-0.5cm}
      \includegraphics[width=0.32\textwidth]{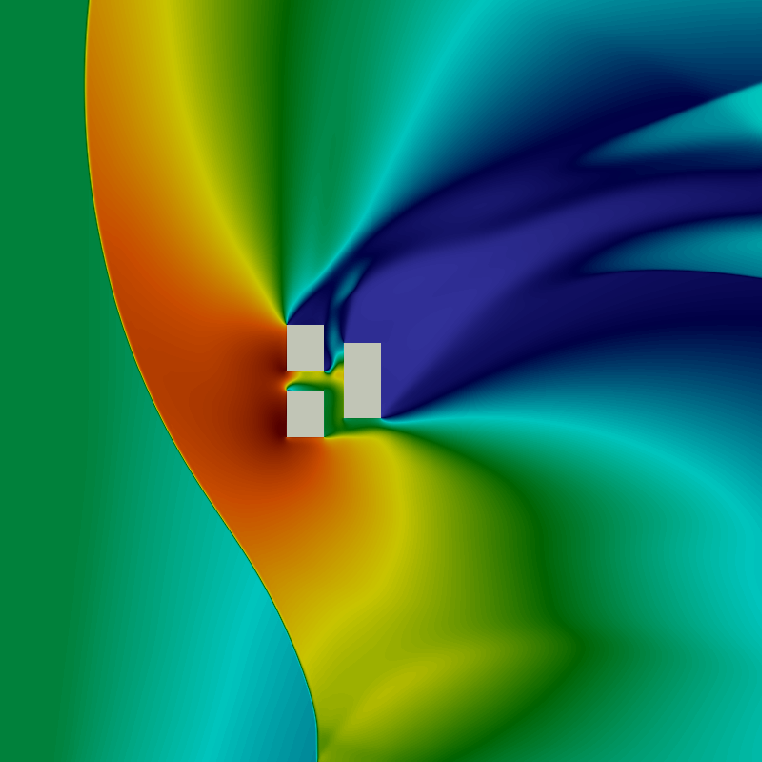} &
      \\   
      \hspace{-0.5cm}
      \includegraphics[width=0.32\textwidth]{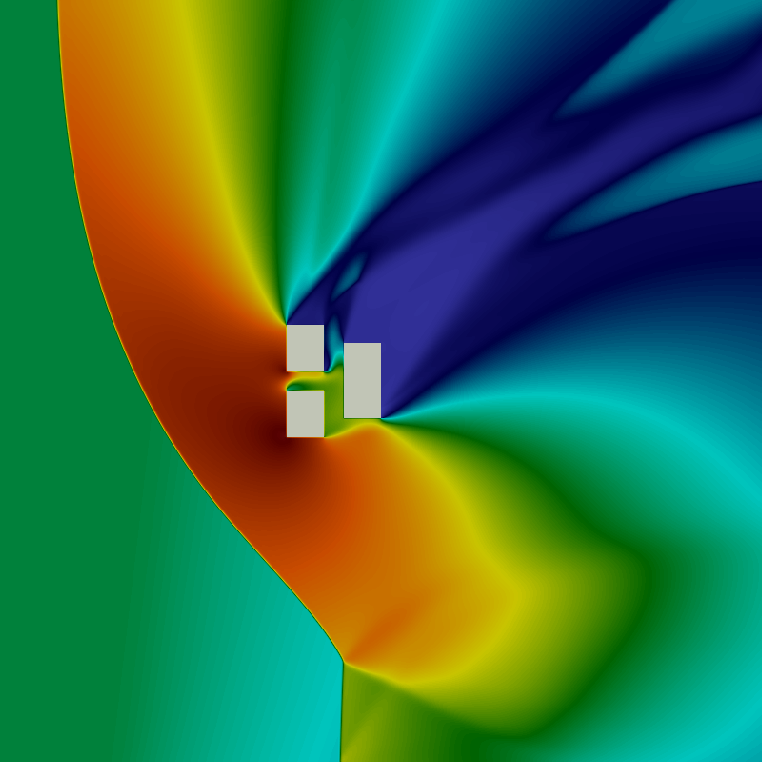} &
      \hspace{-0.5cm}
      \includegraphics[width=0.32\textwidth]{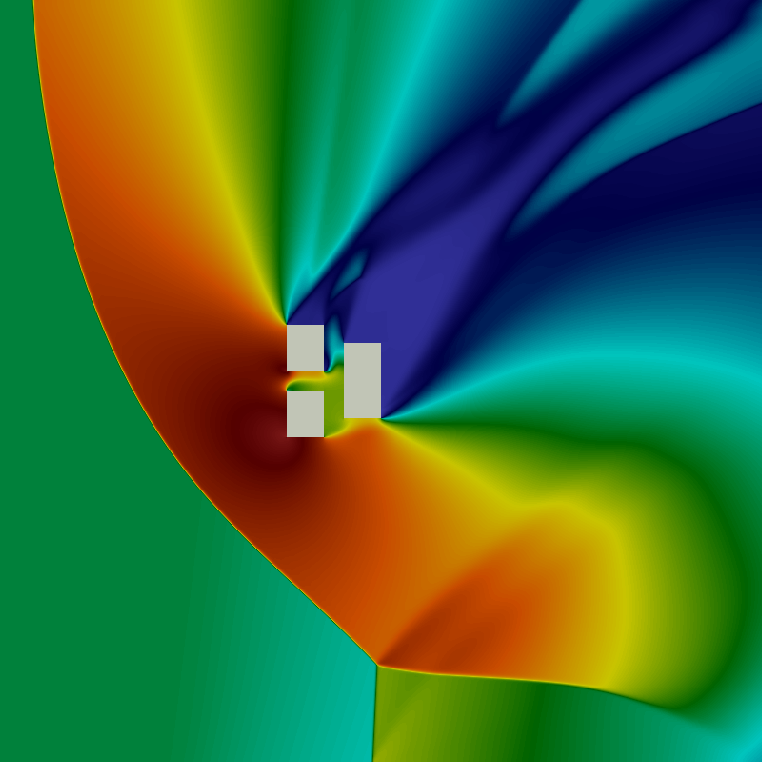} &
      \hspace{-0.5cm}
      \includegraphics[width=0.32\textwidth]{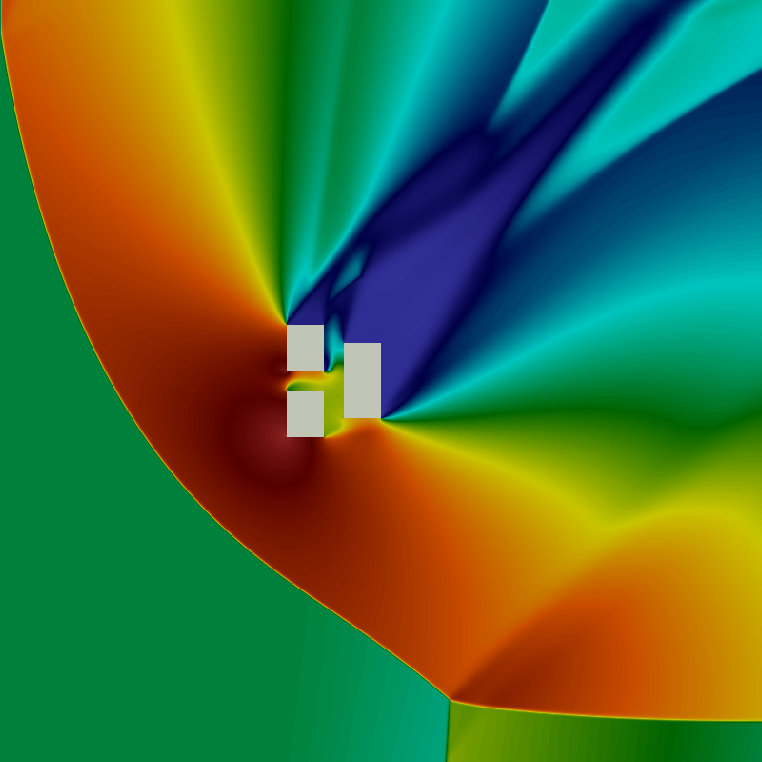} &
    \end{tabular}   
    \caption{
    \label{fig:test7_reentry2Dx2D_BGK0}
    Test 3.2. Two dimensional re-entry test case for $\tau=0$ 
      with $M=800 \times 800$ spatial cells and $N=32^2$ velocity cells. Density profile.      
      Top-left to bottom-right  iterations $500$, $1400$, $4500$, $8000$, $10000$, $12000$, $14000$, $16000$, $20000$.
    }
  \end{center}
\end{figure}

\begin{figure}
  \begin{center}
    \begin{tabular}{ccc} 
      \hspace{-0.5cm}
      \includegraphics[width=0.32\textwidth]{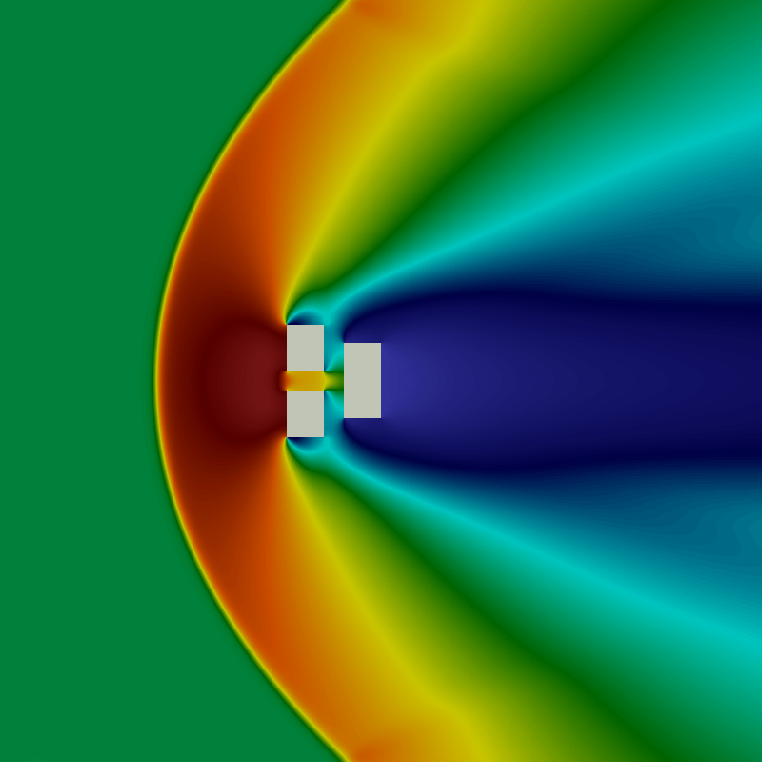} &
      \hspace{-0.5cm}
      \includegraphics[width=0.32\textwidth]{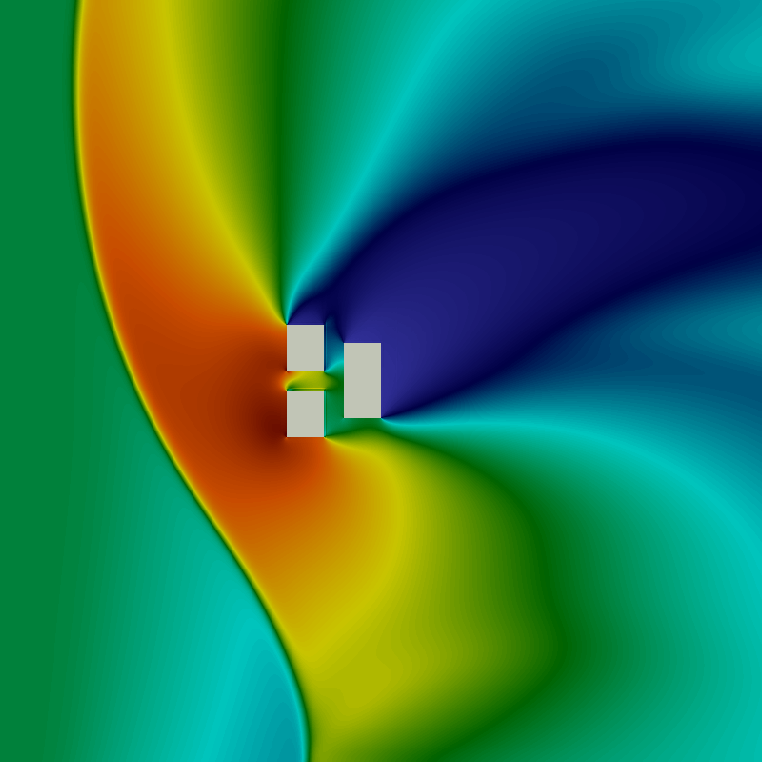} &
      \hspace{-0.5cm}
      \includegraphics[width=0.32\textwidth]{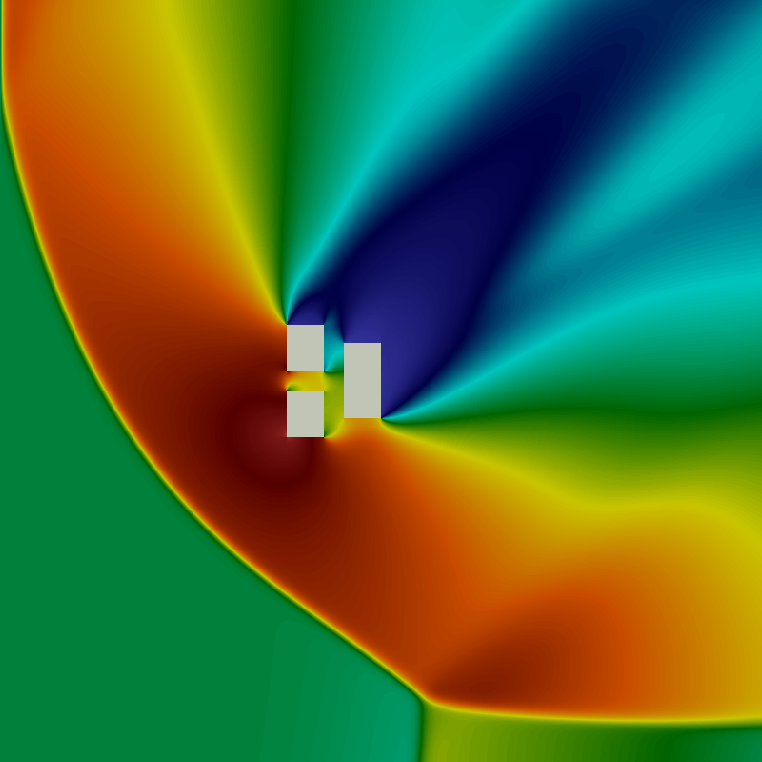} \\
      \hspace{-0.5cm}
      \includegraphics[width=0.32\textwidth]{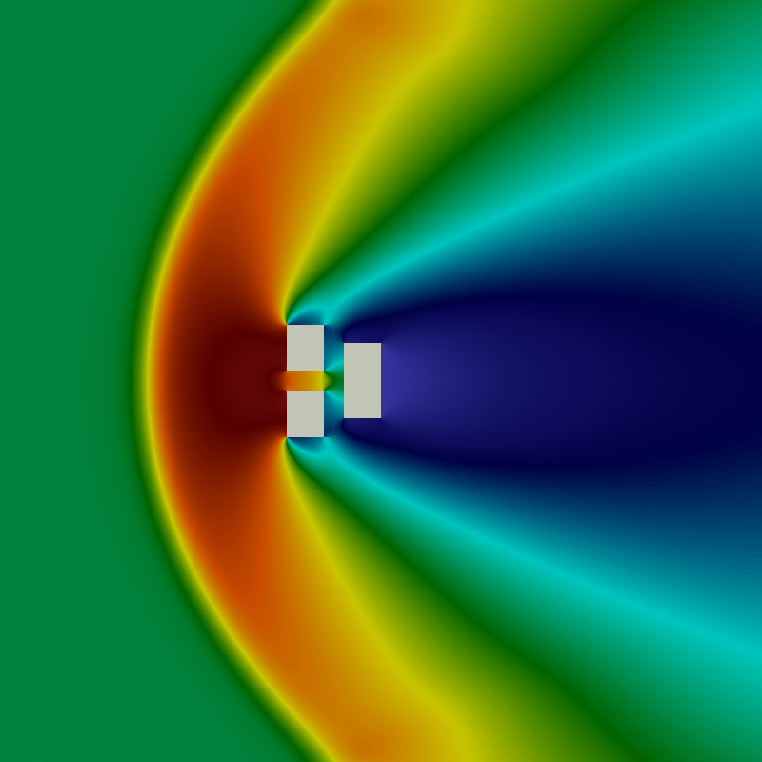} &
      \hspace{-0.5cm}
      \includegraphics[width=0.32\textwidth]{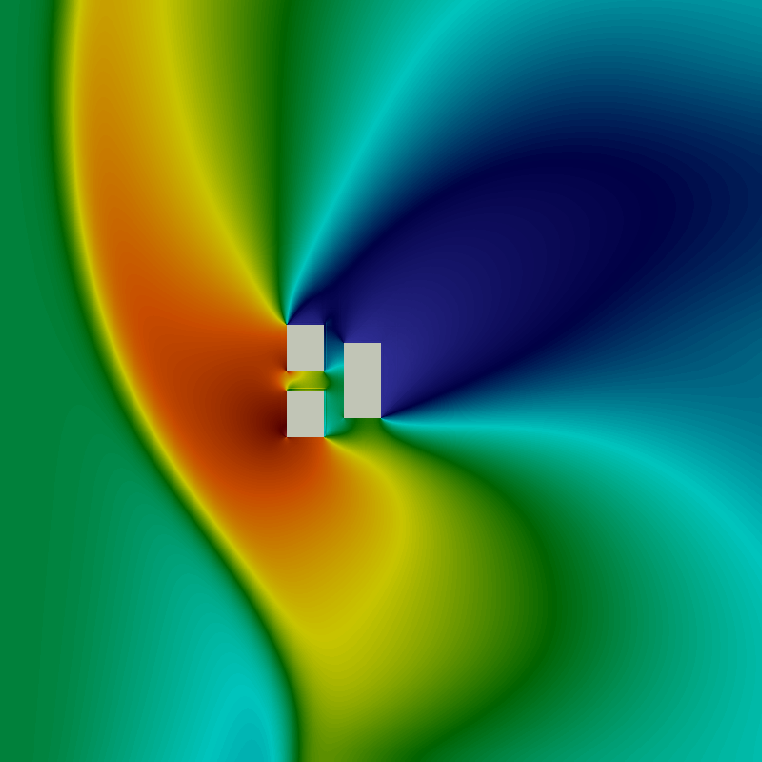} &
      \hspace{-0.5cm}
      \includegraphics[width=0.32\textwidth]{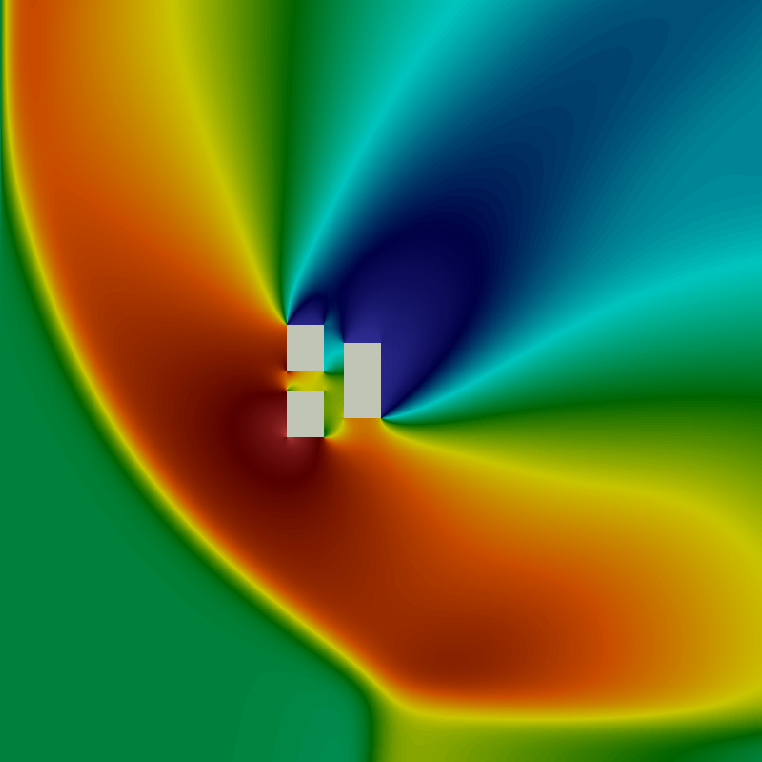} \\
      \hspace{-0.5cm}
      \includegraphics[width=0.32\textwidth]{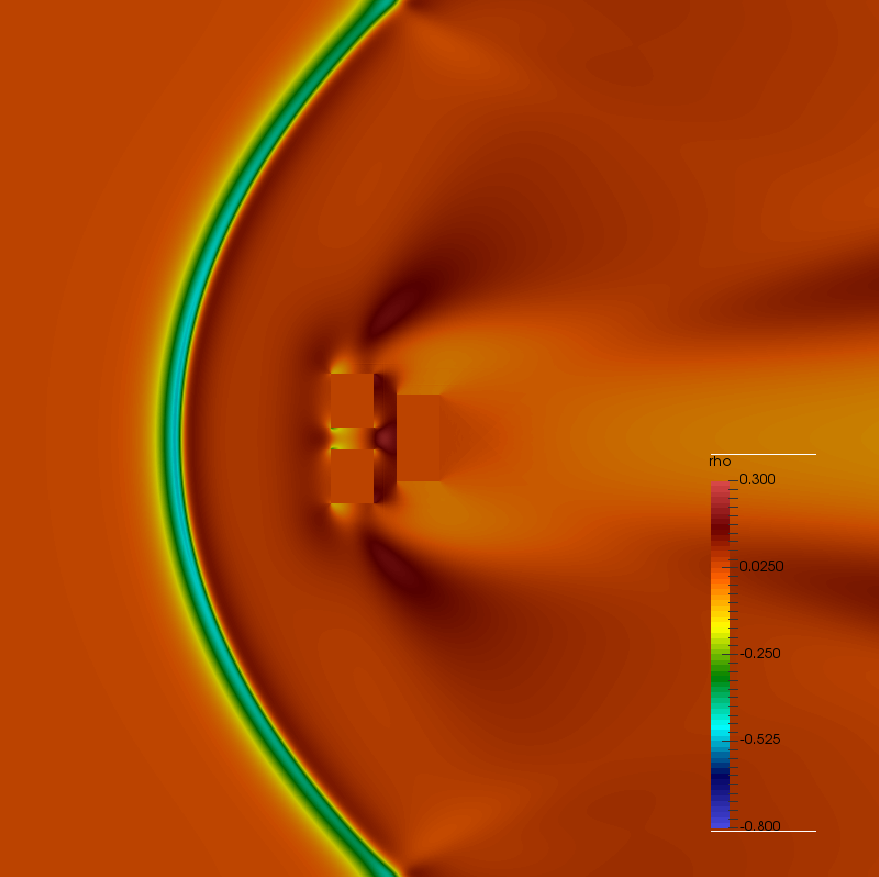} &
      \hspace{-0.5cm}
      \includegraphics[width=0.32\textwidth]{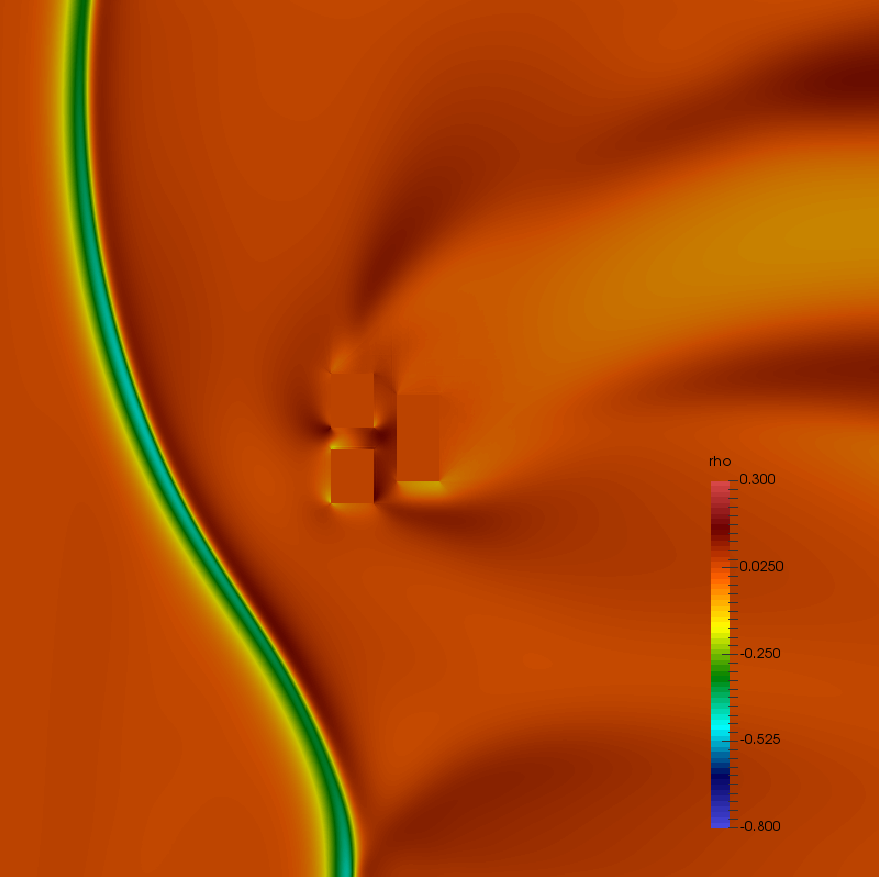} &
      \hspace{-0.5cm}
      \includegraphics[width=0.32\textwidth]{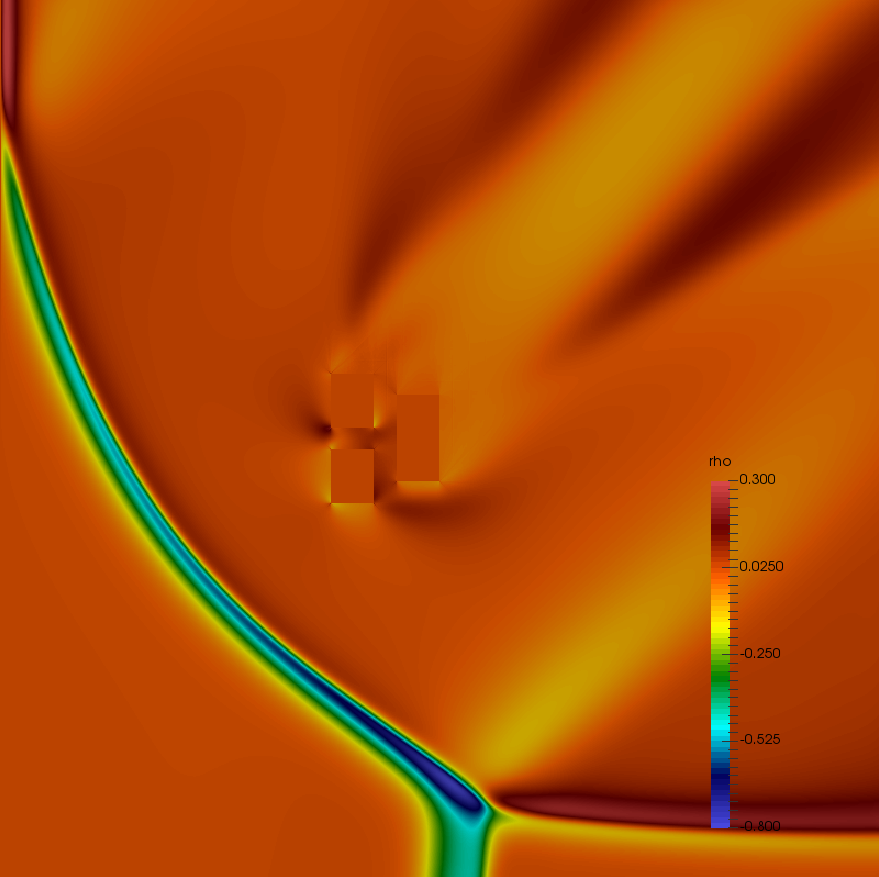} \\
    \end{tabular} 
    \caption{
    \label{fig:test7_reentry2Dx2D_BGK0vsBOZ0}  
    Test 3.2. Two dimensional re-entry test case for $\tau=10^{-2}$ 
      with $M=800 \times 800$ spatial cells and $N=32^2$ velocity cells. BGK model (top row)
      Boltzmann model (middle row), difference between the two models (bottom row) at iterations $4500$,  $12000$ and $20000$.
      The color legend for the first two rows can be found in figure~\ref{fig:test7_reentry2Dx2D_BGK0}. Density profile.
    }
  \end{center}
\end{figure}

\begin{figure}
  \begin{center}
    \begin{tabular}{ccc} 
      \hspace{-0.5cm}
      \includegraphics[width=0.32\textwidth]{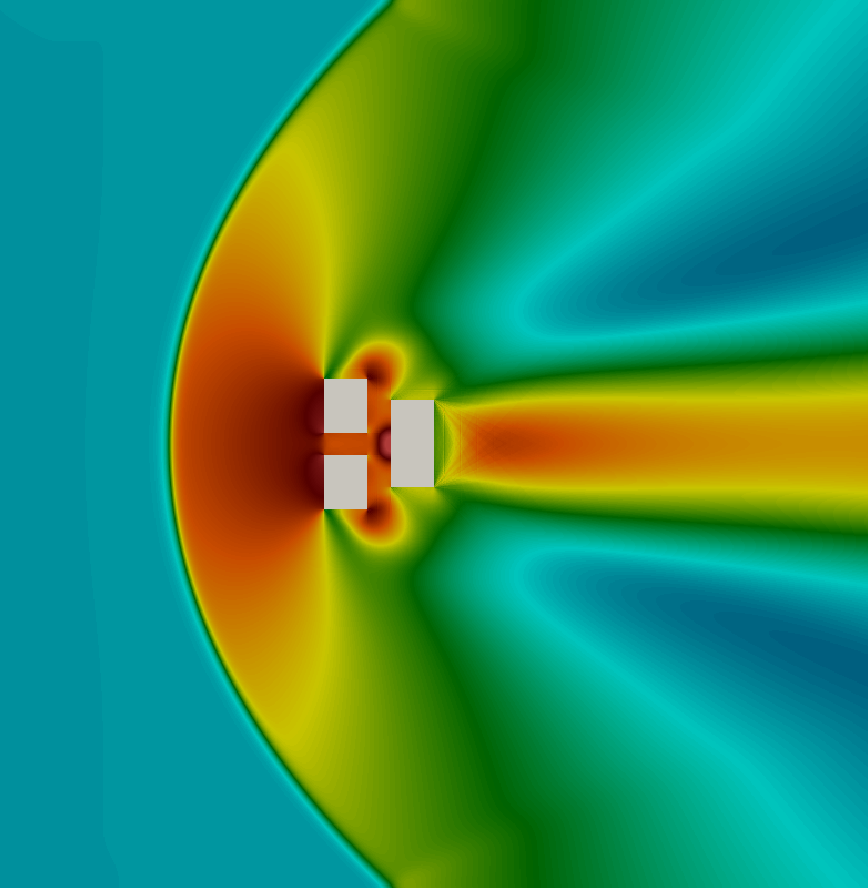} &
      \hspace{-0.5cm}
      \includegraphics[width=0.32\textwidth]{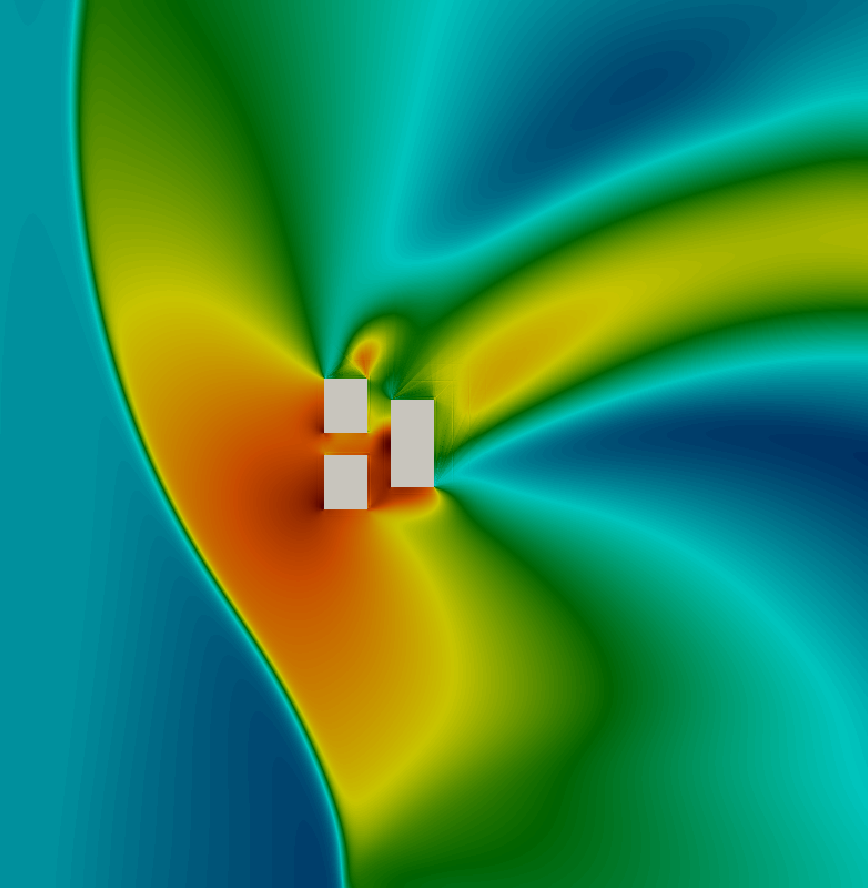} &
      \hspace{-0.5cm}
      \includegraphics[width=0.32\textwidth]{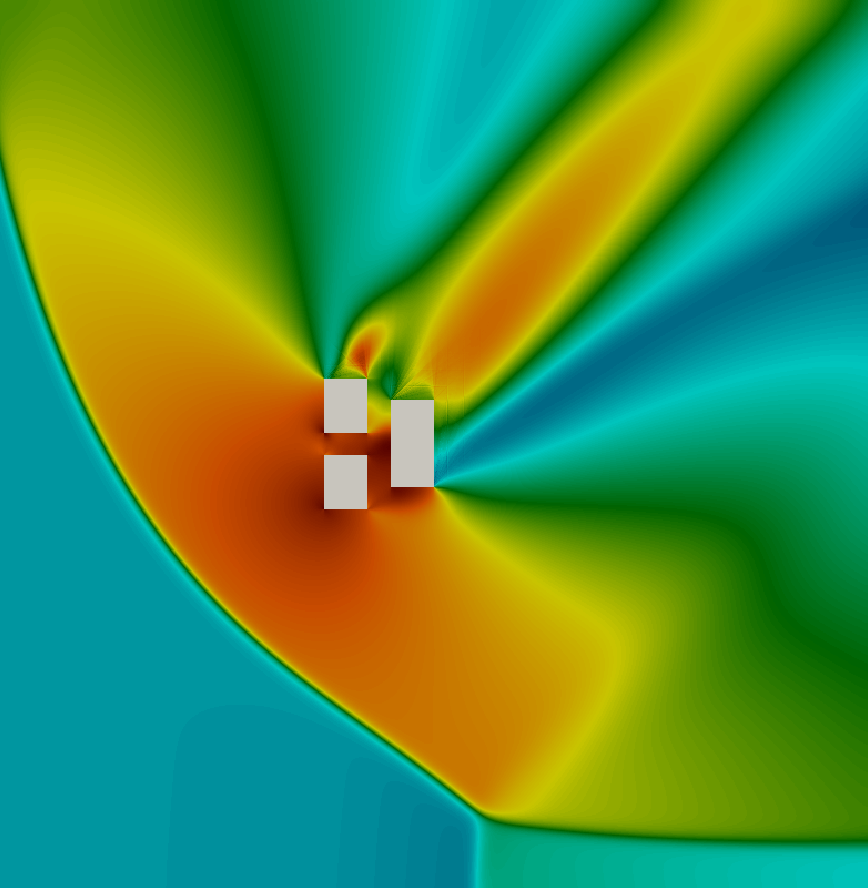} \\
      \hspace{-0.5cm}
      \includegraphics[width=0.32\textwidth]{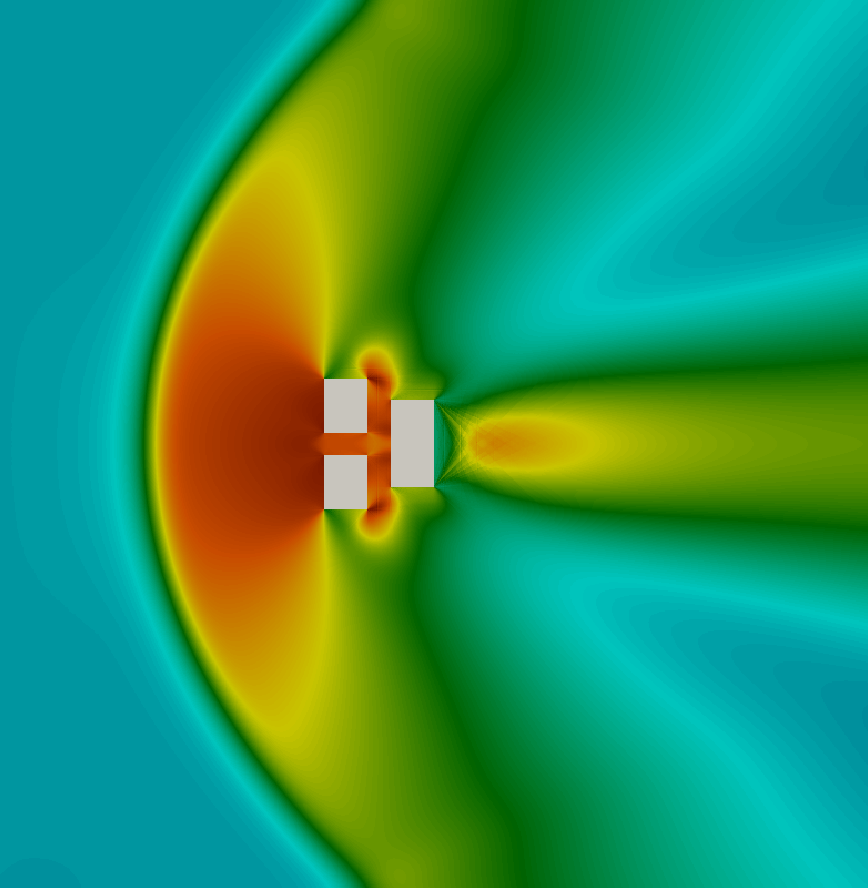} &
      \hspace{-0.5cm}
      \includegraphics[width=0.32\textwidth]{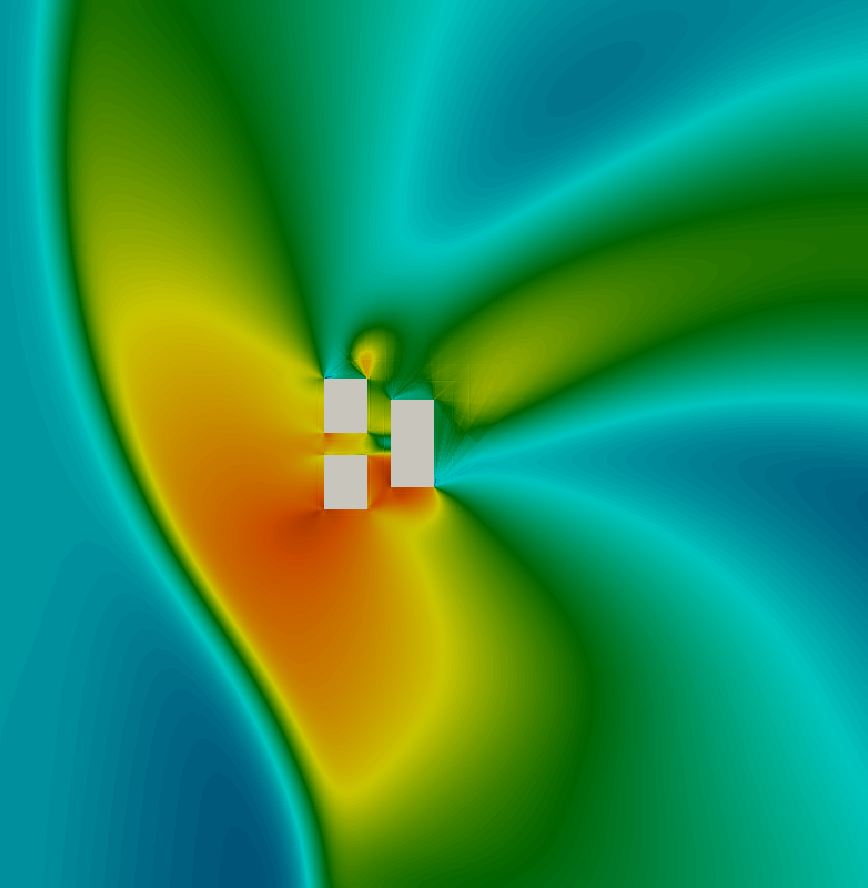} &
      \hspace{-0.5cm}
      \includegraphics[width=0.32\textwidth]{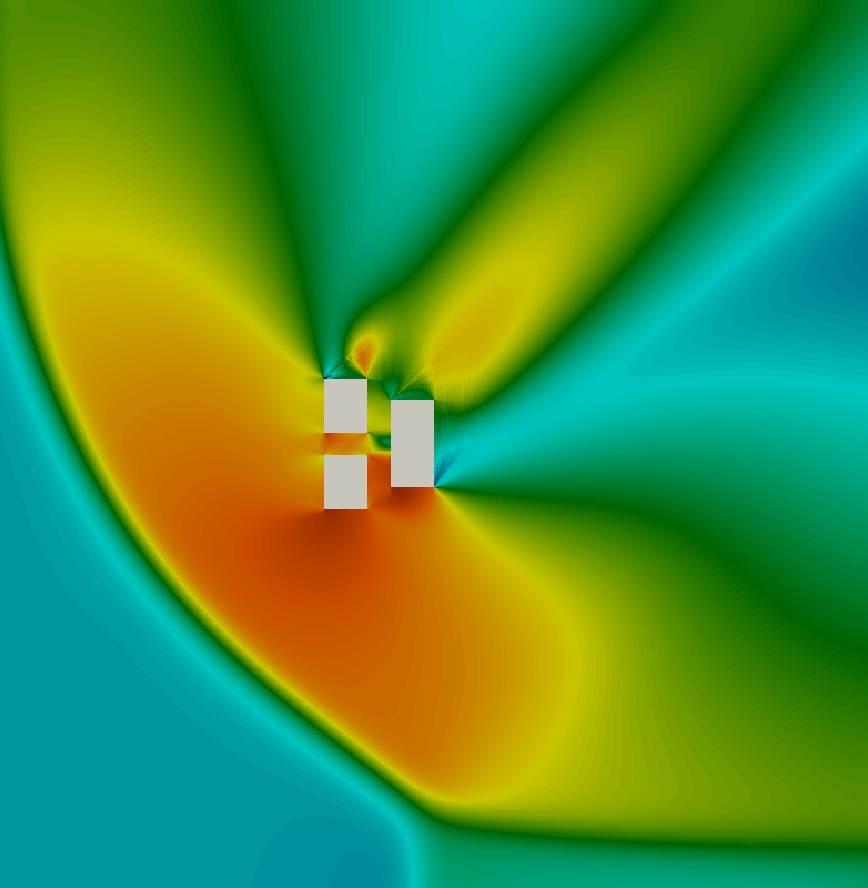} \\
      \hspace{-0.5cm}
      \includegraphics[width=0.32\textwidth]{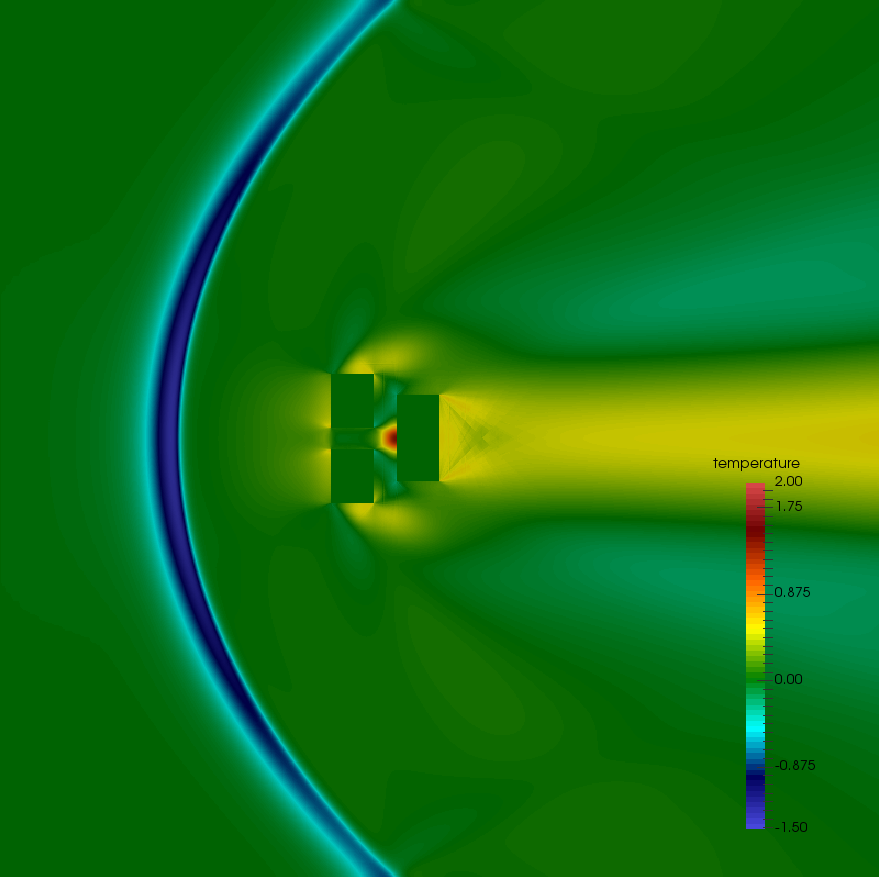} &
      \hspace{-0.5cm}
      \includegraphics[width=0.32\textwidth]{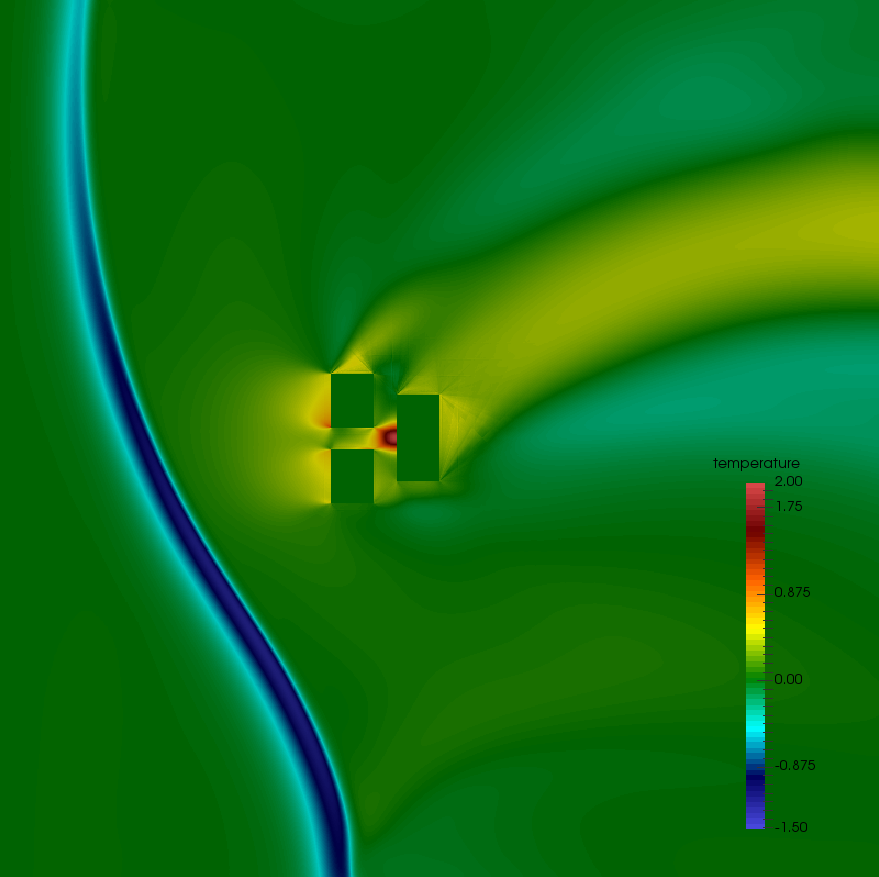} &
      \hspace{-0.5cm}
      \includegraphics[width=0.32\textwidth]{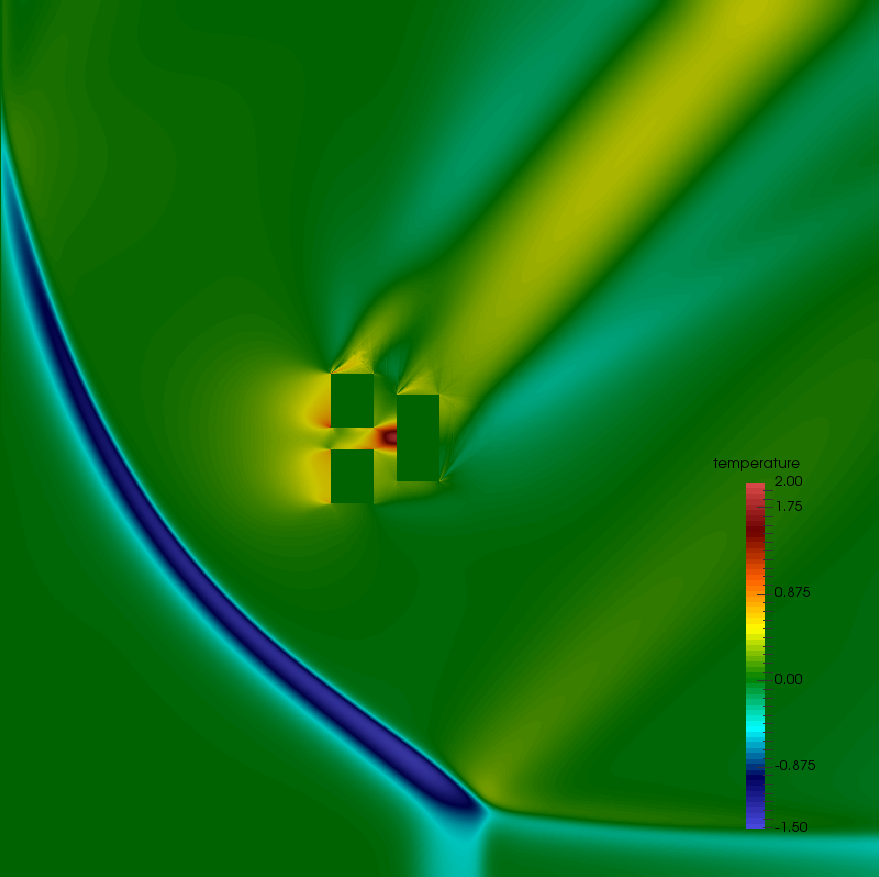} \\
    \end{tabular} 
    \caption{
    \label{fig:test7_reentry2Dx2D_BGK0vsBOZ0_Temp}  
    Test 3.2. Two dimensional re-entry test case for $\tau=10^{-2}$ 
      with $M=800 \times 800$ spatial cells and $N=32^2$ velocity cells. BGK model (top row)
      Boltzmann model (middle row), difference between the two models (bottom row) at iterations $4500$,  $12000$ and $20000$.
      The color legend for the first two rows can be found in figure~\ref{fig:test7_reentry2Dx2D_BGK0}. Temperature profile.
    }
  \end{center}
\end{figure}

In Figure~\ref{fig:test7_reentry2Dx2D_BGK0} we present the results for nine intermediate times
when $\tau=0$ for the density. These results correspond to the fluid limit model that could be obtained
when solving compressible Euler equations.
Next, in Figures~\ref{fig:test7_reentry2Dx2D_BGK0vsBOZ0} and \ref{fig:test7_reentry2Dx2D_BGK0vsBOZ0_Temp}, we compare the BGK results (top row)
with the Boltzmann ones (middle row) for intermediate iterations $4500$, $12000$ and $20000$ for $\tau=10^{-2}$ for respectively the density and the temperature.
In the bottom row, it is shown the differences between the two models in terms of density and temperature. The color legend is the same for the first
two rows for both Figures and it can be found in Figure~\ref{fig:test7_reentry2Dx2D_BGK0}, while for the bottom rows they are shown in the pictures.


In order to conclude this part, let us present some performance data related to those simulations.
The total amount of CPU time needed to compute the $N_{\text{cycle}}=26000$ cycles for the 
BGK solution is $14.5$h, while for Boltzmann model is $304$h$=12.67$d.
The ratio is of the order $21$ in favor of BGK consistently with the previous simulations. 
However, even if Boltzmann results demand a large amount of CPU time, we have seen that discrepancies do
exist with respect to BGK model and in some cases, especially far from equilibrium, they cannot be ignored. 

%
%
\subsection{Part 4. Numerical results for the space three dimensional case.} \label{sec:act3}
In this last part, we present one numerical test in which we compare the relaxation model with the Boltzmann model in three space and velocity dimensions in a kinetic regime. Solving the full Boltzmann equation in three dimensions is extremely resource consuming, even if the fast spectral methods is used and consequently shared memory systems are not sufficient for this kind of problems.
For this reason, only for this last case we adapted the method to distributed memory systems by employing MPI architecture as already stated and described in Section \ref{sec:implementation}.
The results reported for this situation are not to be intended as optimal since we are adapting the scheme to this kind of architecture and improvements in terms of efficiency are attended in the next future.

\subsubsection{Test 4.1. Three dimensional re-entry test case.}
\begin{figure}
  \begin{center}
\begin{tabular}{cc}
          \includegraphics[width=0.37\textwidth]{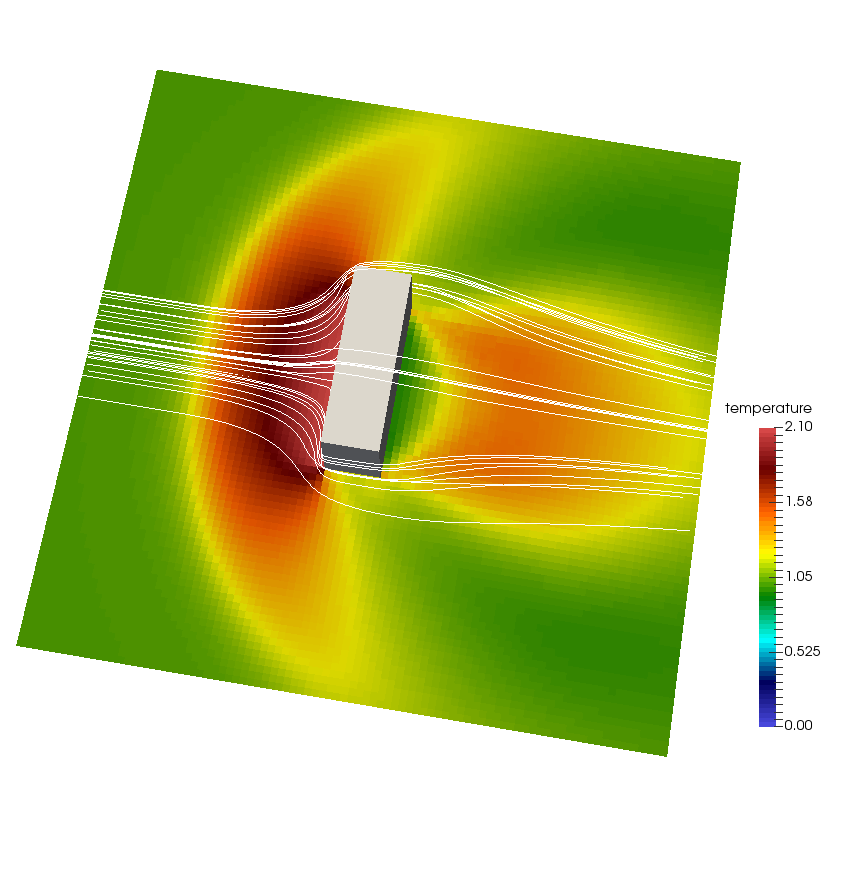}       &       \includegraphics[width=0.37\textwidth]{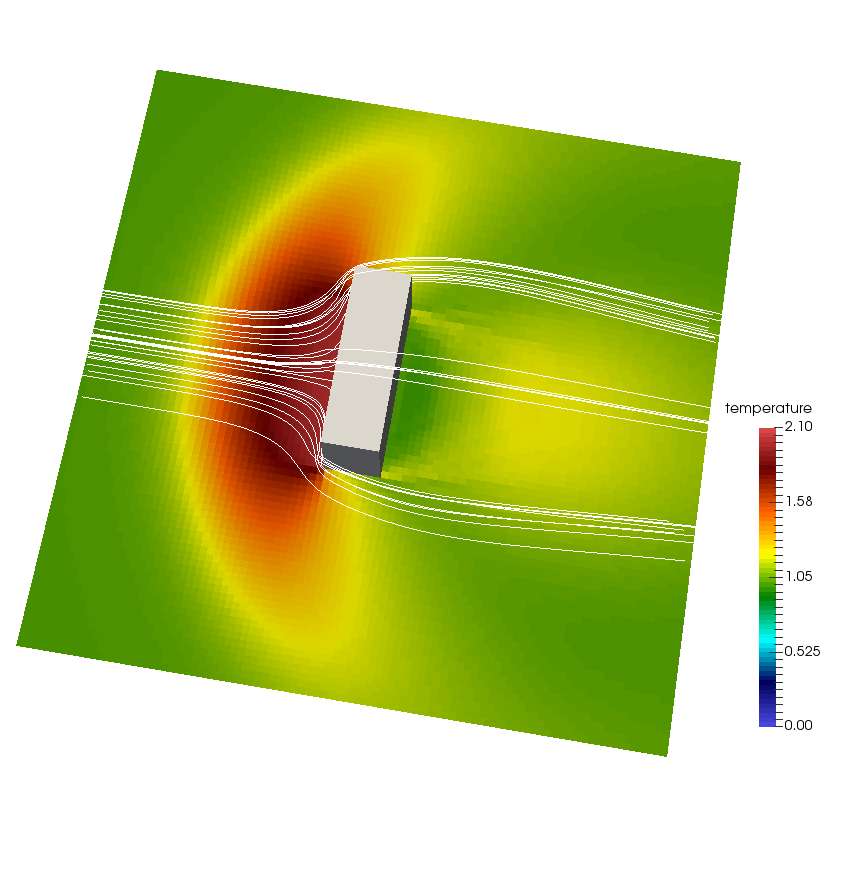}     \\
\multicolumn{2}{c}{  \includegraphics[width=0.37\textwidth]{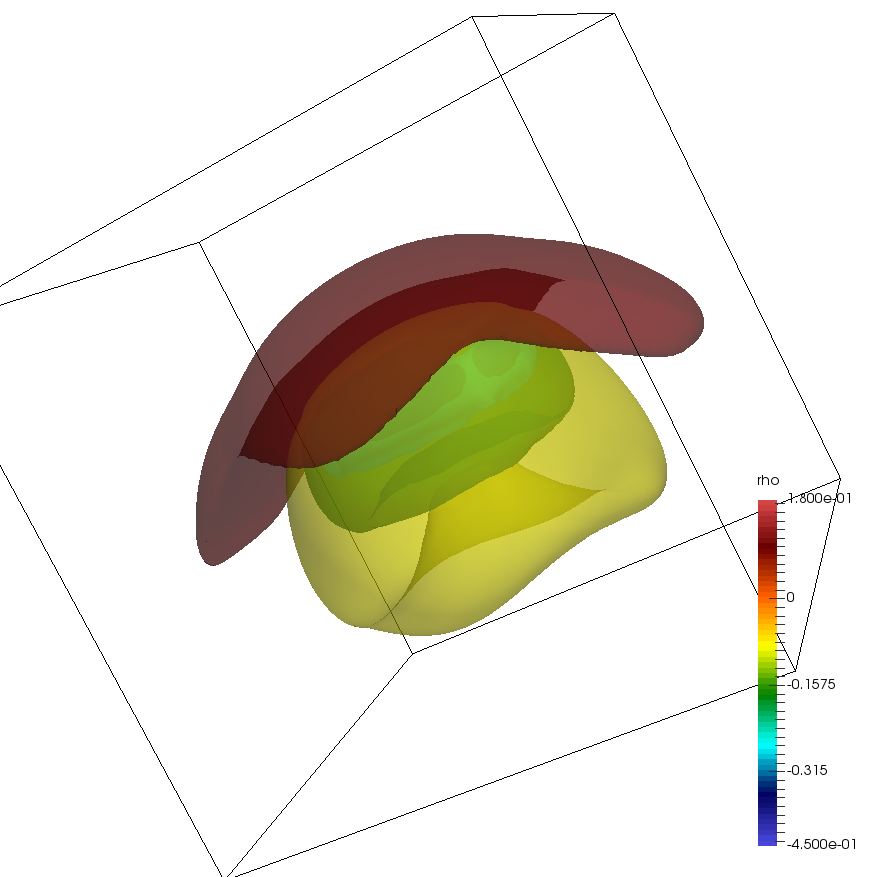} } \\
\multicolumn{2}{c}{BGK-Boltzmann} \\
     \includegraphics[width=0.37\textwidth]{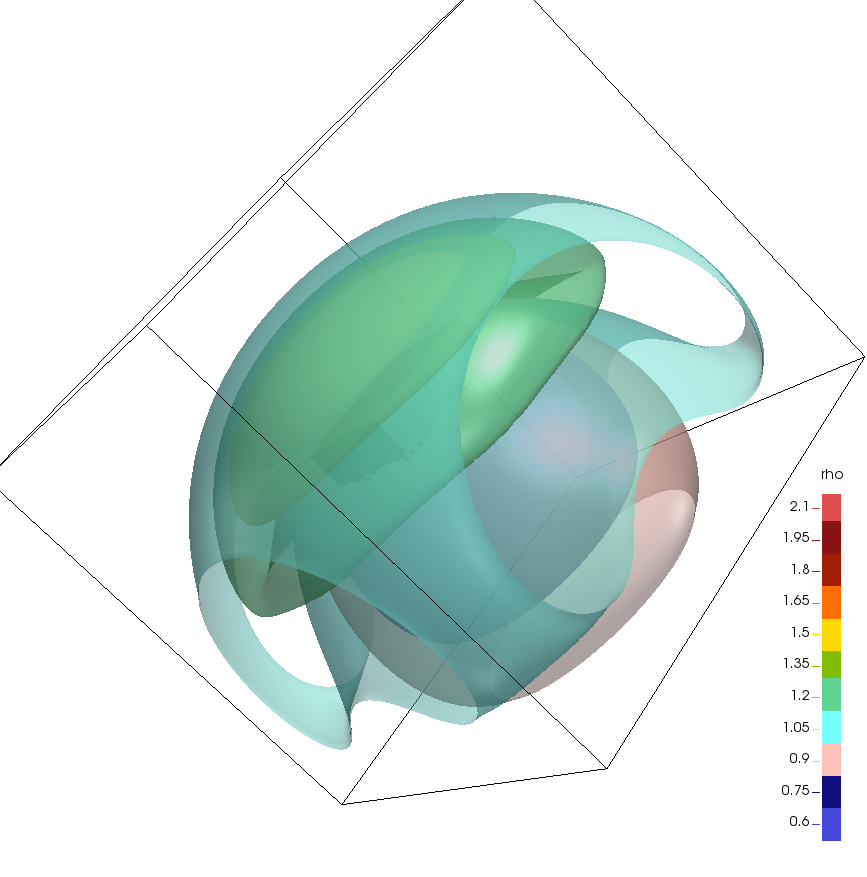}         &        \includegraphics[width=0.37\textwidth]{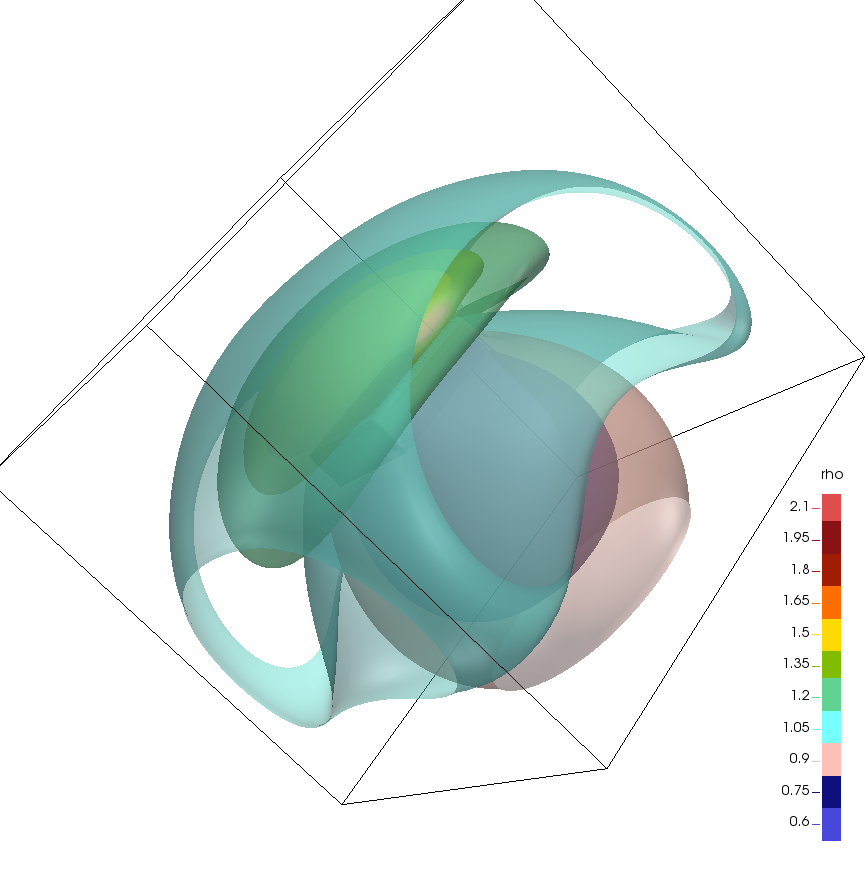}      \\
      BGK &  Boltzmann   
\end{tabular}
  \end{center}
  \caption{
    \label{fig:reentry3D}  
    Test 4.1. Three dimensional re-entry test case for $\tau=0.3$ 
    with $M=90 \times 90 \times 90$ spatial cells and $N=32^3$ velocity cells.  BGK model (left column), Boltzmann model (right column) at time $t_{\text{final}}=0.6$. Top row:
    temperature field with velocity streamlines, bottom row: isosurfaces of the density. Middle row: isosurfaces of the density difference between BGK and Boltzmann models.
  }
\end{figure}

\begin{figure}
  \begin{center}
    \begin{tabular}{cc}
      \includegraphics[width=0.45\textwidth]{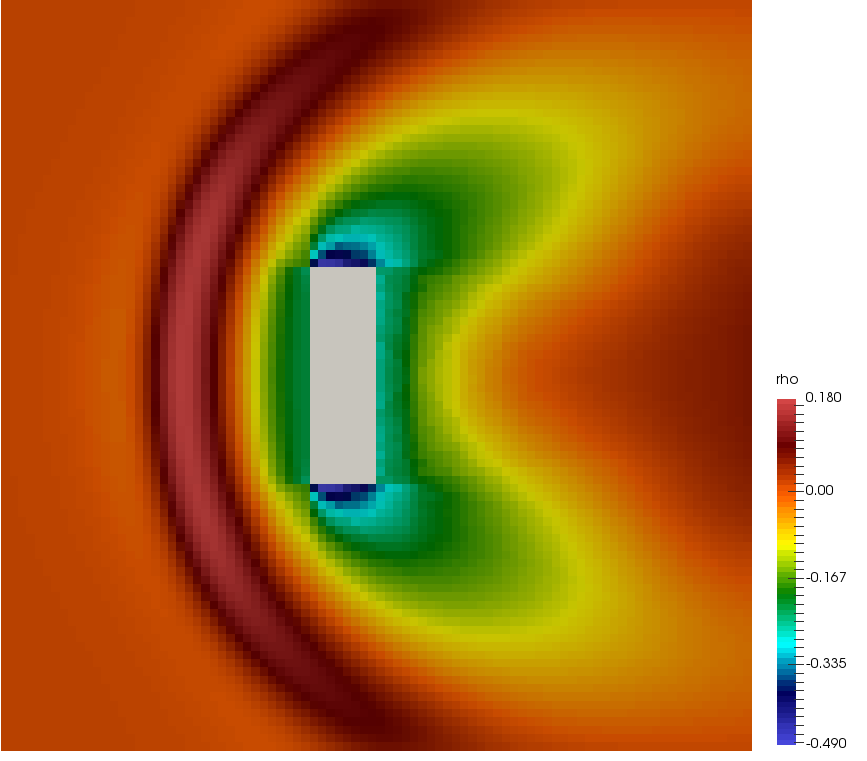} &
      \includegraphics[width=0.45\textwidth]{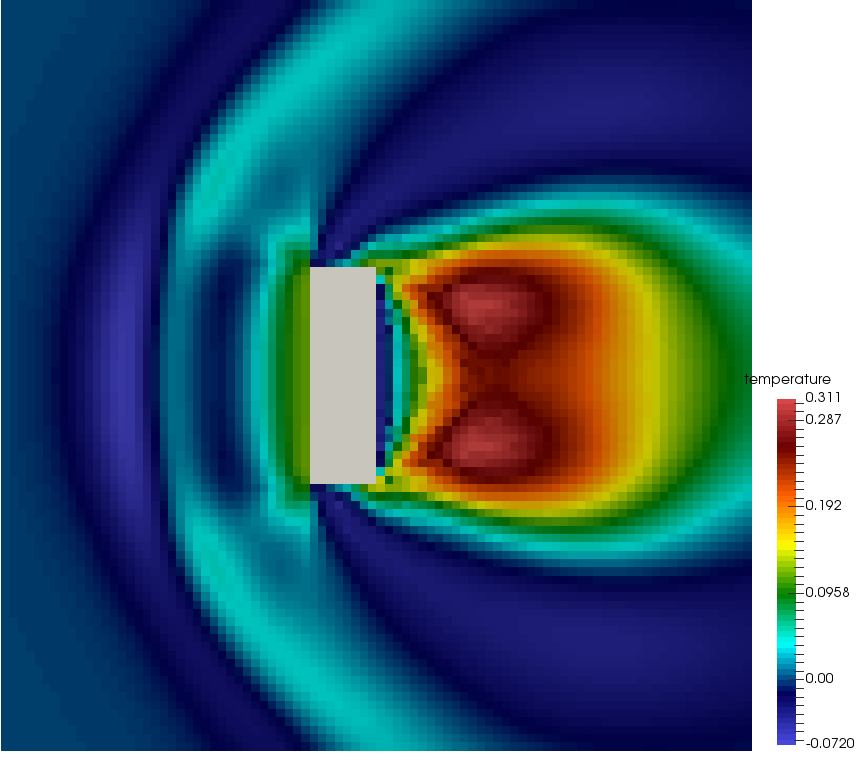}
      \\
      density & temperature \\
    \end{tabular}
  \end{center}
  \caption{
    \label{fig:reentry3D_diff}  
    Test 4.1. Three dimensional re-entry test case for $\tau=0.3$ 
    with $M=90 \times 90 \times 90$ spatial cells and $N=32^3$
    velocity cells. Discrepancies between the BGK solution and the
    Boltzmann solution at time $t_{\text{final}}=0.6$. Density on the left and
    temperature on the right.
  }
\end{figure}

The computational domain is set to $\Omega = [0,2]^3$ with a static cuboid placed in the center (see Fig.\ref{fig:reentry3D}). The velocity space is $[-10,10]^3$ and discretized with $32^3$ points. The
relaxation parameter is set to $\tau =0.3$. The initial density $\rho$ is set to $1$, the temperature $T=1$ while the initial velocity is given by $(u_x,u_y,u_z) = (2,0,0)$. The final time is set to $t_{\text{final}}= 0.6$ leading to $379$ time steps. The inflow boundary conditions are imposed on the left boundary ($x=0$) while outflow boundary conditions on the remaining part of the boundary are imposed. 
Hard sphere molecules are considered for Boltzmann while for the BGK model the frequency $\nu$ is taken equal to $\mu=C_\alpha 4\pi(2\lambda \pi)^\alpha$. For both models the CFL condition considered is consequently given by \be \Delta t\leq \min \left(\frac{\Delta x}{|v_{max}|},\frac{\tau}{\mu}\right).\ee
The results are shown for the temperature and the density in Figure \ref{fig:reentry3D} while the discrepancies between the two in Figure \ref{fig:reentry3D_diff}. 
From the analysis of such results it clearly emerge a difference in the profiles of the macroscopic quantities between the
two models.
We now analyze the performances. This test case was run on the EOS supercomputer at CALMIP, Toulouse France (\url{https://www.calmip.univ-toulouse.fr/}). The supercomputer is equipped with $612$
computational nodes, each of them containing two Intel$^\text{\textregistered}$ Ivybridge $2.8$GHz 10 core CPUs and 64 GB of RAM. Each CPU is equipped with $25$MB of cache memory. The code
was compiled with gcc-5.3.0 and executed on 90 computational nodes. That is to say, on 1800 computational cores in parallel. In the case of the Boltzmann collision operator the Fast Fourier Transforms
were computed by means of the fftw library, version 3.3.4.  The total runtime was equal to $t=93713s$ ($26h$) for the Boltzmann model. This is equivalent to $46000$ computational hours ($3.25$ years) on
a sequential machine. On the other hand, the runtime for the BGK model was only $t = 2174s$ ($0.6h$) with a ratio $43$ in favor of the simpler relaxation model. The profiling data are
summarized in Table \ref{tab:3Dprof}. The communications and MPI synchronization take $38\%$ of the computational time for BGK equation. 
On the other hand, for the Boltzmann model the time spent on communications and synchronization amounts merely to $3.4\%$ of the
total runtime. This is due to extreme computational complexity of the three dimensional Boltzmann collision kernel. As discussed in the previous paragraphs, the costs relatives to the routines \texttt{ ToConservative}, \texttt{ ToPrimitive} which are relative to the computation of the macroscopic variables from the distribution function can be imputed to the transport phase, since this is the only
part of the scheme which modifies these quantities, since collisions preserves density, momentum and energy in the cell.

\begin{table}
  \begin{center}
  \numerikNine
  \begin{tabular}{|c|c|c||c||cc|}
    \hline
    &  \textbf{Cycle} &  \textbf{CPU}  &  \textbf{Main routines} 
    &  \textbf{Cost CPU }  &  \textbf{Cost}  \\
    & & (s) & & vs total (s) & vs total (\%)  \\
    \hline
    \hline
     \multirow{6}{*}{\begin{sideways} \textbf{BGK}  \end{sideways}}
    & \multirow{6}{*}{ 379 }
    & \multirow{6}{*}{ 2174 }
    &Transport          &   0.03 & $<$0.1\%  \\
    \cline{4-6}
    & & &ToConservative &  1126 & 52\%  \\
    \cline{4-6}
    & & &ToPrimitive    &  0.06  & $<$0.1\%  \\
    \cline{4-6}
    & & &Collisions     &  217 & 10\% \\  
    \cline{4-6}
    & & & MPI comm. & 825 & 38\%  \\
    \cline{4-6}
    & & & = & 2168 & 100\% \\
    \hline
    \hline
    \multirow{6}{*}{\begin{sideways} \textbf{Boltzmann}  \end{sideways}}
    & \multirow{6}{*}{ 379 }
    & \multirow{6}{*}{ 93713 }
    &Transport          &    0.5 &  $<$0.1\%  \\
    \cline{4-6}
    & & &ToConservative &   1127 &  1.2\%  \\
    \cline{4-6}
    & & &ToPrimitive    &    0.07 &  $<$0.1\%  \\
    \cline{4-6}
    & & &Collisions     & 89396 & 95.4\%  \\  
    \cline{4-6}
    & & & MPI comm.  & 3190 & 3.4\%  \\
    \cline{4-6}
    & & & = & 93713 & 100\% \\
    \hline
  \end{tabular}
  \caption{ \label{tab:3Dprof} Profiling of the average cost
    for each routine of the code on the test 4.1. (three dimensional
    re-entry test case) simulated using $M=90^3$ and $N=32^3$ points.  }
  \end{center}
\end{table}



\section{Conclusion and perspectives} \label{sec:conclusion}
In this paper we have generalized the Fast Kinetic Scheme \cite{FKS, FKS_HO, FKS_GPU, RFKS} to the challenging case of the
Boltzmann collision integral. We have shown that it is possible to solve the full unsteady three dimensional Boltzmann equation
with variable hard sphere kernel in a reasonable amount of time by using parallel architectures. Up to the author knowledge, the results
reported in this paper represent one of the very first attempts of solving the seven dimensional Boltzmann model with deterministic 
numerical schemes. This has been made possible by combining a fast semi-Lagrangian approach for the transport part with a fast spectral method for the collision dynamic.

We have performed several numerical tests with the aim of detailing the behavior of the method in different situations in order to understand its strengths and weaknesses.
A side scope of the paper has been to show the differences that arise between the Boltzmann model and the popular BGK relaxation model. 

Differences have been observed to be large far from equilibrium situations. 

In the future, we aim in working in the direction of additionally improving the fast spectral method since we have observed it to be one of the main bottleneck of the scheme by using for
instance different grids during the transport and the collision phases. Another fundamental direction we aim to pursue is the optimization of the MPI parallelization which is 
necessary for considering more realistic applications as well as the development of techniques for treating complex boundaries. Finally, the extension of the present scheme
to plasmas is under study.


\section*{Acknowledgments}
The authors would like to thanks Professor Francis Filbet from the University of Toulouse III, Professor R{\'e}mi Abgrall and Dott. Tulin Kaman from the institute of mathematics in Z{\"u}rich, Switzerland
for their suggestions and help.\\
This work has been supported by the Galileo project G14 (Fast Asymptotic-Preserving and semi-Lagrangian schemes for High Performance Computing : applications to plasmas)
from the Franco-Italian University and by the ANR project MOONRISE (MOdels, Oscillations and NumeRIcal SchEmes, 2015-2019). 
Thomas Rey was partially funded by Labex CEMPI (ANR-11-LABX-0007-01).\\
This work was granted access to the HPC resources of CALMIP supercomputing center under the allocation 2016-P1542 and the authors acknowledge the help from CALMIP. \\

\section*{References}
\bibliography{biblio.bib}

\end{document}